\documentclass[10pt]{article}
\usepackage[final]{graphicx}
\usepackage{amsfonts}
\usepackage{amscd}
\usepackage{amsmath}
\usepackage{latexsym, euscript}
\begin{document}

\newcommand{\CLIP}[2]{\left(#1,#2 \right)_{L^2}} 
\newcommand{\DLIP}[2]{\left(#1,#2 \right)_{L^2_h}} 
\newcommand{\DSIP}[2]{\left(#1,#2 \right)_{S_h}} 
\newcommand{\CLN}[1]{\left\| #1 \right\|_{L^2}}  
\newcommand{\CLNtx}[1]{\left\| #1 \right\|_{L^2_t L^2_x}}  
\newcommand{\CLNStx}[1]{\left\| #1 \right\|_{L^2_t L^2_x}^2}  
\newcommand{\CLNxt}[1]{\left\| #1 \right\|_{L^2_x L^2_t}}  
\newcommand{\DLN}[1]{\left\| #1 \right\|_{L^2_h}}  
\newcommand{\DLNt}[1]{\left\| #1 \right\|_{L^2_k}}  
\newcommand{\DLNSxh}[1]{\left\| #1 \right\|_{L^2_{x,h}}^2}  
\newcommand{\DLNStk}[1]{\left\| #1 \right\|_{L^2_{t,k}}^2}  
\newcommand{\DLNtx}[1]{\left\| #1 \right\|_{L^2_{t,k} L^2_{x,h}}}  
\newcommand{\DLNxt}[1]{\left\| #1 \right\|_{L^2_{x,h} L^2_{t,k}}}  
\newcommand{\DSN}[1]{\left\| #1 \right\|_{S_h}}  

\newcommand{\CLNS}[1]{\left\| #1 \right\|_{L^2}^2}  
\newcommand{\CLNSx}[1]{\left\| #1 \right\|_{L^2_x}^2}  
\newcommand{\CLNSt}[1]{\left\| #1 \right\|_{L^2_t}^2}  
\newcommand{\DLNS}[1]{\left\| #1 \right\|_{L^2_h}^2}  
\newcommand{\DSNS}[1]{\left\| #1 \right\|_{S_h}^2}  
\newcommand{\SupNorm}[1]{\left\| #1 \right\|_{\infty}}  
\newcommand{\SupNormtx}[1]{\left\| #1 \right\|_{L^{\infty}_t L^{\infty}_x}}  
\newcommand{\SupNormxt}[1]{\left\| #1 \right\|_{L^{\infty}_x L^{\infty}_t}}  
\newcommand{\SupNormtLNx}[1]{\left\| #1 \right\|_{L^{\infty}_t L^2_x}}
\newcommand{\LNxSupNormt}[1]{\left\| #1 \right\|_{L^2_x L^{\infty}_t}}

\newcommand{\IHN}[1]{\left \langle #1 \right \rangle}  

\centerline{\large \bf Unbounded Solutions of the Modified Korteweg-De Vries Equation}\vspace{0.1truein}

\def\thefootnote{\arabic{footnote}}
\begin{center}
  John Gonzalez  \footnote{Email:  johngonz@gmail.com } \\
  \textit{Northeastern University}
\end{center}

\begin{abstract}
\noindent We prove local existence and uniqueness of solutions of the focusing modified Korteweg - de Vries equation $u_t + u^2u_x + u_{xxx} = 0$ in classes of unbounded functions that admit an asymptotic expansion at infinity in decreasing powers of $x$.  We show that an asymptotic solution differs from a genuine solution by a smooth function that is of Schwartz class with respect to $x$ and that solves a generalized version of the focusing mKdV equation.  The latter equation is solved by discretization methods. \\
\end{abstract}
\section{Introduction}  
In this article we consider the focusing modified Korteweg-De Vries equation \\

\begin{eqnarray}
 \left\{ 
\begin{array}{rcl}
w_t + w^2 w_x + w_{xxx} &= &0 \label{dagger} \\
w_{|t=0}& =& w_0(x) 
\end{array} 
\right. 
\end{eqnarray}

\noindent where the initial data $w_0$ is possibly unbounded at $+\infty$ and/or $-\infty$.  We construct local (in time) solutions to $(\ref{dagger})$ that lie in the spaces $S^{\beta}\left(\mathbb{R}\times I \right)$ introduced in [3] by T. Kappeler, P. Perry, M. Shubin, and P. Topalov.  These spaces are defined as follows:  \\
\indent Let $I \subset \mathbb{R}$ be an interval and $\beta \in \mathbb{R}$ be given.  Denote by $S^{\beta}\left(\mathbb{R}\times I \right)$ the linear space of $C^{\infty}\left(\mathbb{R}\times I \to \mathbb{R}\right)$ functions having asymptotic expansions at $\pm\infty$ given by $f(x,t) \sim \sum_{k=0} ^{\infty} a_k ^{+}(t) x^{\beta_k}$ as $x \to \infty$ and $f(x,t) \sim \sum_{k=0} ^{\infty} a_k ^{-}(t) x^{\beta_k}$ as $x \to -\infty$ where $a_k^{\pm}\in C^{\infty}(I \to \mathbb{R})$ and $\beta = \beta_0 > \beta_1 > \cdots$ with $\lim_{k \to \infty} \beta_k = -\infty$.  By definiton, the asymptotic relation $\sim$ means that for every compact interval $J\subset I$, and integers $N,i,j \geq 0$ there exists $C_{J,N,i,j} > 0$ such that for any $ \pm x \geq 1$ and $t\in J$ we have.  \\
$$\left| \partial^i_t \partial^j_x\left( f(x,t) - \sum_{k=0} ^N a_k ^{\pm}(t)(\pm x)^{\beta_k}\right) \right| \leq C_{J,N,i,j} \left| x \right| ^{\beta_{N+1} - j}$$
\noindent  We denote by $S^{-\infty}\left(\mathbb{R}\times I \right)$ the space of $C^{\infty}\left(\mathbb{R}\times I \to \mathbb{R}\right)$ functions having asymptotic expansions at $\pm \infty$ which are identically zero.  Analogously, we define the spaces $S^{\beta}\left(\mathbb{R} \right)$ and $S^{-\infty}\left(\mathbb{R}\right)$ as the space of functions $f(x) \in C^{\infty}(\mathbb{R} \to \mathbb{R})$ having such asymptotic expansions where the coefficients $a_k^{\pm}$ are constants independent of $t$.  We shall construct solutions $w(x,t) \in S^{\beta}(\mathbb{R}\times I)$ for $(\ref{dagger})$ with initial data $w_0\in S^{\beta}(\mathbb{R})$ when $\beta \leq \frac{1}{2}$. \\
\indent  If $w(x,t) \in S^{\beta}(\mathbb{R}\times I)$ is a solution for $(\ref{dagger})$ then one expects its asymptotic expansions $\sum_{k=0} ^{\infty} a_k ^{\pm}(t) (\pm x)^{\beta_k}$, although not generally convergent, to give formal solutions (see lemma 5.2).  We define a pair of formal power series $\sum_{k=0} ^{\infty} a_k ^{\pm}(t) x^{\beta_k}$ to be a formal solution to $(\ref{dagger})$ if $\sum_{k=0} ^{\infty} a_k ^{+}(t) x^{\beta_k}$ and $\sum_{k=0} ^{\infty} a_k ^{-}(t) (-x)^{\beta_k}$ satisfy $(\ref{dagger})$ for all $t\in I$ when $x$ is taken as a formal variable and differentiation in $x$ is carried out in the ordinary way.  \\
\indent When $\beta > \frac{1}{2}$ one can easily see that there are generally no formal solutions to $(\ref{dagger})$ and hence no solutions in $S^{\beta}(\mathbb{R}\times I)$.  Indeed, if $\sum_{k=0} ^{\infty} a_k ^{+}(t) x^{\beta_k}$ satisfies $(\ref{dagger})$ formally where $\beta>\frac{1}{2}$ and $a_0^+ \ne 0$ then 
\begin{equation*}
\sum_{j=0}^{\infty}\dot{a}_j^+(t) x^{\beta_j} =  - \Big(\sum_{j=0}^{\infty} a_j^+(t) x^{\beta_j}\Big)^2 \cdot\Big(\sum_{j=0}^{\infty}a_j^+(t)\cdot \beta_j \cdot x^{\beta_j-1} \Big) - \Big(\sum_{j=0}^{\infty}a_j^+(t)\cdot \beta_j \cdot( \beta_j -1)\cdot (\beta_j -2) x^{\beta_j-3}\Big)
\end{equation*}
\noindent The largest exponent on the left side is $\beta_0$ and the largest exponent on the right side is $3\beta_0-1$ which is larger than $\beta_0$.  Therefore by equating the coefficients of $x^{3\beta_0 - 1}$ one deduces that $0 = -\beta_0(a_0^+)^3$ which implies that $a_0^+ = 0$, a contradiction.  On the other hand when $\beta \leq \frac{1}{2}$ one has formal solutions defined for $t\in I = \left[\left.-c,\infty \right. \right]$ for some $c>0$ (see lemma 6.2) and therefore one can hope to find solutions in $S^{\beta}(\mathbb{R}\times I)$.  \\
\indent  For an arbitrarily chosen pair of such formal power series $\sum_{k=0} ^{\infty} a_k ^{\pm}(t) (\pm x)^{\beta_k}$ there exists a function $f(x,t)\in C^{\infty}(\mathbb{R}\times I \to \mathbb{R})$ asymptotic to the pair (see for example [7] proposition 3.5).  The function $f$ is not unique but if $\sum_{k=0} ^{\infty} a_k ^{\pm}(t) (\pm x)^{\beta_k}$ is a formal solution then any such $f$ will be an asymptotic solution for $(\ref{dagger})$ (see lemma 6.2).  By definition an asymptotic solution is a function $f\in S^{\beta}(\mathbb{R}\times I)$ such that \\
\[
\left\{ 
\begin{array}{rcl}
f_t + f^2 f_x + f_{xxx} & \in & S^{-\infty} \left(\mathbb{R} \times I\right) \\
f_{|t=0} - w_0 & \in & S^{-\infty}\left(\mathbb{R}\right) 
\end{array} 
\right. 
\]  
\indent  Given an asymptotic solution $f(x,t) \in S^{\beta}(\mathbb{R}\times I)$ for $(\ref{dagger})$ one can attempt to construct a genuine solution $w(x,t) \in S^{\beta}(\mathbb{R}\times I)$ to $(\ref{dagger})$  by constructing $u(x,t) \in S^{-\infty} \left(\mathbb{R} \times I\right)$ such that $w := f+u$ is a genuine solution of $(\ref{dagger})$.  If $w$ satisfies $(\ref{dagger})$ then $u$ must satisfy\\

\begin{eqnarray}
 \left\{ 
\begin{array}{rcl}
u_t + u^2 u_x + u_{xxx}  +  (u^2 f)_x + (f^2 u)_x + g &= &0 \label{daggerdagger}\\
u_{|t=0} &=& u_0(x) 
\end{array} 
\right. 
\end{eqnarray}

\noindent where $u_0 = w_0-f(x,0) \in S^{-\infty}\left(\mathbb{R}\right)$ and $g \in S^{-\infty}\left(\mathbb{R} \times I \right)$ is the result of plugging $f$ into $(\ref{dagger})$. \\
\indent We shall prove existence of finite time solutions $u(x,t) \in S^{-\infty}\left(\mathbb{R} \times \left[ 0,T \right] \right) $ to $(\ref{daggerdagger})$ by using the discretization method introduced by Menikoff in [5] and further developed by Bondareva in [1].  Moreover, uniqueness will also be proven so that we shall show the following theorem: \\

\noindent \textbf{Theorem 1.1} \hspace{1mm} Let $f\in C^{\infty}\left(\mathbb{R}\times \left[\left. 0,\infty \right. \right] \to \mathbb{R}\right)$ be a function satisfying the property that for every compact $J \subset \left[\left. 0,\infty \right. \right]$ we have $f^{(n)}(x,t) = O\left(\left| x \right|^{\frac{1}{2}-n} \right)$ uniformly for $t\in J$ and let $g$ be any function lying in $S^{-\infty}\left(\mathbb{R}\times\left[\left. 0,\infty \right. \right] \right)$.  Suppose $u_0 \in S^{-\infty}\left(\mathbb{R}\right)$.  Then there exists $T>0$ such that  $(\ref{daggerdagger})$ has a solution $u(x,t) \in S^{-\infty}(\mathbb{R}\times \left[0,T\right])$.  Moreover, the solution $u$ is unique in $S^{-\infty}(\mathbb{R}\times \left[0,T\right])$.  \\

\indent The finite-time existence and uniqueness theorem for $(\ref{daggerdagger})$ will enable us to prove finite-time existence and uniqueness for $(\ref{dagger})$ in the space $S^{\beta}(\mathbb{R}\times \left[0,T\right])$ for $\beta \leq \frac{1}{2}$ which can be stated as the following main theorem:  \\

\noindent \textbf{Theorem 1.2} \hspace{1mm}  For any $\beta \leq \frac{1}{2}$ and for any initial condition $w_0 \in S^{\beta}(\mathbb{R})$ there exists a $T>0$ and a unique solution $w(x,t) \in S^{\beta}(\mathbb{R}\times \left[0,T\right])$ of the initial value problem $(\ref{dagger})$.  Moreover, if $w_0 \sim \sum_{k=0}^{\infty} a_k^{\pm} x^{\beta_k}$ and $j$ is the smallest index such that $a_j^{+} \ne 0$ (resp. $a_j^{-} \ne 0$) then the coefficient $a_j^{+}(t)$ (resp. $a_j^{-}(t)$) in the asymptotic expansion of the solution is a nonvanishing continuous function of $t$ and all preceeding coefficients are identically zero.\\

\indent The second statement in theorem 1.2 indicates that the asymptotic growth of the solution is determined throughout its time of existence by the leading exponents in the asymptotic expansion of its initial data.  In particular if $0<\beta_0 \leq \frac{1}{2}$ and $a_0^{\pm}\ne 0$ then the solution $w(x,t)$ for $(\ref{dagger})$ is unbounded in $x$.  \\

\noindent  \textit{Related Work} \hspace{1mm} For the defocusing mKdV equation T. Kappeler, P. Perry, M. Shubin, and P. Topalov constructed global solutions lying in $S^{\beta}(\mathbb{R}\times \mathbb{R})$ (as well as other spaces) for $\beta < \frac{1}{2}$ in [3].  We remark that the methods we use in this article to prove theorems 1.1 and 1.2 can also be applied to the defocusing mKdV equation with no significant changes and it would yield the same results as theorems 1.1 and 1.2 above.  Related results on unbounded solutions for the KdV equation were obtained by I.N. Bondareva and M. Shubin in [1,2] and by A. Menikoff in [5] where the authors construct global in time solutions in certain classes of functions whose spatial growth is of order $|x|^{\beta}$ for $\beta \leq 1$.  In [4] Kenig, Ponce, and Vega constructed unbounded solutions for the KdV equation lying in certain spaces where the initial growth is polynomial but the solutions instantaneously (i.e. for $t>0$) exhibit linear growth in $x$.    \\
\indent In section two we introduce the discretization method of Menikoff and give some general lemmas.  The goal of section two is to prove finite time existence for a discretized version of $(\ref{daggerdagger})$.  Section three contains various estimates which are necessary in order to pass from discrete solutions to smooth solutions.  In section four we show how to pass from discrete solutions to smooth solutions by using a smoothing operator $I_h$, introduced by Stummel in [8].  The existence statements of theorems 1 and 2 are proved in section four and the uniqueness results are proved in section five. 
\section{Discretization of the Generalized mKdV Equation}
\subsection{Definitions and General Setup}
For now let us fix two mesh size numbers $0 < h,k <1 $ and let us denote $x_n := nh$ and $t_j := jk$ for each $n,j \in \mathbb{Z}$.   We shall let $\mathbb{R}_h$ and $\mathbb{R}_k$ denote the (discrete) collection of real numbers of the form $x_n$ and $t_j$ respectively and we shall refer to those sets and the cartesian product $\mathbb{R}_h \times \mathbb{R}_k$ as meshes.  If $\rho$ is any real-valued function defined on a mesh then we will refer to $\rho$ as a mesh function.  Obviously any real valued function defined on a continuum $\mathbb{R}$ or $\mathbb{R} \times \mathbb{R}$ (which we may call continuum functions) can also be considered as a mesh function by restricting its domain to the mesh.  If $\rho$ is a mesh function on $\mathbb{R}_h \times \mathbb{R}_k$ then we will ease some notation by writing $\rho_{n,j} :=  \rho(x_n,t_j)$ and $\rho_j := \rho(\cdot, t_j)$.    \\
\indent We also introduce three discrete derivative operators $D_+$, $D_-$ and $D_0$ that "differentiate" mesh functions $\rho$ defined on $\mathbb{R}_h$ (and hence they can also differentiate continuum functions $\rho$ defined on $\mathbb{R}$).  The operators are given by  
\begin{eqnarray*}
D_+ \rho(x) = \frac{\rho(x+h) - \rho(x)}{h} \hspace{1cm}
D_- \rho(x) = \frac{\rho(x) - \rho(x-h)}{h} \hspace{1cm}
D_0 \rho(x) = \frac{\rho(x+h) - \rho(x-h)}{2h}
\end{eqnarray*} 
\noindent We will also sometimes use shifting operators $E$ and $E^{-1}$ given by 
\begin{eqnarray*}
(E\rho)(x)  =  \rho(x+h) \hspace{3cm} (E^{-1}\rho)(x) = \rho(x-h)
\end{eqnarray*}
\noindent These four operators will only act on the $x$ variable of our functions $\rho(x,t)$.  We will also use the operator $D_{t,+}\eta(t) = \frac{\eta(y+k) - \eta(t)}{k}$.  The following properties of $D_+$, $D_-$, $D_0$, $D_{t,+}$, $E$, and $E^{-1}$ are immediate consequences of their definition:

\begin{enumerate}
\item $D_0 \rho= \frac{1}{2}\left(D_+ \rho + D_- \rho\right)$ 
\item If $\rho = \rho(x,t)$ then the operators $D_+$, $D_-$, $D_0$, $D_{t,+}$ and $E$ all commute when acting on $\rho(x,t)$.
\item $D_+(\nu \cdot \rho)(x_n) = \nu(x_n) D_+ \rho(x_n) + (E \rho)(x_n)D_+\nu(x_n)$
\item $D_-(\nu \cdot \rho)(x_n) = \nu(x_n) D_- \rho(x_n) + (E^{-1}\rho)(x_n)D_-\nu(x_n) $
\item $D_0(\nu \cdot \rho)(x_n) = (E^{-1}\nu)(x_n) D_0 \rho(x_n) + (E\rho)(x_n)D_0\nu(x_n) $
\item If $\rho$ is a continuum function differentiable on $\mathbb{R}$ then for each $x_n \in \mathbb{R}_h$ there exists $x\in\mathbb{R}$ where $x_n \leq x \leq x_{n+1}$ such that we have $D_+ \rho(x_n) = \frac{d}{dx}\rho(x)$.
\end{enumerate}

\indent In order to solve equation $(\ref{daggerdagger})$ we shall consider the following difference scheme, which is a discretization of $(\ref{daggerdagger})$

\begin{eqnarray}
 \left\{ 
\begin{array}{rcl}
D_{t,+}u_{n,j} +  u^2_{n,j} D_0 u_{n,j+1} + D^2_+ D_- u_{n,j+1} + 2f_{n,j} u_{n,j} D_0 u_{n,j+1}\hspace{1cm} \\
 + 2 f_{n,j} (f_x)_{n,j} u_{n,j+1} + (f_x)_{n,j} u_{n,j} u_{n,j+1}  + f^2_{n,j} D_0 u_{n,j+1} + g_{n,j} = 0 \label{star1}\\
\end{array} 
\right. 
\end{eqnarray}

\noindent The advantage behind using this particular discretization lies in the fact that the equation can then be written in a more convenient and concise form.  In order to rewrite this difference scheme we introduce a linear operator $Q_j$ on mesh functions given by 

\begin{equation}
(Q_j \rho)(x_n) := u^2_{n,j} D_0 \rho(x_n) + D^2_+ D_- \rho(x_n) + 2f_{n,j} u_{n,j} D_0 \rho(x_n) + 2 f_{n,j} (f_x)_{n,j} \rho(x_n) + (f_x)_{n,j} u_{n,j} \rho(x_n)  + f^2_{n,j} D_0 \rho(x_n) \label{Q_j}
\end{equation}

\noindent Then $(\ref{star1})$ can be written in a shorter form as 

\begin{eqnarray}
\begin{array}{rcl}
 (I+kQ_j)u_{j+1} &=& u_j - kg_j \label{star2}
\end{array} 
\end{eqnarray}

\noindent where $g_j \in S^{-\infty}\left(\mathbb{R}\right)$ is considered as a mesh function.  The task behind solving $(\ref{star2})$ then is to show that one can invert the operator $I+k Q_j$, at least for some finite amount of time.  The invertibility will be possible only in certain function spaces, therefore we will now introduce an appropriate space.  \\
\indent First we consider the discrete inner products 
\begin{eqnarray*}
\DLIP{u}{v} &=& \sum_{-\infty} ^{\infty}u(x_n)v(x_n) h  \\
\DSIP{u}{v} &=& \DLIP{ \IHN{x}  u}{ \IHN{x} v} + \DLIP{ \IHN{x} D_+^3 u}{ \IHN{x} D_+^3 v} + \DLIP{D_+^5 u}{D_+^5 v} 
\end{eqnarray*} 

\noindent where $\IHN{x} = \sqrt{x^2+1}$, and we define the corresponding norms and Hilbert spaces,
$$\DLNS{u} = \DLIP{u}{u} \hspace{1cm} L^2_h = \left\{ u(x_n)\hspace{1mm} mesh \hspace{1mm} functions \hspace{1mm}on \hspace{1mm}\mathbb{R}_h: \DLN{u} < \infty \right\}  $$ 
$$\DSNS{u} = \DSIP{u}{u} \hspace{1cm} S_h = \left\{ u(x_n) \hspace{1mm} mesh \hspace{1mm} functions \hspace{1mm} on \hspace{1mm} \mathbb{R}_h : \DSN{u} < \infty \right\} $$

\noindent  From the definitions of the $L^2_h$ inner product and its norm we have the properties: 
\begin{enumerate}
\item $\DLN{E\rho}=\DLN{\rho}$
\item If $\rho, \nu, D_- \rho, D_+ \nu \in L^2_h$ then $\DLIP{D_+ \nu}{\rho} = - \DLIP{\nu}{D_-\rho} $
\item For any $j,k,l \in \mathbb{N}$ we have $\DLN{D_+^j D_-^k D_0^l \rho} \leq C\DLN{\rho}$ where $C$ is a constant depending on $h$.
\item For any $j\in\mathbb{N}$ we have $\DLN{\IHN{x}^j D_0 \rho}\leq C \DLN{\IHN{x}^j D_+ \rho}$ where $C$ is independent of $h$.
\item For any $j,l\in\mathbb{N}$ we have $\DLN{D_+^j\IHN{x}^N E^l \rho} \leq C \DLN{\IHN{x}^{N-j} \rho}$ where $C$ is independent of $h$.
\end{enumerate}

\indent The first and second properties follow from simply reindexing and/or rearranging terms in the summation. The third property just requires use of the triangle inequality on each summand of $D_+^j D_-^k D_0^l u$.  For the fourth property we write $D_0$ in terms of $D_+$ and $D_-$ and use the fact that $h\in\left[0,1 \right]$.  For the fifth property we write the definition of $D_+$ and $\IHN{x}$ and use the fact that $h\in \left[0,1\right]$.
\subsection{Preliminary Lemmas}

\noindent The following two propositions provide more basic facts regarding the operators $D_+, D_-, D_0$ and certain basic inequalities that we will use throughout the article sometimes even without reference.  The proofs are found in appendix B.  \\

\noindent \textbf{Proposition 2.1} \hspace{1mm} For any mesh functions $\rho(x_n), \nu(x_n), \xi(x_n)$ we have the following:
\begin{enumerate}
\item If $\rho, D_+^3 \rho \in L^2_h$ then $\DLIP{D_+^2D_- \rho}{\rho} \geq 0$
\item For any $n\in\mathbb{N}$ we have $$D_+^n (\rho \cdot \nu \cdot \xi) = \sum_{i_1+i_2 + i_3 = n}c_{i_1,i_2,i_3} (E^{i_2 + i_3}D_+^{i_1}\rho)\cdot(E^{i_3}D_+^{i_2}\nu)\cdot(D_+^{i_3}\xi) $$ for some constants $c_{i_1,i_2,i_3}\in\mathbb{N}$.
\item If $\nu,\rho\nu, D_0 \nu, (E\nu)( D_+ \rho) \in L^2_h$ then $\DLIP{\rho\nu}{D_0 \nu} = -\frac{1}{2}\DLIP{\nu}{E\nu D_+\rho}$ 
\item If $\nu, \rho \nu, \rho D_+ \nu, \rho D_- \nu, \rho E D_+ \nu, D_+^2 D_-\nu, \nu D_-^2D_+ \rho,D_+\nu, (D_0 \rho)(D_+\nu),(D_- \rho)(D_+\nu) \in L^2_h$ and $\rho(x_n)\geq 0$ for all $x_n \in \mathbb{R}_h$ then $$\DLIP{\rho \nu}{D_+^2D_- \nu} \geq -\frac{1}{2}\DLIP{\nu}{\nu D_-^2D_+ \rho} + \DLIP{D_0 \rho D_+ \nu}{D_+ \nu} + \frac{1}{2}\DLIP{D_- \rho D_+ \nu}{D_+ \nu}$$
\end{enumerate}

\vspace{2mm}

\noindent \textbf{Proposition 2.2} \hspace{1mm} Let $N,n \in \mathbb{N}$, $T>0$, $0<h_1,k_1<1$ and let $g\in S^{-\infty}(\mathbb{R} \times \left[ \left.-c,\infty \right)\right.)$ for some $c>0$.  There exists $C_{N,n} > 0$ such that $\DLN{\IHN{x}^N D_+^n g_j} < C_{N,n}$ for each $0<h\leq h_1$, $0<k\leq k_1$, and $0\leq t_j \leq T$. \\

\indent Another simple but important fact that we will frequently use is the following:\\
\noindent If $c>0$ and $f\in C^{\infty}\left(\mathbb{R}\times \left[ \left.-c,\infty \right)\right. \to \mathbb{R} \right)$ satisfies the property that for every $n\in \mathbb{N}$ and for every compact interval $J \subset \left[ \left.-c,\infty \right)\right.$ we have $\frac{d^n}{dx^n}f(x,t) = O\left(\left| x \right|^{\frac{1}{2}-n} \right)$ uniformly for $t \in J$ then for each $t_j\in J$ we have 
$$\left| \frac{d^n}{dx^n}f(x, t_j ) \right| \leq C \left|\IHN{x}^{\frac{1}{2}-n}\right|$$ where $C>0$ is independent of $k$ and $j$ (but C might depend on $J$).  This statement follows directly from the definitions of $O$ and $\IHN{x}$. \\
\indent The Sobolev inequalities stated below will allow us to prove that the operators $I+kQ_j$ for $j \in \mathbb{N}$ are bounded below and are thus invertible.  These inequalities are stated in [5] but we will state them here and prove them in appendix B for the sake of completeness.  \\
\indent  As a notational remark, from now on we will let $C$ denote a constant whose value might change between consecutive inequalities but the variables that it depends on will often be noted by its indices for example as $C_{n,j,h}$ means some constant depending on $n,j,$ and $h$.\\ 

\noindent \textbf{Lemma 2.3} \hspace{1mm} (Sobolev's inequalities, discrete version)\\
For every $n\in \mathbb{N} $ there exists $C_n>0$ such that for every $h>0$ and for every mesh function $u(x_n)$ defined on $\mathbb{R}_h$ we have
\begin{enumerate}
\item   $\DLN{ D^k_+ u } \leq C_n\left( \DLN{ u} + \DLN{ D^n_+u} \right) $  \hspace{1cm} for $0 \leq k \leq n$
\item   $\SupNorm{ D^k_+ u } \leq C_n \left( \DLN{ u} + \DLN{ D^n_+u} \right)$  \hspace{1cm} for $0 \leq k < n$ 
\end{enumerate}
\vspace{2mm}
\noindent \textbf{Corollary 2.4} \hspace{1mm} For all $N,k,j \in \mathbb{N}$, there exists $C_{N,j,k} > 0$ such that for any $h \in (0,1)$ and for all mesh functions $u(x_n)$ defined on $\mathbb{R}_h$ we have
\begin{enumerate}
\item   $\DLN{ \IHN{x}^N D^j_+ u } \leq C_{N,j,k} \left( \DLN{ \IHN{x}^N u} + \DLN{\IHN{x}^N D^{j+k}_+u} \right) $ \hspace{1cm} for $k\geq 0$
\item  $\SupNorm{ \IHN{x}^N D^j_+ u } \leq C_{N,j,k} \left( \DLN{ \IHN{x}^N u} + \DLN{\IHN{x}^N D^{j+k}_+u} \right) $ \hspace{1cm} for $k\geq 1$
\end{enumerate}

\vspace{2mm}

\indent The next lemma will allow us to prove that the solutions stay bounded for finite time with respect to the Schwartz semi-norms.  The proof can be found in [5].  \\

\noindent \textbf{Lemma 2.5} \hspace{1mm} Suppose $P,Q$ are $C^1\left(\mathbb{R}\right)$, nondecreasing, positive functions, $\Delta t >0$, and for each $j \in \mathbb{N}$ we have $t_j := j \Delta t$.  Let $\eta:\left[ 0, T_0 \right] \rightarrow \mathbb{R}$ be an arbitrary function satisfying \\
\begin{equation}
\frac{\eta_{j+1}-\eta_j}{\Delta t} \leq P(\eta_j) \eta_{j+1} + Q(\eta_j) \label{3.1}
\end{equation}
for each $t_j,t_{j+1} \in \left[0, T_0\right] $ where $\eta_j := \eta(t_j)$, and suppose that $\eta_0 \leq K$ for some $K >0$.  Then there exists $0<T \leq T_0$ and $L,\epsilon >0$ all three depending on $K, P, $ and $Q$ such that if $\Delta t < \epsilon$ then $\eta_j \leq L$ for each $j$ where $t_j \leq T$.  Moreover, if $P$ and $Q$ are constants then we may take $T=T_0$. 
\subsection{Finite Time Existence for Discrete Generalized mKdV in $S_h$}
We will now prove finite time existence for $(\ref{star2})$.  The following lemma is the key estimate for establishing invertibility of the operator $I+kQ_j$ in the space $S_h$.  \\

\noindent \textbf{Lemma 2.6} \hspace{1mm} Suppose $T>0$, $h,k \in \left(0,1 \right)$ and that for each $j$ where $t_j\in \left[0,T\right]$ we have a given mesh function $u_j(x_n)$ defined on $\mathbb{R}_h$.  Define the operators $Q_j$ as in (\ref{Q_j}). Then there exists $C>0$ depending only on $f$ and $T$ but not on $h$, $k$, $j$, or the mesh functions $u_j$ such that for any mesh function $u=u(x_n)$ defined on $\mathbb{R}_h$ the inequality \\
\begin{equation}
\DSIP{Q_j u}{u} \geq -C \DSNS{u } \left(1+ \DSNS{ u_j }\right) \label{4.1}
\end{equation}
holds for each $j$ where $0 \leq t_j \leq T$.\\

\noindent \textbf{Proof of Lemma 2.6} \hspace{1mm} By the definition of $\DSIP{\cdot}{\cdot}$ and $Q_j$ we have that,
\begin{eqnarray*}
\DSIP{Q_j u}{u} & = & \DLIP{\IHN{x}^2 Q_j u}{u} + \DLIP{\IHN{x}^2 D_+^3 (Q_j u)}{D_+^3 u} + \DLIP{D_+^5 (Q_j u)}{D_+^5 u} \nonumber \\ 
                & = & \DLIP{\IHN{x}^2 u_j^2 D_0 u}{u} + \DLIP{\IHN{x}^2 D_+^2 D_- u}{u} + 2 \DLIP{\IHN{x}^2 f_j u_j D_0 u}{u} + \nonumber \\
                &   & + 2\DLIP{\IHN{x}^2 f_j (f_x)_j u}{u} + \DLIP{\IHN{x}^2 (f_x)_j u_j u}{u} + \DLIP{\IHN{x}^2 f_j^2 D_0 u}{u} + \nonumber \\
                &   & + \DLIP{\IHN{x}^2 D_+^3(u_j^2 D_0 u)}{D_+^3 u} + \DLIP{\IHN{x}^2 D_+^3(D_+^2D_- u)}{D_+^3 u} + 2\DLIP{\IHN{x}^2 D_+^3(f_j u_j D_0 u)}{D_+^3 u} + \nonumber \\
                &   & + 2\DLIP{\IHN{x}^2 D_+^3(f_j (f_x)_j u)}{D_+^3 u} + \DLIP{\IHN{x}^2 D_+^3((f_x)_j u_j u)}{D_+^3 u} + \DLIP{\IHN{x}^2 D_+^3(f_j^2 D_0 u)}{D_+^3 u} + \nonumber \\
                &   & + \DLIP{D_+^5(u_j^2 D_0 u)}{D_+^5 u} + \DLIP{D_+^5 D_+^2D_- u}{D_+^5 u} + 2\DLIP{D_+^5(f_j u_j D_0 u)}{D_+^5 u} + \nonumber \\
                &   &  + 2\DLIP{D_+^5(f_j (f_x)_j u)}{D_+^5 u} + \DLIP{D_+^5 ((f_x)_j u_j u)}{D_+^5 u} + \DLIP{D_+^5 (f_j^2 D_0 u)}{D_+^5 u} 
\end{eqnarray*}
We will now show how to bound each term above by the right side of (\ref{4.1}) for some appropriate constant $C$.  Upon adding all the inequalities we will obtain inequality (\ref{4.1}).  For conciseness we shall only write estimates for the first several terms.  The other estimates can be obtained by using the same ideas.  \\

\noindent \underline{\textit{Estimate for Term}} $\DLIP{\IHN{x}^2 u_j^2 D_0 u}{u}$: \\
\begin{eqnarray*}
\DLIP{\IHN{x}^2 u_j^2 D_0 u}{u} & \geq & -\SupNorm{u_j}^2 \DLN{\IHN{x} D_0 u} \DLN{\IHN{x}u} \\
                                & \geq & -C \big(\DLN{u_j} + \DLN{D_+^5 u_j} \big)^2 \DLN{\IHN{x}D_+ u} \DLN{\IHN{x}u} \geq  -C \DSNS{u_j} \DSNS{u}
\end{eqnarray*}

\noindent \underline{\textit{Estimate for Term}} $\DLIP{\IHN{x}^2 D_+^2 D_- u}{u}$: \\
\begin{eqnarray*}
\DLIP{\IHN{x}^2 D_+^2 D_- u}{u} & \geq & -\DLN{\IHN{x}D_+^2D_- u} \DLN{\IHN{x} u} \geq  -C \DLN{\IHN{x}D_+^3 u} \DLN{\IHN{x} u} \geq  -C \DSNS{u} 
\end{eqnarray*}

\noindent \underline{\textit{Estimate for Term}} $\DLIP{\IHN{x}^2 f_j u_j D_0 u}{u}$:\\ 
\begin{eqnarray*}
\DLIP{\IHN{x}^2 f_j u_j D_0 u}{u}  & \geq & -\SupNorm{\IHN{x} u} \DLN{\IHN{x}D_0 u} \DLN{f_j u_j} \geq -C\left( \DLN{\IHN{x}u} + \DLN{\IHN{x}D_+^3 u} \right)^2 \DLN{\IHN{x}u_j} \\
                                   & \geq & -C \DSNS{u}\DSN{u_j} \geq -C \DSNS{u}\left(1+\DSN{u_j}^2 \right)
\end{eqnarray*}

\noindent \underline{\textit{Estimate for Term}} $\DLIP{\IHN{x}^2 f_j (f_x)_j u}{u}$: \\
\begin{eqnarray*}
\DLIP{\IHN{x}^2 f_j (f_x)_j u}{u}  & \geq & -\SupNorm{f_j (f_x)_j} \DLNS{\IHN{x} u} \geq -C\SupNorm{\IHN{x}^{\frac{1}{2}} \IHN{x}^{-\frac{1}{2}}} \DLNS{\IHN{x} u} \geq -C \DSNS{u}
\end{eqnarray*}

\noindent \underline{\textit{Estimate for Term}} $\DLIP{\IHN{x}^2 (f_x)_j u_j u}{u}$: \\ 
\begin{eqnarray*}
\DLIP{\IHN{x}^2 (f_x)_j u_j u}{u}  & \geq & -\DLNS{\IHN{x} u} \SupNorm{(f_x)_j}\SupNorm{u_j} \geq -C\DSNS{u} \SupNorm{\IHN{x}^{-\frac{1}{2}}} \SupNorm{ u_j} \\
                                   & \geq & -C \DSNS{u} \left( \DLN{u_j}+\DLN{D_+^5 u_j }\right) \geq -C \DSNS{u} \DSN{u_j} \geq -C\DSNS{u} \big( 1+\DSNS{u_j} \big)
\end{eqnarray*}

\noindent \underline{\textit{Estimate for Term}} $\DLIP{\IHN{x}^2 f_j^2 D_0 u}{u}$: \\
\begin{eqnarray*}
\DLIP{\IHN{x}^2 f_j^2 D_0 u}{u} \geq  -\frac{1}{2}\SupNorm{\frac{1}{\IHN{x}^2} D_+\left(\IHN{x}^2 f_j^2\right)}\DLNS{\IHN{x} u} \geq -C \left( \SupNorm{\frac{ f_j^2}{\IHN{x}^2} D_+\IHN{x}^2 } + \SupNorm{\frac{\IHN{Ex}^2}{\IHN{x}^2}} \SupNorm{D_+ f_j^2} \right) \DSNS{u} \hspace{1cm} \\
\geq  -C \left( \SupNorm{\frac{\IHN{x}}{\IHN{x}^2} D_+\left(\IHN{x}^2\right) } + \SupNorm{\frac{d}{dx}(f_j^2)} \right) \DSNS{u} \geq -C \left(1+ \SupNorm{(f_j)(f_x)_j} \right)\DSNS{u} \geq -C\DSNS{u} \hspace{13mm}
\end{eqnarray*}

\noindent \underline{\textit{Estimate for Term}} $\DLIP{\IHN{x}^2 D_+^3(u_j^2 D_0 u)}{D_+^3 u}$:  \\  
By the product rule for $D_+$ we have,
\begin{equation*}
\DLIP{\IHN{x}^2 D_+^3 (u_j^2 D_0 u)}{D_+^3 u} = \sum_{i_1+i_2+i_3 = 3} c_{i_1,i_2,i_3}\DLIP{\IHN{x}(E^{i_2 + i_3}D_+^{i_1} u_j)(D_+^{i_2}u_j)(E^{i_2} D_+^{i_3}D_0 u)}{D_+^3 u}
\end{equation*}
and if we assume that $i_1 \leq i_2$ then we see that,
\begin{eqnarray*}
\DLIP{\IHN{x}(E^{i_2 + i_3}D_+^{i_1} u_j)(D_+^{i_2}u_j)(E^{i_2} D_+^{i_3}D_0 u)}{D_+^3 u} \geq  -\SupNorm{E^{i_2}D_+^{i_3}D_0 u} \SupNorm{E^{i_2 + i_3}D_+^{i_1}u_j} \DLN{\IHN{x}D_+^{i_2} u_j} \DLN{\IHN{x} D_+^3 u} \\
                                              \geq  -C\big( \DLN{u} + \DLN{D_+^5 u} \big) \big( \DLN{u_j} + \DLN{D_+^5 u_j} \big) \left(\DLN{\IHN{x} u_j} + \DLN{\IHN{x} D_+^3 u_j} \right) \DSN{u}
    \geq  -C \DSNS{u} \DSNS{u_j}
\end{eqnarray*}

\noindent \underline{\textit{Estimate for Term}} $\DLIP{\IHN{x}^2 D_+^3(D_+^2D_- u)}{D_+^3 u}$: 
\begin{eqnarray*}
\DLIP{\IHN{x}^2 D_+^3(D_+^2D_- u)}{D_+^3 u}   \geq  -\frac{1}{2} \DLIP{D_+^3 u}{D_+^3 uD_-^2D_+ \IHN{x}^2} + \DLIP{D_0 \IHN{x}^2 D_+^4 u}{D_+^4 u} + \frac{1}{2} \DLIP{D_-\IHN{x}^2 D_+^4 u}{D_+^4 u} \hspace{.5cm} \\
																						  \geq  -\frac{1}{2} \DLNS{D_+^3 u} \SupNorm{D_-^2D_+ \IHN{x}^2} - \DLIP{D_-(D_0 \IHN{x}^2 D_+^4u)}{D_+^3 u} - \frac{1}{2} \DLIP{D_-(D_-\IHN{x}^2 D_+^4 u)}{D_+^3 u} \hspace{11mm}\\ 
																						  \geq  -C  \DLNS{\IHN{x}D_+^3 u} -C \SupNorm{D_+^2 \IHN{x}^2} \DLN{D_+^4 u} \DLN{D_+^3 u} -C \DLN{D_+^5 u}\DLN{D_+\IHN{x}^2 D_+^3 u}\hspace{23mm}\\
																						  \geq  -C \DSNS{u} -C \left(\DLN{u} + \DLN{D_+^5 u} \right)^2 -C\left(\DLN{u} + \DLN{D_+^5 u} \right) \DLN{\IHN{x}D_+^3 u} \geq -C \DSNS{u} \hspace{2cm}
\end{eqnarray*}

\noindent \underline{\textit{Estimate for Term}} $\DLIP{\IHN{x}^2 D_+^3(f_j u_j D_0 u)}{D_+^3 u}$: \\ 
By the product rule for $D_+$ we have,
\begin{equation*}
\DLIP{\IHN{x}^2 D_+^3 (f_j u_j D_0) }{D_+^3 u} = \sum_{i_1 + i_2 + i_3 = 3} c_{i_1,i_2,i_3} \DLIP{\IHN{x}^2 (E^{i_2 + i_3}D_+^{i_1} f_j)(D_+^{i_2} u_j)(E^{i_2}D_+^{i_3}D_0 u)}{D_+^3 u}
\end{equation*}
and now we will bound each term.\\
For $0 \leq i_3 \leq 1$ we have,
\begin{eqnarray*}
\DLIP{\IHN{x}^2 (E^{i_2 + i_3}D_+^{i_1} f_j)(D_+^{i_2} u_j)(E^{i_2}D_+^{i_3}D_0 u)}{D_+^3 u} \geq  -\SupNorm{\IHN{x}E^{i_2}D_+^{i_3}D_0 u} \DLN{ (E^{i_2 + i_3}D_+^{i_1} f_j)(D_+^{i_2}u_j)} \DLN{\IHN{x} D_+^3 u} \\
                                                                         \geq  -C \SupNorm{\IHN{x}D_+^{i_3 +1}u} \DLN{\IHN{x}D_+^{i_2} u_j} \DSN{u}  
                                                                          \geq  -C \left( \DLN{\IHN{x} u} + \DLN{\IHN{x}D_+^3 u} \right) \left(\DLN{\IHN{x} u_j} + \DLN{\IHN{x} D_+^3 u_j} \right) \DSN{u}\\
                                                                          \geq -C\DSNS{u}\DSN{u_j} \geq -C\DSNS{u}\big( 1 + \DSNS{u_j} \big)\hspace{106mm}
\end{eqnarray*}
For $i_3 = 2$ we have,
\begin{eqnarray*}
\DLIP{\IHN{x}^2 (E^{i_2+2}D_+^{i_1} f_j)(D_+^{i_2} u_j)(E^{i_2}D_+^{2}D_0 u)}{D_+^3 u}  \geq  -\SupNorm{ (E^{i_2+2}D_+^{i_1} f_j)(D_+^{i_2}u_j) } \DLN{\IHN{x}E^{i_2}D_+^2 D_0 u} \DLN{\IHN{x}D_+^3 u} \hspace{1cm}\\
                                                                          \geq  -C\SupNorm{\IHN{x}D_+^{i_2} u_j} \DLNS{\IHN{x}D_+^3 u} 
                                                                          \geq   -C\big( \DLN{\IHN{x} u_j} + \DLN{\IHN{x}D_+^3 u_j} \big) \DSNS{u} 
                                                                          \geq -C \DSNS{u}\big(1 + \DSNS{u_j} \big)\hspace{1cm}
\end{eqnarray*}
For $i_3 = 3$ we have,
\begin{eqnarray*}
\DLIP{\IHN{x}^2 (E^3 f_j) u_j(D_+^3 D_0 u)}{D_+^3 u}  \geq  -\frac{1}{2}\SupNorm{\frac{1}{\IHN{x}^2}D_+(\IHN{x}^2 (E^3 f_j) u_j)} \DLNS{\IHN{x}D_+^3 u} \hspace{2cm}\\
                                                                          \geq  -C \left( \SupNorm{\IHN{x}^{-1} (E^3 f_j) u_j} + \SupNorm{(E^3 D_+ f_j)u_j } + \SupNorm{(E^4 f_j) D_+ u_j} \right) \DSNS{u}\\
                                                                          \geq  -C \left( \SupNorm{\IHN{x}^{-1} \IHN{x}^{\frac{1}{2}} u_j} + \SupNorm{(f_x)_j} \SupNorm{u_j} + \SupNorm{\IHN{x} u_j} \right) \DSNS{u} \hspace{1.5cm}\\
                                                                          \geq -C \SupNorm{\IHN{x} u_j}\DSNS{u} 
                                                                          \geq -C \DSNS{u} \left(1 + \DSNS{u_j} \right) \hspace{39mm}
\end{eqnarray*}

$\square$ \\

\noindent \textbf{Lemma 2.7} \hspace{1mm} Suppose $K >0$ and $u_0\in S^{-\infty}(\mathbb{R})$ satisfies the property that $\DSN{u_0 } \leq K$ for each $h \in (0,1)$.  Then there exists $T,L, \epsilon >0$ depending only on $K$ such that if $k \in (0, \epsilon)$ and $h\in (0,1)$ then the difference scheme $(\ref{star2})$  may be solved for each mesh function $u_j$ with $t_j\in \left[0,T\right]$.  Moreover, we have that $\DSN{ u_j } \leq L$ for each $j$ where $t_j \in \left[0, T\right]$.  \\

\noindent \textbf{Proof of Lemma 2.7} \hspace{1mm} Choose $T_0 >0$ arbitrarily.  Assume for now that the mesh functions $u_j$ are known for each $h,k\in \left(0,1\right)$ and for $0 \leq t_j \leq T_0$.  We will first construct the aformentioned $T,L$ and an $ \epsilon_0 >0$ and show that the mesh functions $u_j$ whose time mesh size satisfies $k \in (0, \epsilon_0) $ will satisfy the inequality $\DSN{ u_j } \leq L$ for each $0\leq t_j \leq T$. \\
\indent Assume that $0\leq t_j, t_{j+1} \leq T$.  Taking inner product of $(\ref{star2})$ with $u_{j+1}$ we obtain 
$$\DSIP{(I+kQ_j)u_{j+1}}{u_{j+1}} = \DSIP{u_{j+1}}{u_j} - k\DSIP{g_j}{u_{j+1}}$$
and we may use Cauchy-Schwarz inequality on the right side and simply rewrite the left side to obtain 
\begin{equation}
\DSNS{ u_{j+1} } + k\DSIP{Q_j u_{j+1}}{u_{j+1}} \leq \DSN{ u_{j+1}} \cdot \DSN{ u_j } + k \DSN{ g_j } \cdot \DSN{ u_{j+1} } \label{L5.1}
\end{equation}
By lemma 2.6 we may choose $C>0$ such that 
\begin{equation}
\DSNS{ u_{j+1}} + k\DSIP{Q_j u_{j+1} }{u_{j+1}} \geq \DSNS{ u_{j+1} } \left[ 1-kC(\DSNS{ u_j } +1) \right] \nonumber
\end{equation}
Also, since $g \in S^{-\infty}\left(\mathbb{R}\times \left[\left. -c,\infty \right) \right. \right)$ we have that $\sup_{t\in \left[ 0, T_0 \right] } \DSN{ g(\cdot, t) } < \infty$.  We may then enlarge $C$ so that $C> \sup_{t\in \left[ 0, T_0 \right] } \DSN{ g(\cdot,t ) }$, which is clearly still independent of $h$ and $j$.  By combining $( \ref{4.1} )$ and $( \ref{L5.1} )$ we thus obtain for $t_j,t_{j+1} \in \left[ 0,T_0 \right]$ that 
\begin{eqnarray*}
\DSNS{ u_{j+1} } \big[ 1-kC(\DSNS{ u_j } +1) \big] & \leq &  \DSN{ u_{j+1} } \cdot \DSN{ u_j } + k \DSN{g_j } \cdot \DSN{ u_{j+1} } \leq \DSN{ u_{j+1} } \left( \DSN{ u_j} + k C\right)
\end{eqnarray*}
or equivalently, 
\begin{equation}
\frac{\DSN{ u_{j+1} } - \DSN{ u_j }}{k} \leq C \left( \DSNS{ u_j } + 1\right) \DSN{ u_{j+1} } + C \nonumber
\end{equation}
 which is an inequality of the form (\ref{3.1}).  Then by lemma 2.5 there exists $0<T\leq T_0$ and $L,\epsilon_0 >0$ depending  on $K \geq \DSN{u_0 }$ such that if $k \in (0, \epsilon_0) $ then $\DSN{ u_j } \leq L$ for each $j$ where $0 \leq t_j\leq T$.  The $T,L,\epsilon_0$ are independent of $h\in (0,1)$ because the constant $C$ is independent of $h$.  Moreover, the $T, L,$ and $\epsilon_0$ depend only on $K, C, $ and $P(v)$ :=$v+1$ by lemma 2.5.  Since $K$ is given and $C$ is determined by the given functions $f$ and $g$ and on the value of $T_0$ we may construct $T,L,$ and $\epsilon_0$ without assuming that $u_j$ is constructed for $0 \leq t_j \leq T_0$.  \\
\indent Given $T,L$ as constructed above it suffices to show that there exists $ 0 < \epsilon \leq \epsilon_0$ such that for any $h\in(0,1)$ and $k \in (0,\epsilon)$ the difference scheme $(\ref{star2})$ may be solved for $u(x_n,t_j)$ where $ (x_n,t_j)\in \mathbb{R}_h\times \left(\mathbb{R}_k \cap \left[0,T\right] \right)$ - the desired bound $\DSN{ u(\cdot,t_j) } \leq L$ would follow automatically by our construction of $T,L,$ and $\epsilon_0$ given above.\\
\indent Choose $\epsilon >0$ so that $\epsilon C (L^2 +1) < \frac{1}{2}$ and $0<\epsilon \leq \epsilon_0$ and fix values for $h \in (0,1)$ and $k\in (0,\epsilon)$.  
\noindent  Suppose $u_0,u_1,\ldots,u_j$ are known for some $j \geq 0$.  We will show that one may construct $u_{j+1}$ as long as $t_{j+1} \leq T$.  Define an operator $P_j$ := $I+kQ_j$.  Then by lemma 2.6 we have for any mesh function $u=u(x_n)$ 
\begin{equation}
   \DSIP{P_j u}{u}  \geq  \left[1-kC(\DSNS{ u_j }+1) \right]\DSNS{ u } 
     \geq  \left[ 1-\epsilon C(L^2+1)\right]\DSNS{ u } 
     \geq  \frac{1}{2} \DSNS{ u } \label{5.4}
\end{equation}
from which it easily follows that $P_j$ is injective as an operator on mesh functions.  By choosing an appropriate domain $\mathcal{D}$ we may consider $P_j$ as a linear (possibly unbounded) operator on $S_h$.    To this end we define $\mathcal{D}$ := $\left\{u \in S_h : P_j u \in S_h \right\}$ and so we have 
$$\begin{CD}
S_h \supset \mathcal{D} @>P_j>> S_h
\end{CD}$$ 
$P_j$ would be an unbounded operator on $S_h$ if $f$ is unbounded.  \\
\indent In order to solve the difference scheme $(\ref{star2})$  it is enough to show that the operators $P_j$ are surjective because then $P_j :\mathcal{D} \rightarrow S_h$ would be a bijection so that we could define $u_{j+1} := P_j^{-1}(u_j-kg_j)$ for $t_{j+1} \leq T$.  Let $\phi(x) \in S_h$.  We will construct a preimage of $\phi$ with the aid of the operators $P^R_j$ which are bounded versions of $P_j$ that are bijective.\\
\indent We introduce (bounded) operators $P_j^R$: $S_h \rightarrow S_h $ for each $R>0$ given by  
\[
\begin{array}{rcl}
P_j^R \rho &:=& \rho +k\left(u^2_j D_0 \rho + D^2_+ D_- \rho + 2f_j u_j D_0 \rho + 2 f_j (f_x)_j \rho + (f_x)_j u_j \rho  + (f^2_R)_j D_0 \rho \right)
\end{array} 
\]
where $(f_R)_j$ is a bounded function agreeing with $f_j$ on $\left[ -R,R \right]$ but truncated to zero outside $\left[ -2R,2R \right]$.  Formally, we take $\psi(x) \in C_c^{\infty}\left(\mathbb{R}\right)$ where $\psi(x) = 1$ for $\left| x \right|<1$ and $\psi(x) = 0$, $\left| x \right|>2$, and $\psi(x) \leq 1$ for all $x$ and define $(f_R)_j :=f_j \cdot \psi(x/R)$.  Then clearly $|(f_R)_j(x)| \leq |f_j(x)|$ for each $x$.  We will now prove that operators $P_j^R$ are bounded and bijective on $ S_h$.  \\
 
\noindent \textit{Claim 1} \hspace{1mm} The maps $P_j^R$: $S_h \rightarrow S_h $ are bounded for each $R>0$ \\
\noindent \textit{proof of claim 1} \hspace{1mm} The boundedness of $P_j^R$ follows from the below estimates: \\
\begin{equation*}
\DSN{P_j^R \rho} \leq  \DSN{\rho} +k\left( \DSN{u^2_j D_0 \rho} + \DSN{D^2_+ D_- \rho} + 2\DSN{f_j u_j D_0 \rho} 
+ 2 \DSN{f_j (f_x)_j \rho} + \DSN{(f_x)_j u_j \rho}  + \DSN{(f^2_R)_j D_0 \rho} \right) 
\end{equation*}

\noindent The terms in parenthesis can each be bounded by lemma 2.3 and the triangle inequality.  For the first term we have
\begin{eqnarray*}
\DSN{u^2_j D_0 \rho} & \leq &  \DLN{\left \langle x\right \rangle  u^2_j D_0 \rho} + \DLN{\left \langle x\right \rangle D_+^3(u^2_j D_0 \rho)} + \DLN{D_+^5( u^2_j D_0 \rho)}  \\
                     & \leq &  \DLN{\left \langle x\right \rangle  u^2_j D_0 \rho} + C_h\DLN{\left \langle x\right \rangle u^2_j D_0 \rho}   + C_h \DLN{u_j^2 D_0 \rho} \leq  C_h \DLN{\left \langle x\right \rangle  u^2_j D_0 \rho}  \\
                     & \leq &  C_h \SupNorm{u_j}^2 \DLN{\left \langle x\right \rangle D_0 \rho} \leq  C_h \left(\DLN{u_j}+ \DLN{D_+ u_j} \right)^2 \DLN{\left \langle x\right \rangle D_+ \rho} \leq C_{h,j} \DSN{\rho} 
\end{eqnarray*}
\noindent and for the last term we have,

\begin{eqnarray*}
\DSN{(f^2_R)_j D_0 \rho} & \leq & \DLN{\IHN{x} (f^2_R)_j D_0 \rho} + \DLN{\IHN{x} D_+^3((f^2_R)_j D_0 \rho)} +\DLN{ D_+^5((f^2_R)_j D_0 \rho)} \\
                         & \leq & \DLN{\IHN{x} (f^2_R)_j D_0 \rho} + C_h \DLN{\IHN{x}(f^2_R)_j D_0 \rho} + C_h \DLN{(f^2_R)_j D_0 \rho} \leq C_h \DLN{\IHN{x}(f^2_R)_j D_0 \rho}  \\ 
                         & \leq & C_h \SupNorm{(f^2_R)_j }\DLN{\IHN{x} D_0 \rho} \leq C_{h,R}\DLN{\IHN{x} D_+ \rho} \leq C_{h,R}\DSN{\rho} 
\end{eqnarray*}

\noindent  The other terms can be bounded in a similar way.  This concludes the proof of claim 1.\\

\indent By tracing the estimates of lemma 2.6 we also see that $P^R_j$ satisfies estimate $ ( \ref{5.4} )$.  We omit the details here but the important point here is that $C$ can be taken independent of $R$.  Thus we would obtain 
\begin{equation}
\DSIP{P^R_j u}{u} \geq \frac{1}{2} \DSNS{u} \label{5.5}
\end{equation}

\noindent Therefore we see that the operators $P^R_j$ are bounded injective operators. \\

\noindent \textit{Claim 2} \hspace{1mm} The image of $P_j^R$: $S_h \rightarrow S_h $ is closed.   \\
\noindent \textit{proof of claim 2} \hspace{1mm} Suppose $P^R_j u_n \rightarrow v$ as $n\rightarrow \infty$.  We will construct $u\in S_h$ such that $P^R_j u = v$.  \\
\noindent Inequality (\ref{5.5}) with $u_n-u_m$ implies, by using the Cauchy-Schwarz inequality on the left side, that for any $m,n\in \mathbb{N}$ we have $$ \DSN{P^R_j (u_n-u_m)} \geq \frac{1}{2} \DSN{u_n-u_m} $$ which implies that the sequence $u_n$ is Cauchy, therefore by completeness of $S_h$ the sequence $u_n$ converges to some $u\in S_h$.  Moreover,
\begin{eqnarray*}
\DSN{P^R_j u - v} & \leq & \DSN{P^R_j u - P^R_j u_n} + \DSN{P^R_j u_n - v} \leq \left\|P^R_j \right\| \DSN{u-u_n}+ \DSN{P^R_j u_n - v}
\end{eqnarray*}  
\noindent and hence for $n$ sufficiently large the right side can be made arbitrarily small so that $P^R_j u = v$.  This concludes the proof of claim 2.  \\

\noindent \textit{Claim 3} \hspace{1mm}  The operators $P_j^R$: $S_h \rightarrow S_h $ are bijective.   \\
\noindent \textit{proof of claim 3} \hspace{1mm} Since injectivity follows from (\ref{5.5}) it is enough to prove surjectivity.  Since the image of $P_j^R$ is closed we have that $S_h = Im (P_j^R) \oplus Im (P_j^R)^{\bot_{S_h}}$.  Suppose there exists $v \in Im (P_j^R)^{\bot_{S_h}}$.  Then $0 =  \DSIP{P_j^R v}{v}   \geq \frac{1}{2} \DSNS{v}  $ which implies that $v\equiv 0$.  This proves surjectivity.   \\  
 
\indent Now we define functions $u_R(x) := (P^R_j)^{-1}\phi$.  Then for each $x_n$ we have by (\ref{5.5}) 

\begin{equation}
\left| u_R(x_n)\right|   \leq  \SupNorm{u_R} \leq \DSN{u_R} \leq  \DSN{\left(P_j^R\right)^{-1}}\cdot \DSN{\phi} \leq  2 \DSN{\phi} \label{5.6}
\end{equation}

\noindent so that the set $\left\{u_R(x_n)\right\}$ is uniformly bounded in $x_n$ and in $R$.  Since the domain of $u_R$ is countable we may find a subsequence $\left\{u_{R_i}\right\}_{i=0}^{\infty}$ (where of course $R_i \rightarrow \infty$ as $i \rightarrow \infty $) that converges pointwise to some function $u$, i.e. for each $x_n\in \mathbb{R}_h$ we have $u_{R_i}(x_n)\rightarrow u(x_n)$ as $i \rightarrow \infty$.  \\

\noindent \textit{Claim 4} \hspace{1mm} The limit function $u$ belongs to the space $S_h$. \\
\noindent \textit{proof of claim 4} \hspace{1mm} By Fatou's Lemma and (\ref{5.6}) we have

$$\sum_n u(x_n)^2 h  = \sum_n \liminf_{i\rightarrow \infty} u_{R_i}(x_n)^2 h \leq \liminf_{i \rightarrow \infty} \DLNS{u_{R_i}} \leq \liminf_{i \rightarrow \infty} \DSNS{u_{R_i}} \leq  4 \DSNS{\phi} $$

\noindent similar inequalities hold for $\IHN{x}u$, $\IHN{x}D_+^3u$ and $D^5_+ u$ and thus we may conclude by adding these estimates that the limit function $u$ lies in $S_h$.  \\

\indent We will now show that $P_j u = \phi$.\\
\indent Fix a point $x_n \in \mathbb{R}_h $.  Then by definitions of $P_j$ and $P^R_j$, the triangle inequality, and the fact that $|(f_R)_j(x)| \leq |f_j(x)|$ we have the following estimate 
\begin{eqnarray}
\left|P_j u(x_n) - \phi(x_n) \right|  & \leq & \left| P_j u(x_n) - P_j u_{R_i}(x_n)\right| + \left| P_j u_{R_i}(x_n) - P_j^{R_l}u_{R_i}(x_n) \right|
                                         + \left| P_j^{R_l}u_{R_i}(x_n) - P_j^{R_l}u_{R_l}(x_n)\right| \label{5.7} \\
                                      & \leq &  \Big[ \left| u(x_n) - u_{R_i}(x_n)\right|+k \left|u^2_jD_0 \left( u(x_n) - u_{R_i}(x_n)\right)\right| +  
                                          k \left|D^2_+D_-\left( u(x_n) - u_{R_i}(x_n)\right)\right|   \nonumber \\
                                      &      & + 2k \left|f_j u_jD_0\left( u(x_n) - u_{R_i}(x_n)\right) \right| + 
                                                2k \left|f_j (f_x)_j\left( u(x_n) - u_{R_i}(x_n)\right)\right| \nonumber \\
                                      &      &  + k \left|(f_x)_j u_j D_0\left( u(x_n) - u_{R_i}(x_n)\right) \right| + 
                                         k \left|(f^2)_j D_0\left(u(x_n) - u_{R_i}(x_n)\right) \right| \nonumber \Big] \\
                                      &      & +\Big[ \left|\left((f^2_{R_l})_j - f^2_j\right) D_0 u_{R_i}(x_n)\right| \Big] \nonumber \\
                                      &      & + \Big[ \left| u_{R_i}(x_n)-u_{R_l}(x_n)\right| + 
                                             k \left|u^2_j D_0\left( u_{R_i}(x_n) - u_{R_l}(x_n)\right) \right| + 
                                         k \left|D^2_+D_-\left( u_{R_i}(x_n) - u_{R_l}(x_n)\right)\right| \nonumber \\
                                      &      & + 2k \left|f_j u_j D_0\left( u_{R_i}(x_n) - u_{R_l}(x_n)\right)\right|+ 
                                         2k \left|f_j (f_x)_j \left( u_{R_i}(x_n) - u_{R_l}(x_n)\right) \right| \nonumber \\
                                      &      &  + k \left|(f_x)_j u_j D_0\left(u_{R_i}(x_n)- u_{R_l}(x_n)\right)\right| 
                                         +k \left|f^2_j D_0\left( u_{R_i}(x_n) - u_{R_l}(x_n)\right) \right| \Big] \nonumber              \end{eqnarray}
                           
\noindent Since $ \left\{ u_{R_i}(x_n) \right\}_{i=0}^{\infty} $ is a Cauchy sequence converging to $u(x_n)$ it follows that each term contained in the first and third pair of brackets can be made arbitrarily small for sufficiently large values of $i$ and $l$.  Therefore, given $\delta >0$ we may find $R_0$ such that if $R_i,R_l\geq R_0$ then the first and third terms of $( \ref{5.7} )$ are less than $\delta/3$.  Moreover, since $f_{R_l}(x_n) \rightarrow f(x_n)$ as $l \rightarrow \infty$ we may fix $R_i\geq R_0$ and find $R_1$ such that if $R_l \geq R_1$ then the middle term of (\ref{5.7}) $\left| (f(x_n) - f_{R_l}(x_n))D_0 u_{R_i}(x_n) \right|$ is less than $\delta/3$.  Therefore we obtain that $P_ju(x_n) = \phi(x_n)$ for each $x_n \in \mathbb{R}_h$ and since $\phi$ lies in $S_h$ we have by definition of $\mathcal{D}$ that $u \in \mathcal{D}$.  This shows that $P_j : \mathcal{D} \rightarrow S_h$ is a bijection.  Therefore we may solve the difference scheme $(\ref{star2})$ by defining $u_{j+1} = P_j^{-1}(u_j - kg_j)$ as long as $t_j \leq T$ (recall that $(\ref{5.4})$ fails for $t_j > T$ so that $P_j$ would not be invertible after time $T$) however the desired bound $\DSN{u_{j+1}}\leq L$ would be true only if $t_{j+1} \leq T$ . $\square$ 

\section{Estimates for the Discrete Solutions}
\subsection{Schwartz Boundedness of Discrete Solutions}
\noindent In this section we will show that the solutions to the discrete equation $(\ref{star2})$ constructed in lemma 2.7 are bounded in all discrete Schwartz norms $\DLN{\IHN{x }^N D^n_+ \cdot}$.  This will follow by some induction arguments shown in next three lemmas.  These lemmas are analogous to those presented by Bondareva in [1].  \\

\noindent \textbf{Lemma 3.1} \hspace{1mm}  Let $u_0\in S^{-\infty}(\mathbb{R})$ and $K,T,\epsilon>0$ all be given as in lemma 2.7 and let $N,n \in \mathbb{N}$.   Suppose that there exists a constant $C_{N,n}>0$ independent of $h \in (0,1)$, $k \in (0,\epsilon)$, and $t_j\in \left[0,T\right]$ such that if $u(x_n,t_j)$ is the solution of the difference scheme $(\ref{star2})$ defined on $\mathbb{R}_h \times \left( \mathbb{R}_k \cap \left[0,T\right] \right)$ with initial condition $u_0 $ then for $0 \leq t_j,t_{j+1} \leq T$ we have 

\begin{eqnarray}
 \DLIP{ \IHN{x}^{2N} D^n_+ (Q_j u_{j+1})}{D^n_+ u_{j+1} } & \geq &  -C_{N,n} \left( \DLNS{ \IHN{x}^N D^n_+ u_{j+1}} + \DLNS{\IHN{x}^N D^n_+ u_j} + 1 \right) \hspace{.5cm} \label{L6.1}
\end{eqnarray}
Then there exists $0 < k_0 \leq \epsilon$ depending on $N$ and $n$ and there exists $C_{N,n}>0$ such that if $h \in (0,1)$, $k\in (0,k_0)$, and $u(x_n,t_j)$ is the solution of the difference scheme $(\ref{star2})$ defined on $\mathbb{R}_h \times \left( \mathbb{R}_k \cap \left[0,T\right] \right)$ with initial condition $u_0 $ then we have $\DLN{\IHN{x}^N D^n_+u_j}<C_{N,n}$ for $t_j \in \left[0, T\right]$.  Moreover, the constant $C_{N,n}$ is independent of the choice of $h$ and $k$. \\

\noindent \textbf{Proof of Lemma 3.1} \hspace{1mm} First let us fix values for $h\in (0,1)$ and $k \in (0,\epsilon)$ and let $u(x_n,t_j)$ be the solution of the difference scheme $(\ref{star2})$ defined on $\mathbb{R}_h \times \left( \mathbb{R}_k \cap \left[0,T\right] \right)$ with initial condition $u_0 $.  If $t_j, t_{j+1}\in \left[0,T\right]$ then we may apply $D_+^n$ to both sides of $(\ref{star2})$ and take the inner product $\DLIP{\IHN{x}^{2n}\cdot}{\cdot}$ of the resulting equation with $D_+^n u_{j+1}$ and use proposition 2.2 to obtain 
\begin{eqnarray}
\DLIP{\IHN{x}^{2N}D_+^n \big(I+kQ_j \big) u_{j+1}}{D_+^n u_{j+1}}  =   \DLIP{\IHN{x}^{2N}D_+^n u_j}{D_+^n u_{j+1}} - k\DLIP{\IHN{x}^{2N}D_+^n g_j}{D_+^n u_{j+1}} \nonumber \\
                                                                   \leq \DLN{\IHN{x}^N D_+^n u_{j+1}} \DLN{\IHN{x}^N D_+^n u_j} + k \DLN{\IHN{x}^N D_+^n u_{j+1}} \DLN{\IHN{x}^N D_+^n g_j} \hspace{14mm} \nonumber \\
                                                                  \leq \frac{1}{2}\DLNS{\IHN{x}^N D_+^n u_{j+1}} + \frac{1}{2} \DLNS{\IHN{x}^N D_+^n u_j} + k C_{N,n} \left( \DLNS{\IHN{x}^N D_+^n u_{j+1}} +1\right) \hspace{10.5mm} \label{L6.2}                                                              
\end{eqnarray}

\noindent We may then continue the left side of $(\ref{L6.2})$ by using $(\ref{L6.1})$ to obtain 
\begin{eqnarray}
\DLIP{\IHN{x}^{2N}D_+^n \big(I+kQ_j \big) u_{j+1}}{D_+^n u_{j+1}} = \DLNS{\IHN{x}^N D_+^n u_{j+1}} + k\DLIP{\IHN{x}^{2N}D_+^n Q_j u_{j+1}}{D_+^n u_{j+1}} \hspace{22mm} \nonumber \\  
  \geq \DLNS{\IHN{x}^N D_+^n u_{j+1}} - k C \left( \DLNS{\IHN{x}^N D_+^n u_j} + \DLNS{\IHN{x}^N D_+^n u_{j+1}} + 1 \right) \label{L6.3}
\end{eqnarray}
By combining $(\ref{L6.2})$ and $(\ref{L6.3})$ we obtain
\begin{eqnarray*}
 \frac{1}{k}\left(\DLNS{\IHN{x}^N D_+^n u_{j+1}} - \DLNS{\IHN{x}^N D_+^n u_j} \right) 
 \leq  C \left( \DLNS{\IHN{x}^N D_+^n u_{j+1}} + \DLNS{\IHN{x}^N D_+^n u_j}  + 1 \right)
\end{eqnarray*}
which is an inequality of the form (\ref{3.1}).  Now we may invoke lemma 2.5 to obtain the existence of $k_0 >0$ and $C_{N,n}>0$ such that if $k \in (0,k_0)$ and $t_j \in \left[0,T\right]$ then $$\DLN{\IHN{x}^N D_+^n u_j} < C_{N,n} $$ and since $C$ is independent of $h$ and $k$ we also have that $C_{N,n}$ is independent of $h$ and $k$.  Thus we have proven that for each $0<h\leq 1$, $0<k\leq k_0$, $0 \leq t_j \leq T$ and solution $u(x_n,t_j)$ of (\ref{star2}) defined on $\mathbb{R}_h\times\left( \mathbb{R}_k \cap \left[0,T\right]\right)$ with initial data $u_0$ we have the inequality $\DLN{\IHN{x}^N D_+^n u_j} < C_{N,n} $.  Moreover since $C$ depends on $N$ and $n$ it follows that $k_0$ depends also on $N$ and $n$. $\square$ \\

\noindent \textbf{Lemma 3.2} \hspace{1mm}  Let $u_0\in S^{-\infty}(\mathbb{R})$ and $K,T,L,\epsilon>0$ all be given as in lemma 2.7 and let $n \in \mathbb{N}$.  Then there exists $0 < k_0 \leq \epsilon$ depending on $n$ and there exists $C_n >0$ such that if $h \in (0,1)$, $k\in (0,k_0)$, and $u(x_n,t_j)$ is the solution of the difference scheme $(\ref{star2})$ defined on $\mathbb{R}_h \times \left( \mathbb{R}_k \cap \left[0,T\right] \right)$ with initial condition $u_0 $ then we have $\DLN{\IHN{x} D^n_+u_j}<C_n$ for $t_j \in \left[0, T\right]$.  Moreover, the constant $C_n$ is independent of the choice of $h$ and $k$. \\

\noindent \textbf{Proof of Lemma 3.2} \hspace{1mm} By lemma 2.7 the statement is true for $0 \leq n \leq 3$ by taking $k_0 = \epsilon$ and $C_n = L$.  We shall prove the statement by induction on $n$.  Assume it is true for all $l\leq n-1$.  We will prove that it is true for $l=n$.  First we shall use the inductive hypothesis to prove some slightly weaker claims which we give below. \\

\noindent \textit{Claim 1} \hspace{1mm} There exists $0<\epsilon_0 \leq \epsilon$ and $C_n>0$ such that if $h\in (0,1)$ and $k\in (0,\epsilon_0)$ and $u(x_n,t_j)$ is the solution of the difference scheme $(\ref{star2})$ defined on $\mathbb{R}_h \times \left( \mathbb{R}_k \cap \left[0,T\right] \right)$ with initial condition $u_0 $ then we have $\DLN{D_+^n u_j} \leq C_n$ for $t_j \in \left[0,T\right]$. \\
\noindent \textit{proof of claim 1} \hspace{1mm} By lemma 3.1 it suffices to prove that there exists a constant $C_n>0$ independent of $h \in (0,1)$, $k \in (0,\epsilon)$, and $t_j\in \left[0,T\right]$ such that if $u(x_n,t_j)$ is the solution of the difference scheme $(\ref{star2})$ defined on $\mathbb{R}_h \times \left( \mathbb{R}_k \cap \left[0,T\right] \right)$ with initial condition $u_0 $ then for $0 \leq t_j,t_{j+1} \leq T$ we have 
\begin{eqnarray*}
 \DLIP{ D^n_+ (Q_j u_{j+1})}{D^n_+ u_{j+1} }  \geq  -C_n \left( \DLNS{ D^n_+ u_{j+1}} + \DLNS{ D^n_+ u_j} + 1 \right) 
\end{eqnarray*}

\noindent To this end we shall fix values for $h\in (0,1)$ and $k\in (0,\epsilon)$ and a solution $u(x_n,t_j)$ of the difference scheme $(\ref{star2})$ defined on $\mathbb{R}_h \times \left( \mathbb{R}_k \cap \left[0,T\right] \right)$ with initial condition $u_0 $.  In order to prove that the above estimate holds for some constant $C_n>0$ we simply prove that the estimate can be made for each term of $Q_j$ and then by adding all these estimates we will obtain the estimate for $Q_j$.  These estimates are analogous to those given in lemma 2.6 however we will also use the inductive hypothesis.  For simplicity we will ignore all occurences of the shift $E$.  Here are the necessary estimates. \\

\noindent \underline{\textit{Estimate for Term}} $\DLIP{D_+^n (u_j^2 D_0 u_{j+1}) }{D_+^n u_{j+1}}$ : \\ 
By the product rule we obtain,
\begin{equation*}
\DLIP{D_+^n (u_j^2 D_0 u_{j+1}) }{D_+^n u_{j+1}} = \sum_{i_1 + i_2 + i_3=n} c_{i_1,i_2,i_3} \DLIP{ (D_+^{i_1}u_j)(D_+^{i_2}u_j)(D_+^{i_3}D_0 u_{j+1}) }{D_+^n u_{j+1}}
\end{equation*}
For the terms where $2 \leq i_1 + i_2 \leq n-2$ we have the estimate,
\begin{eqnarray*}
\DLIP{ (D_+^{i_1}u_j)(D_+^{i_2}u_j)(D_+^{i_3}D_0 u_{j+1}) }{D_+^n u_{j+1}}  \geq  -\SupNorm{D_+^{i_1}u_j} \SupNorm{D_+^{i_2}u_j} \DLN{D_+^{i_3}D_0 u_{j+1}} \DLN{D_+^n u_{j+1}} \\
                                                                            \geq  -C \left( \DLN{u_j} + \DLN{D_+^{n-1}u_j} \right)^2 \left( \DLN{u_{j+1}}+ 
                                                                                   \DLN{D_+^{n-1}u_{j+1}} \right) \DLN{D_+^n u_{j+1}} 
                                                                             \geq -C \left(1+\DLNS{D_+^n u_{j+1}} \right)
\end{eqnarray*}
For the terms where $ i_1 + i_2 = 0$ we have the estimate,
\begin{eqnarray*}
\DLIP{u_j^2 D_+^n D_0 u_{j+1}}{D_+^n u_{j+1}} & \geq & -\frac{1}{2} \SupNorm{D_+ u_j^2} \DLNS{D_+^n u_{j+1}} \geq -C \SupNorm{u_j} \SupNorm{D_+ u_j} \DLNS{D_+^n u_{j+1}} \\
                                              & \geq & -C \left( \DLN{u_j} + \DLN{D_+^{n-1}u_j} \right)^2 \DLNS{D_+^n u_{j+1}} 
                                               \geq  -C \DLNS{D_+^n u_{j+1}} 
\end{eqnarray*}
For the terms where $i_1=1$, $i_2=0$, and $i_3=n-1$ we have the estimate,
\begin{eqnarray*}
\DLIP{u_j D_+ u_j D_+^{n-1}D_0 u_{j+1}}{D_+^n u_{j+1}} & \geq & -\SupNorm{u_j} \SupNorm{D_+ u_j} \DLN{D_+^{n-1}D_0 u_{j+1}} \DLN{D_+^n u_{j+1}} \\
                                                       & \geq & -C \left( \DLN{u_j} + \DLN{D_+^{n-1}u_j} \right)^2 \DLNS{D_+^n u_{j+1}} \geq -C \DLNS{D_+^n u_{j+1}}
\end{eqnarray*}
For the terms where $n-1 \leq i_1 + i_2 \leq n$ and $i_1 \leq i_2$ we have the estimate,
\begin{eqnarray*}
\DLIP{ (D_+^{i_1}u_j)(D_+^{i_2}u_j)(D_+^{i_3}D_0 u_{j+1}) }{D_+^n u_{j+1}}  \geq  -\SupNorm{D_+^{i_3}D_0 u_{j+1}} \SupNorm{D_+^{i_1}u_j} \DLN{D_+^{i_2}u_j} \DLN{D_+^n u_{j+1}} \hspace{1.5cm}\\
                                                                            \geq  -C\left( \DLN{u_{j+1}} + \DLN{D_+^{n-1}u_{j+1}} \right) \left( \DLN{u_j} + \DLN{D_+^n u_j} \right) \DLN{D_+^n u_{j+1}} 
                                                                            \geq  -C \left(1+\DLN{D_+^n u_j} \right) \DLN{D_+^n u_{j+1}}\\
                                                                             =    -C \left(\DLN{D_+^n u_{j+1}} + \DLN{D_+^n u_{j+1}} \DLN{D_+^n u_j} \right) 
                                                                            \geq  -C \left(1+ \DLNS{D_+^n u_{j+1}} + \DLNS{D_+^n u_j} \right) \hspace{3.5cm}
\end{eqnarray*}

\noindent \underline{\textit{Estimate for Term}} $\DLIP{D_+^n (D_+^2 D_- u_{j+1}) }{D_+^n u_{j+1}}$ : \\ 
\begin{equation*}
\DLIP{D_+^n (D_+^2 D_- u_{j+1}) }{D_+^n u_{j+1}} = \DLIP{ D_+^2 D_- (D_+^n u_{j+1}) }{D_+^n u_{j+1}} \geq 0
\end{equation*}

\noindent \underline{\textit{Estimate for Term}} $\DLIP{D_+^n (f_j u_j D_0 u_{j+1}) }{D_+^n u_{j+1}}$ : \\ 
By the product rule we have
\begin{equation*}
\DLIP{D_+^n (f_j u_j D_0 u_{j+1}) }{D_+^n u_{j+1}}  =  \sum_{i_1 + i_2 + i_3 = n} c_{i_1,i_2,i_3} \DLIP{ (D_+^{i_1}f_j)(D_+^{i_2}u_j)(D_+^{i_3}D_0 u_{j+1}) }{D_+^n u_{j+1}}
\end{equation*}
and for the terms where $1 \leq i_1 \leq n$ and $i_3 \leq n-2$ we have the estimate
\begin{eqnarray*}
\DLIP{ (D_+^{i_1}f_j)(D_+^{i_2}u_j)(D_+^{i_3}D_0 u_{j+1})}{D_+^n u_{j+1}} & \geq & -C \SupNorm{D_+^{i_1}f_j} \DLN{D_+^{i_3}D_0 u_{j+1}} \SupNorm{D_+^{i_2}u_j}  \DLN{D_+^n u_{j+1}} \\
                                                                          & \geq & -C \left(\DLN{u_j} + \DLN{D_+^n u_j} \right) \DLN{D_+^n u_{j+1}} \\
                                                                          & \geq & -C \left(1 + \DLN{D_+^n u_j} \right) \DLN{D_+^n u_{j+1}} \\
                                                                          & \geq & -C \left(1+ \DLNS{D_+^n u_{j+1}} + \DLNS{D_+^n u_j} \right) 
\end{eqnarray*}

\noindent and for the term where $i_1 = 1$ and $i_3 = n-1$ we have the estimate
\begin{eqnarray*}
\DLIP{ (D_+ f_j) u_j (D_+^{n-1}D_0 u_{j+1})}{D_+^n u_{j+1}} & \geq & -C \SupNorm{D_+ f_j} \SupNorm{u_j} \DLN{D_+^{n-1}D_0 u_{j+1}} \DLN{D_+^n u_{j+1}} \\
                                                                          & \geq & -C \DLNS{D_+^n u_{j+1}}\\
\end{eqnarray*}

\noindent and for the terms where $i_1 = 0$ and $1 \leq i_2 \leq n-2$ we have the estimate
\begin{eqnarray*}
\DLIP{f_j (D_+^{i_2}u_j)(D_+^{i_3}D_0 u_{j+1}) }{D_+^n u_{j+1}}  \geq  -\SupNorm{f_j D_+^{i_2}u_j} \DLN{D_+^{i_3}D_0 u_{j+1}} \DLN{D_+^n u_{j+1}} \\
                                                                 \geq  -C \SupNorm{\IHN{x} D_+^{i_2}u_j} \left( \DLN{u_{j+1}} + \DLN{D_+^n u_{j+1}} \right) \DLN{D_+^n u_{j+1}} \hspace{6cm}\\
                                                                 \geq  -C \left(\DLN{\IHN{x} u_j} + \DLN{\IHN{x} D_+^{n-1}u_j} \right) \left(1+\DLN{D_+^n u_{j+1}} \right) \DLN{D_+^n u_{j+1}} 
                                                                  \geq -C \left(1+ \DLNS{D_+^n u_{j+1}} \right)\hspace{4mm}
\end{eqnarray*}
\noindent and for the terms where $i_1 = 0$ and $i_2 = 0$ we have the estimate
\begin{eqnarray*}
\DLIP{f_j u_j D_+^n D_0 u_{j+1}}{D_+^n u_{j+1}}  \geq  -\frac{1}{2} \SupNorm{D_+ (f_j u_j)} \DLNS{D_+^n u_{j+1}} \geq -C \left(\DLN{f_j u_j} + \DLN{D_+^2 (f_j u_j)} \right) \DLNS{D_+^n u_{j+1}} \\
                                                 \geq  -C\left( \DLN{\IHN{x} u_j} + \SupNorm{D_+^2 f_j} \DLN{u_j} + \SupNorm{D_+ f_j} \DLN{D_+ u_j}  + \DLN{\IHN{x}D_+^2 u_j} \right) \DLNS{D_+^n u_{j+1}} \geq -C \DLNS{D_+^n u_{j+1}}
\end{eqnarray*}
\noindent and for the terms where $i_1 = 0$ and $n-1 \leq i_2 \leq n$ we have the estimate
\begin{eqnarray*}
\DLIP{f_j (D_+^{i_2}u_j)(D_+^{i_3}D_0 u_{j+1}) }{D_+^n u_{j+1}}  \geq  -\SupNorm{f_j D_+^{i_3}D_0 u_{j+1}} \DLN{D_+^{i_2}u_j} \DLN{D_+^n u_{j+1}} \hspace{5cm}\\
                                                                 \geq  -C \SupNorm{\IHN{x} D_+^{i_3 +1}u_{j+1}} \left( \DLN{u_j} + \DLN{D_+^n u_j} \right) \DLN{D_+^n u_{j+1}}\hspace{27mm} \\
                                                                 \geq  -C\left( \DLN{\IHN{x} u_{j+1}} + \DLN{\IHN{x}D_+^{n-1}u_{j+1}} \right) \left( 1 + \DLN{D_+^n u_j} \right) \DLN{D_+^n u_{j+1}} \hspace{8mm}\\
                                                                 \geq  -C \left( 1 + \DLN{D_+^n u_j} \right) \DLN{D_+^n u_{j+1}} 
                                                                 \geq  -C \left( 1 + \DLNS{D_+^n u_{j+1}} + \DLNS{D_+^n u_j} \right)\hspace{4mm}
\end{eqnarray*}
\noindent \underline{\textit{Estimate for Term}} $\DLIP{D_+^n (f_j (f_x)_j u_{j+1}) }{D_+^n u_{j+1}}$ : \\ 
By the product rule we get
\begin{equation*}
\DLIP{D_+^n (f_j (f_x)_j u_{j+1}) }{D_+^n u_{j+1}}  =  \sum_{i=0}^n c_i \DLIP{D_+^i (f_j (f_x)_j)D_+^{n-i}u_{j+1}}{D_+^n u_{j+1}}
\end{equation*}
and we can bound each term by the estimate,
\begin{eqnarray*}
\DLIP{D_+^i (f_j (f_x)_j)D_+^{n-i}u_{j+1}}{D_+^n u_{j+1}}  & \geq &  -\SupNorm{D_+^i (f_j (f_x)_j)} \DLN{D_+^{n-i}u_{j+1}} \DLN{D_+^n u_{j+1}} \\
                                                           & \geq & -\SupNorm{\frac{d^i}{dx^i}(f_j (f_x)_j)} \left( \DLN{u_{j+1}} + \DLN{D_+^n u_{j+1}} \right) \DLN{D_+^n u_{j+1}} \\
                                                           & \geq &  -C \left( 1+ \DLN{D_+^n u_{j+1}} \right) \DLN{D_+^n u_{j+1}} 
                                                           \geq  -C \left(1+ \DLNS{D_+^n u_{j+1}} \right)
\end{eqnarray*}
\noindent \underline{\textit{Estimate for Term}} $\DLIP{D_+^n (u_j (f_x)_j u_{j+1}) }{D_+^n u_{j+1}}$ : \\ 
By the product rule we get
\begin{equation*}
\DLIP{D_+^n (u_j (f_x)_j u_{j+1})}{D_+^n u_{j+1}}  =  \sum_{i_1 + i_2 + i_3=n} c_{i_1,i_2,i_3} \DLIP{ (D_+^{i_1}u_j)(D_+^{i_2}(f_x)_j)(D_+^{i_3}u_{j+1}) }{D_+^n u_{j+1}}
\end{equation*}
and for the terms where $ i_1 \leq i_3$ we can use the estimate,
\begin{eqnarray*}
\DLIP{ (D_+^{i_1}u_j)(D_+^{i_2}(f_x)_j)(D_+^{i_3}u_{j+1}) }{D_+^n u_{j+1}}  & \geq & -\SupNorm{D_+^{i_2}(f_x)_j} \SupNorm{D_+^{i_1}u_j} \DLN{D_+^{i_3}u_{j+1}} \DLN{D_+^n u_{j+1}} \hspace{1cm}\\
                                                                           & \geq & -C \left( \DLN{u_{j+1}} + \DLN{D_+^n u_{j+1}} \right) \DLN{D_+^n u_{j+1}} \\
                                                                           & \geq & -C \left( 1+\DLN{D_+^n u_{j+1}} \right) \DLN{D_+^n u_{j+1}} 
                                                                            \geq  -C \left( 1+\DLNS{D_+^n u_{j+1}} \right)
\end{eqnarray*}
and for the terms where $i_3 \leq i_1$ we can use the estimate,
\begin{eqnarray*}
\DLIP{ (D_+^{i_1}u_j)(D_+^{i_2}(f_x)_j)(D_+^{i_3}u_{j+1}) }{D_+^n u_{j+1}} & \geq & -\SupNorm{D_+^{i_2}(f_x)_j} \SupNorm{D_+^{i_3}u_{j+1}} \DLN{D_+^{i_1}u_j}  \DLN{D_+^n u_{j+1}} \\
                                                                            &\geq & -C \left( \DLN{u_j} + \DLN{D_+^n u_j} \right) \DLN{D_+^n u_{j+1}} \\
                                                                           & \geq & -C \left( 1+\DLNS{D_+^n u_{j+1}} + \DLNS{D_+^n u_j}\right)
\end{eqnarray*}
\noindent \underline{\textit{Estimate for Term}} $\DLIP{D_+^n (f_j^2 D_0 u_{j+1}) }{D_+^n u_{j+1}}$ : \\ 
By the product rule we get
\begin{equation*}
\DLIP{D_+^n (f_j^2 D_0 u_{j+1}) }{D_+^n u_{j+1}}  =  \sum_{i=0}^n c_i \DLIP{ (D_+^i f_j^2)(D_+^{n-i}D_0 u_{j+1}) }{D_+^n u_{j+1}}
\end{equation*}
and for the terms where $1 \leq i \leq n$ we can use the estimate
\begin{eqnarray*}
\DLIP{ (D_+^i f_j^2)(D_+^{n-i}D_0 u_{j+1}) }{D_+^n u_{j+1}} & \geq & -\SupNorm{D_+^i f_j^2} \DLN{D_+^n u_{j+1}} \DLN{D_+^{n-i}D_0 u_{j+1}} \\
                                                            & \geq & -C \SupNorm{\frac{d^i}{dx^i}f_j^2 } \left( \DLN{u_{j+1}} + \DLN{D_+^n u_{j+1}} \right)^2 \\
                                                            & \geq & -C \left(1+ \DLN{D_+^n u_{j+1}} \right)^2 \geq -C \left(1+ \DLNS{D_+^n u_{j+1}} \right)
\end{eqnarray*}
and for the term $i=0$ we can use the estimate
\begin{eqnarray*}
\DLIP{f_j^2 D_+^n D_0 u_{j+1}}{D_+^n u_{j+1}} \geq -\frac{1}{2} \SupNorm{D_+ f_j^2} \DLNS{D_+^n u_{j+1}} \geq -C \DLNS{D_+^n u_{j+1}}
\end{eqnarray*}
\noindent This concludes the proof of claim 1.  \\

\noindent \textit{Claim 2} \hspace{1mm}   There exists $0<\epsilon_1 \leq \epsilon_0$ and $C_n>0$ such that if $h\in (0,1)$ and $k\in (0,\epsilon_1)$ and $u(x_n,t_j)$ is the solution of the difference scheme $(\ref{star2})$ defined on $\mathbb{R}_h \times \left( \mathbb{R}_k \cap \left[0,T\right] \right)$ with initial condition $u_0 $ then we have $\DLN{D_+^{n+1} u_j} \leq C_n$ for $t_j \in \left[0,T\right]$. \\
\noindent \textit{proof of claim 2} \hspace{1mm} As in the proof of claim 1 we see that by lemma 3.1 it suffices to prove that there exists a constant $C_n>0$ independent of $h \in (0,1)$, $k \in (0, \epsilon_0 )$, and $t_j\in \left[0,T\right]$ such that if $u(x_n,t_j)$ is the solution of the difference scheme $(\ref{star2})$ defined on $\mathbb{R}_h \times \left( \mathbb{R}_k \cap \left[0,T\right] \right)$ with initial condition $u_0 $ then for $0 \leq t_j,t_{j+1} \leq T$ we have 
\begin{eqnarray*}
 \DLIP{ D^{n+1}_+ (Q_j u_{j+1})}{D^{n+1}_+ u_{j+1} }  \geq  -C_n \left( \DLNS{ D^{n+1}_+ u_{j+1}} + \DLNS{ D^{n+1}_+ u_j} + 1 \right) 
\end{eqnarray*}

\noindent To this end we shall again fix values for $h\in (0,1)$ and $k\in (0,\epsilon_0)$ and a solution $u(x_n,t_j)$ to the solution of the difference scheme $(\ref{star2})$ defined on $\mathbb{R}_h \times \left( \mathbb{R}_k \cap \left[0,T\right] \right)$ with initial condition $u_0 $.  In order to prove that the above estimate holds for some constant $C_n>0$ we simply prove that the estimate can be made for each term of $Q_j$ and then by adding all these estimates we will obtain the estimate for $Q_j$.  By using the inductive hypothesis and the result of claim 1, it is easily seen that the estimates shown in the proof of claim 1 can all be applied with $n$ replaced by $n+1$ except for the term $f_j u_j D_0 u_{j+1}$.  We estimate this term below. \\
\indent By the product rule we have,
\begin{equation*}
\DLIP{D_+^{n+1}(f_j u_j D_0 u_{j+1})}{D_+^{n+1}u_{j+1}} = \sum_{i_1 + i_2 + i_3 = n+1} c_{i_1,i_2,i_3} \DLIP{ (D_+^{i_1}f_j)(D_+^{i_2}u_j)(D_+^{i_3}D_0 u_{j+1}) }{D_+^{n+1}u_{j+1}}
\end{equation*}
\noindent The terms corresponding to (a) $1 \leq i_1 \leq n+1$, $i_3 \leq n-1$, (b) $ i_1 =1$, $i_3 = n$ (c) $i_1=0=i_2$, and (d) $i_1=0$ , $n \leq i_2 \leq n+1$ can all be estimated in the same way as those similar terms in the proof of claim 1 with $n$ replaced by $n+1$.  For the terms e) $i_1 = 0$, $1 \leq i_2 \leq n-2$ we may use the same estimate with $n$ replaced by $n+1$ but it doesn't work for $i_2=n-1$.  For this term we have 
\begin{eqnarray*}
\DLIP{f_j D_+^{n-1}u_j D_+^2 D_0 u_{j+1}}{D_+^{n+1}u_{j+1}}  & \geq & -\SupNorm{D_+^{n-1}u_j} \DLN{f_j D_+^2 D_0 u_{j+1}} \DLN{D_+^{n+1}u_{j+1}}  \\
                                                            & \geq & -C \left( \DLN{u_j} + \DLN{D_+^n u_j} \right) \DLN{\IHN{x}D_+^3 u_{j+1}} \DLN{D_+^{n+1}u_{j+1}} \\
                                                           &  \geq &-C \left(1+ \DLNS{D_+^{n+1}u_{j+1}} \right)
\end{eqnarray*}
\noindent This concludes the proof of claim 2.\\

\noindent \textit{Claim 3} \hspace{1mm}  There exists $0<\epsilon_2 \leq \epsilon_1 $ and $C_n>0$ such that if $h\in (0,1)$ and $k\in (0,\epsilon_2)$ and $u(x_n,t_j)$ is the solution of the difference scheme $(\ref{star2})$ defined on $\mathbb{R}_h \times \left( \mathbb{R}_k \cap \left[0,T\right] \right)$ with initial condition $u_0 $ then we have $\DLN{D_+^{n+2} u_j} \leq C_n$ for $t_j \in \left[0,T\right]$. \\
\noindent \textit{proof of claim 3} \hspace{1mm} We again use the same reasoning as given in the above two claims.  In this case it then suffices by lemma 3.1 to prove the estimate.
\begin{eqnarray*}
 \DLIP{ D^{n+2}_+ (Q_j u_{j+1})}{D^{n+2}_+ u_{j+1} }  \geq  -C_n \left( \DLNS{ D^{n+2}_+ u_{j+1}} + \DLNS{ D^{n+2}_+ u_j} + 1 \right) 
\end{eqnarray*}

\noindent As in the proof of claim 2, we may use the previous estimates shown in claim 1 by replacing $n$ by $n+2$ except for the term $f_j u_j D_0 u_{j+1}$.  Again the estimates for the terms (a) $1 \leq i_1 \leq n+2$,$i_3 \leq n$, (b) $i_1 = 1$, $i_3 = n$ (c) $i_1=0=i_2$, and (d) $i_1=0$, $n+1 \leq i_2 \leq n+2$ can all be written in the same way as those similar terms shown in the proof of claim 1 with $n$ replaced by $n+2$.  For the terms (e) $i_1=0$, $1\leq i_2 \leq n-2$ we may use the same estimate with $n+2$ as its corresponding term but it doesn't work for $i_2=n-1$ and $i_2 = n$.  For the terms where $i_1=0$, $i_2=n-1$ we may use the estimate
\begin{eqnarray*}
\DLIP{f_j D_+^{n-1}u_j D_+^3 D_0 u_{j+1}}{D_+^{n+2}u_{j+1}}  &\geq & -\SupNorm{D_+^3 D_0 u_{j+1}} \DLN{f_j D_+^{n-1}u_j} \DLN{D_+^{n+2}u_{j+1}} \\
                                                            & \geq & -C \left( \DLN{u_{j+1}} + \DLN{D_+^5 u_{j+1}} \right) \DLN{\IHN{x}D_+^{n-1}u_j} \DLN{D_+^{n+2}u_{j+1}} \\
                                                          &  \geq &-C \left( 1+\DLNS{D_+^{n+2}u_{j+1}} \right)
\end{eqnarray*}
and for the terms where $i_1=0$, $i_2=n$ we may use the estimate
\begin{eqnarray*}
\DLIP{f_j D_+^n u_j D_+^2 D_0 u_{j+1}}{D_+^{n+2}u_{j+1}} & \geq & -\SupNorm{D_+^n u_j} \DLN{D_+^{n+2}u_{j+1}} \DLN{f_j D_+^2 D_0 u_{j+1}} \\
                                                         & \geq & -C \left(\DLN{u_j} + \DLN{D_+^{n+1}u_j} \right) \DLN{D_+^{n+2}u_{j+1}} \DLN{\IHN{x}D_+^3 u_{j+1}} \\ 
                                                         &\geq &-C \left(1+ \DLNS{D_+^{n+2}u_{j+1}} \right)
\end{eqnarray*}
\noindent This concludes the proof of claim 3.\\

\indent  Now we return to the proof of the lemma in the case $l=n$.  By using the same reasoning as in the above claims we see by lemma 3.1 that in order to construct $k_0 \in (0,\epsilon_2)$ and $C_n$ it suffices to prove the estimate
\begin{eqnarray*}
 \DLIP{ \IHN{x}D^n_+ (Q_j u_{j+1})}{\IHN{x}D^n_+ u_{j+1} }  \geq   -C \left( \DLNS{ \IHN{x}D^n_+ u_{j+1}} + \DLNS{\IHN{x}D^n_+ u_j} + 1 \right) 
\end{eqnarray*} 
\noindent  This estimate will follow by adding all of the below estimates.  In the below estimates we will use our inductive hypothesis and the estimates $\DLN{D_+^{n+i}u_j}\leq C_{n+i}$ for $i=0,1,2$.  For conciseness we shall only formulate estimates for a few terms because the rest can be bounded similarly.

\noindent \underline{\textit{Estimate for Term}} $\DLIP{\IHN{x}^2 D_+^n (u_j^2 D_0 u_{j+1}) }{D_+^n u_{j+1}}$ : \\ 
By the product rule we have
\begin{equation*}
\DLIP{\IHN{x}^2 D_+^n (u_j^2 D_0 u_{j+1}) }{D_+^n u_{j+1}}  =  \sum_{i_1 + i_2 + i_3 = n} c_{i_1,i_2,i_3} \DLIP{\IHN{x}^2 (D_+^{i_1}u_j)(D_+^{i_2}u_j)(D_+^{i_3}D_0 u_{j+1}) }{D_+^n u_{j+1}}
\end{equation*}
For the terms where $i_3 = n-1$ we have
\begin{eqnarray*}
\DLIP{\IHN{x}^2 u_j D_+ u_j(D_+^{n-1}D_0 u_{j+1}) }{D_+^n u_{j+1}} & \geq & -\SupNorm{u_j}\SupNorm{D_+ u_j} \DLN{\IHN{x}D_+^n u_{j+1}} \DLN{\IHN{x}D_+^{n-1}u_{j+1}} \\
                                                                                    & \geq &  -C \left(\DLN{u_j}+ \DLN{D_+^{n-1}u_j} \right)^2 \DLN{\IHN{x}D_+^n u_{j+1}} \\
                                                                                    & \geq &  -C \DLNS{\IHN{x}D_+^n u_{j+1}} \hspace{1cm}
\end{eqnarray*}
\noindent For the terms where $i_3 = n$ we have
\begin{eqnarray*}
\DLIP{\IHN{x}^2 u_j^2 (D_+^n D_0 u_{j+1}) }{D_+^n u_{j+1}} & \geq & -\frac{1}{2}\SupNorm{\frac{1}{\IHN{x}^2}D_+(\IHN{x}^2 u_j^2) } \DLNS{\IHN{x}D_+^n u_{j+1}} \\
                                                           & \geq & -C \DLNS{\IHN{x}D_+^n u_{j+1}} \left( \SupNorm{D_+ u_j^2} + \SupNorm{u_j^2} \SupNorm{\frac{D_+ \IHN{x}^2}{\IHN{x}^2}} \right)  \\
                                                           & \geq & -C \DLNS{\IHN{x}D_+^n u_{j+1}} \left( \SupNorm{D_+ u_j}\SupNorm{u_j} + \SupNorm{u_j}^2 \cdot 1 \right)  \\
                                                           & \geq & -C \DLNS{\IHN{x}D_+^n u_{j+1}} \left(\DLN{u_j} + \DLN{D_+^n u_j} \right)^2
                                                            \geq  -C \DLNS{\IHN{x}D_+^n u_{j+1}}
\end{eqnarray*}
\noindent For the terms where $1 \leq i_3 \leq n-2$ and $i_1 \leq i_2$ we have
\begin{eqnarray*}
\DLIP{\IHN{x}^2 (D_+^{i_1}u_j)(D_+^{i_2}u_j)(D_+^{i_3}D_0 u_{j+1}) }{D_+^n u_{j+1}} & \geq & -\SupNorm{D_+^{i_1} u_j} \SupNorm{D_+^{i_2}u_j} \DLN{\IHN{x}D_+^n u_{j+1}} \DLN{\IHN{x}D_+^{i_3}D_0 u_{j+1}} \\
                                                                                    & \geq & -C \DLN{\IHN{x}D_+^n u_{j+1}}\DLN{\IHN{x} D_+^{i_3 +1}u_{j+1}} \\
                                                                                     & \geq &-C\left(1+ \DLNS{\IHN{x}D_+^n u_{j+1}} \right) 
\end{eqnarray*}
\noindent For the terms where $i_3 = 0$ and $i_1 \leq i_2 \leq n-1$ we have
\begin{eqnarray*}
\DLIP{\IHN{x}^2 (D_+^{i_1}u_j)(D_+^{i_2}u_j)(D_0 u_{j+1}) }{D_+^n u_{j+1}} & \geq & -\SupNorm{D_0 u_{j+1}} \SupNorm{D_+^{i_1}u_j} \DLN{\IHN{x}D_+^n u_{j+1}} \DLN{\IHN{x}D_+^{i_2}u_j} \\
                                                                           & \geq & -C \DLN{\IHN{x}D_+^n u_{j+1}} \geq -C \left(1+ \DLNS{\IHN{x}D_+^n u_{j+1}} \right)
\end{eqnarray*}
\noindent For the terms where $i_2 = n$ we have
\begin{eqnarray*}
\DLIP{\IHN{x}^2 u_j(D_+^n u_j)(D_0 u_{j+1}) }{D_+^n u_{j+1}} & \geq & -\SupNorm{u_j} \SupNorm{D_0 u_{j+1}} \DLN{\IHN{x}D_+^n u_{j+1}} \DLN{\IHN{x}D_+^n u_j} \\
                                                             & \geq & -C \DLN{\IHN{x}D_+^n u_{j+1}} \DLN{\IHN{x}D_+^n u_j} \\
                                                             & \geq & -C \left( \DLNS{\IHN{x}D_+^n u_{j+1}} + \DLNS{\IHN{x}D_+^n u_j} \right)
\end{eqnarray*}
\noindent  This concludes the proof of the lemma. \hspace{2cm} $\square$ \\

\noindent \textbf{Lemma 3.3} \hspace{1mm} Let $u_0\in S^{-\infty}(\mathbb{R})$ and $K,T,\epsilon>0$ all be given as in lemma 2.7 and let $N,n \in \mathbb{N}$.  Then there exists $0 < k_0 \leq \epsilon$ depending on $N$ and $n$ and there exists $C_{N,n} >0$ such that if $h \in (0,1)$, $k\in (0,k_0)$, and $u(x_n,t_j)$ is the solution of the difference scheme $(\ref{star2})$ defined on $\mathbb{R}_h \times \left( \mathbb{R}_k \cap \left[0,T\right] \right)$ with initial condition $u_0 $ then we have $\DLN{\IHN{x}^N D^n_+u_j}<C_{N,n}$ for $t_j \in \left[0, T\right]$.  Moreover, the constant $C_{N,n}$ is independent of the choice of $h$ and $k$. \\

\noindent \textbf{Proof of Lemma 3.3} \hspace{1mm} By lemma 3.2 the assertion is true for $N=1$ and for all $n\in\mathbb{N}$, therefore it is clearly true for $N\in \left\{0,\frac{1}{2},1\right\}$ and for all $n\in\mathbb{N}$.  We prove the assertion by induction on $N$ with increment $\frac{1}{2}$.  Assume it is true for all $M\leq N-\frac{1}{2}$ and for all $n\in\mathbb{N}$.  We will show that for $M=N$ the statement is satisfied by all $n\in \mathbb{N}$ by induction on $n$.  \\
\indent  Let us denote the value of $k_0$ corresponding to a particular value of $N$ and $n$ by $k_0(N,n)$.  By construction of $k_0(N,n)$ for its known values we see that if $N_1 < N_2$ then $k_0(N_1,n)\leq k_0(N_2,n)$ and if $n_1 < n_2$ then $k_0(N,n_1)\leq k_0(N,n_2)$.  Let $\epsilon_0 := \min \left\{ k_0(1,2),k_0(N-\frac{1}{2},0)\right\}$.  By use of lemma 3.1 our statement for $n=0$ will follow if we prove that there exists a constant $C_{N,0}>0$ independent of $h \in (0,1)$, $k \in (0,\epsilon_0)$, and $t_j\in \left[0,T\right]$ such that if $u(x_n,t_j)$ is the solution of the difference scheme $(\ref{star2})$ defined on $\mathbb{R}_h \times \left( \mathbb{R}_k \cap \left[0,T\right] \right)$ with initial condition $u_0 $ then for $0 \leq t_j,t_{j+1} \leq T$ we have

\begin{eqnarray*}
\DLIP{ \IHN{x}^{2N}  Q_j u_{j+1}}{ u_{j+1} } & \geq &  -C_{N,0} \left( \DLNS{ \IHN{x}^N u_{j+1}} + \DLNS{\IHN{x}^N u_j} + 1 \right)
\end{eqnarray*}

\noindent To this end we shall fix values for $h\in (0,1)$ and $k\in (0,\epsilon_0)$ and a solution $u(x_n,t_j)$ of the difference scheme $(\ref{star2})$ defined on $\mathbb{R}_h \times \left( \mathbb{R}_k \cap \left[0,T\right] \right)$ with initial condition $u_0 $.  In order to prove that the above estimate holds for some constant $C_{N,0}>0$ we simply prove that the estimate can be made for each term of $Q_j$ and then by adding all these estimates we will obtain the estimate for $Q_j$.  Here are the necessary estimates.  As in the proof of lemma 3.2 we will ignore all occurences of the shift $E$ for simplicity.\\

\noindent \underline{\textit{Estimate for Term}} $\DLIP{\IHN{x}^{2N} u_j^2 D_0 u_{j+1} }{u_{j+1}}$ : \\ 
\begin{eqnarray*}
\DLIP{\IHN{x}^{2N} u_j^2 D_0 u_{j+1} }{u_{j+1}} & \geq & -\SupNorm{u_{j+1}} \SupNorm{D_0 u_{j+1}} \DLNS{\IHN{x}^N u_j} \geq -C \DLNS{\IHN{x}^N u_j}
\end{eqnarray*}

\noindent \underline{\textit{Estimate for Term}} $\DLIP{\IHN{x}^{2N} D_+^2 D_- u_{j+1} }{u_{j+1}}$ : \\ 
By proposition 2.1 part 4 we have
\begin{eqnarray*}
\DLIP{\IHN{x}^{2N} D_+^2 D_- u_{j+1} }{u_{j+1}} & \geq & -\frac{1}{2}\DLIP{u_{j+1}D_-^2D_+ \IHN{x}^{2N}}{u_{j+1}} + \DLIP{D_0 \IHN{x}^{2N}D_+ u_{j+1}}{D_+ u_{j+1}}\\
                                                &      &  + \frac{1}{2}\DLIP{D_-\IHN{x}^{2N} D_+ u_{j+1}}{D_+ u_{j+1}} \\
                                                & \geq & -C \DLNS{\IHN{x}^{N-\frac{3}{2}}u_{j+1}} - C\DLNS{\IHN{x}^{N-\frac{1}{2}}D_+ u_{j+1}} \geq -C
\end{eqnarray*}
\noindent \underline{\textit{Estimate for Term}} $\DLIP{\IHN{x}^{2N} f_j u_j D_0 u_{j+1}}{u_{j+1}}$ : \\ 
\begin{eqnarray*}
\DLIP{\IHN{x}^{2N} f_j u_j D_0 u_{j+1}}{u_{j+1}} & \geq & -\SupNorm{f_j D_0 u_{j+1}} \DLN{\IHN{x}^N u_{j+1}} \DLN{\IHN{x}^N u_j} \\
                                                 & \geq & -C \SupNorm{\IHN{x} D_+ u_{j+1}}\DLN{\IHN{x}^N u_{j+1}} \DLN{\IHN{x}^N u_j} \\ 
                                                 & \geq & -C \DLN{\IHN{x}^N u_{j+1}} \DLN{\IHN{x}^N u_j}\geq -C \left(\DLNS{\IHN{x}^N u_{j+1}} + \DLNS{\IHN{x}^N u_j} \right) 
\end{eqnarray*}
\noindent \underline{\textit{Estimate for Term}} $\DLIP{\IHN{x}^{2N} f_j (f_x)_j u_{j+1}}{u_{j+1}}$ : \\ 
\begin{eqnarray*}
\DLIP{\IHN{x}^{2N} f_j (f_x)_j u_{j+1}}{u_{j+1}} \geq -\SupNorm{f_j (f_x)_j} \DLNS{\IHN{x}^N u_{j+1}} \geq -C \DLNS{\IHN{x}^N u_{j+1}}
\end{eqnarray*}
\noindent \underline{\textit{Estimate for Term}} $\DLIP{\IHN{x}^{2N} u_j (f_x)_j u_{j+1}}{u_{j+1}}$ : \\ 
\begin{eqnarray*}
\DLIP{\IHN{x}^{2N} u_j (f_x)_j u_{j+1}}{u_{j+1}} \geq -\SupNorm{(f_x)_j} \SupNorm{u_j} \DLNS{\IHN{x}^N u_{j+1}} \geq -C \DLNS{\IHN{x}^N u_{j+1}}
\end{eqnarray*}
\noindent \underline{\textit{Estimate for Term}} $\DLIP{\IHN{x}^{2N} f_j^2 D_0 u_{j+1}}{u_{j+1}}$ : \\ 
By proposition 2.1 part 3 we have
\begin{eqnarray*}
\DLIP{\IHN{x}^{2N} f_j^2 D_0 u_{j+1}}{u_{j+1}} \geq -\frac{1}{2}\DLIP{D_+ (\IHN{x}^{2N}f_j^2 ) E u_{j+1}}{u_{j+1}} \geq -C \DLNS{\IHN{x}^N u_{j+1}}
\end{eqnarray*}

\noindent This concludes the proof for the case $n=0$.  \\

\indent Now we will assume that the statement is true for $l\leq n-1$ and we will prove that it is true for $l=n$.  Let $\epsilon_1 := \min \left\{ k_0(1,n),k_0(1,2),k_0(N-\frac{1}{2}),k_0(N,n-1)\right\}$.  The proof will again follow from lemma 3.1 if we can prove that there exists a constant $C_{N,n}>0$ independent of $h \in (0,1)$, $k \in (0,\epsilon_1)$, and $t_j\in \left[0,T\right]$ such that if $u(x_n,t_j)$ is the solution of the difference scheme $(\ref{star2})$ defined on $\mathbb{R}_h \times \left( \mathbb{R}_k \cap \left[0,T\right] \right)$ with initial condition $u_0 $ then for $0 \leq t_j,t_{j+1} \leq T$ we have

\begin{eqnarray*}
 \DLIP{ \IHN{x}^{2N} D^n_+ (Q_j u_{j+1})}{D^n_+ u_{j+1} }  \geq   -C_{N,n} \left( \DLNS{ \IHN{x}^N D^n_+ u_{j+1}} + \DLNS{\IHN{x}^N D^n_+ u_j} + 1 \right) 
\end{eqnarray*}
\noindent To this end we shall fix values for $h\in (0,1)$ and $k\in (0,\epsilon_1)$ and a solution $u(x_n,t_j)$ of the difference scheme $(\ref{star2})$ defined on $\mathbb{R}_h \times \left( \mathbb{R}_k \cap \left[0,T\right] \right)$ with initial condition $u_0 $.  In order to prove that the above estimate holds for some constant $C_{N,n}>0$ we again prove that the estimate can be made for each term of $Q_j$ and then by adding all these estimates we will obtain the estimate for $Q_j$.  Here are the necessary estimates.  For conciseness we shall only formulate a few of these estimates because the rest can be carried out similarly. \\

\noindent \underline{\textit{Estimate for Term}} $\DLIP{\IHN{x}^{2N} D_+^n (u_j^2 D_0 u_{j+1}) }{D_+^n u_{j+1}}$ : \\ 
By the product rule we have
\begin{equation*}
\DLIP{\IHN{x}^{2N} D_+^n (u_j^2 D_0 u_{j+1}) }{D_+^n u_{j+1}} = \sum_{i_1 + i_2 + i_3 = n} c_{i_1,i_2,i_3} \DLIP{\IHN{x}^{2N} (D_+^{i_1}u_j)(D_+^{i_2}u_j)(D_+^{i_3}D_0 u_{j+1}) }{D_+^n u_{j+1}}
\end{equation*}
\noindent For the terms where $1 \leq i_1 + i_2 \leq n-2$ we have,
\begin{eqnarray*}
\DLIP{\IHN{x}^{2N} (D_+^{i_1}u_j)(D_+^{i_2}u_j)(D_+^{i_3}D_0 u_{j+1}) }{D_+^n u_{j+1}} & \geq & -\SupNorm{D_+^{i_1}u_j} \SupNorm{D_+^{i_2}u_j} \DLN{\IHN{x}^N D_+^n u_{j+1}} \DLN{\IHN{x}^N D_+^{i_3}D_0 u_{j+1}} \\
                                                                                       & \geq & -C \DLN{\IHN{x}^N D_+^n u_{j+1}} \DLN{\IHN{x}^N D_+^{i_3 +1} u_{j+1}} \\
                                                                                       & \geq & -C \left(\DLN{\IHN{x}^N u_{j+1}} + \DLN{\IHN{x}^N D_+^n u_{j+1}} \right)^2 \\
                                                                                       & \geq & -C \left(1+\DLN{\IHN{x}^N D_+^n u_{j+1}} \right)^2 \geq -C \left(1+ \DLNS{\IHN{x}^N D_+^n u_{j+1}} \right)
\end{eqnarray*}
\noindent For the terms where $i_1 + i_2 = 0$ we have,
\begin{eqnarray*}
\DLIP{\IHN{x}^{2N} u_j^2(D_+^nD_0 u_{j+1}) }{D_+^n u_{j+1}} & \geq & -\frac{1}{2}\SupNorm{\frac{1}{\IHN{x}^{2N}}D_+(\IHN{x}^{2N}u_j^2)} \DLNS{\IHN{x}^N D_+^n u_{j+1}} \\
                                                            & \geq & -C \left( \SupNorm{\IHN{x}^{-1}}\SupNorm{ u_j}^2 + \SupNorm{u_j}\SupNorm{D_+ u_j} \right) \DLNS{\IHN{x}^N D_+^n u_{j+1}} \\
                                                            & \geq & -C \DLNS{\IHN{x}^N D_+^n u_{j+1}}  
\end{eqnarray*}
\noindent For the terms where $n-1 \leq i_1 + i_2 \leq n$ and $i_1 \leq i_2$ we have,
\begin{eqnarray*}
\DLIP{\IHN{x}^{2N} (D_+^{i_1}u_j)(D_+^{i_2}u_j)(D_+^{i_3}D_0 u_{j+1}) }{D_+^n u_{j+1}} & \geq & -\SupNorm{D_+^{i_1}u_j} \DLN{D_+^{i_3}D_0 u_{j+1}} \DLN{\IHN{x}^N D_+^n u_{j+1}} \DLN{\IHN{x}^N D_+^{i_2}u_j} \\
                                                                                       & \geq & -C \DLN{\IHN{x}^N D_+^n u_{j+1}}\DLN{\IHN{x}^N D_+^{i_2} u_j} \\
                                                                                       & \geq & -C \DLN{\IHN{x}^N D_+^n u_{j+1}} \left( \DLN{\IHN{x}^N u_j}+\DLN{\IHN{x}^N D_+^n u_j}\right) \\
                                                                                       & \geq & -C \DLN{\IHN{x}^N D_+^n u_{j+1}} \left(1+ \DLN{\IHN{x}^N D_+^n u_j}\right)\\
                                                                                       & \geq & -C \left(1+ \DLNS{\IHN{x}^N D_+^n u_{j+1}} + \DLNS{\IHN{x}^N D_+^n u_j} \right)
\end{eqnarray*}

\noindent \underline{\textit{Estimate for Term}} $\DLIP{\IHN{x}^{2N} D_+^n (D_+^2 D_- u_{j+1}) }{D_+^n u_{j+1}}$ : \\ 
By proposition 2.1 part 4 we have
\begin{eqnarray*}
\DLIP{\IHN{x}^{2N} D_+^n (D_+^2 D_- u_{j+1}) }{D_+^n u_{j+1}} & \geq & -\frac{1}{2}\DLIP{D_+^n u_{j+1}}{ (D_+^n u_{j+1})(D_-^2D_+ \IHN{x}^{2N})} +\\ & & + \DLIP{ (D_0 \IHN{x}^{2N})(D_+^{n+1}u_{j+1}) }{D_+^{n+1}u_{j+1}} + \\ & & + \frac{1}{2}\DLIP{ (D_- \IHN{x}^{2N})(D_+^{n+1}u_{j+1}) }{D_+^{n+1}u_{j+1}} \\
                                                              & \geq &-\frac{1}{2} \DLN{\IHN{x}^{N-2}D_+^n u_{j+1}} \DLN{\IHN{x}^{N-1}D_+^n u_{j+1}} - \DLNS{\IHN{x}^{N-\frac{1}{2}}D_+^n u_{j+1}} \\
                                                              & \geq & -C \DLNS{\IHN{x}^{N-\frac{1}{2}}D_+^n u_{j+1}} \geq -C \cdot 1
\end{eqnarray*}

\noindent This concludes the proof of the lemma. $\square$ \\

\noindent \textbf{Corollary 3.4} \hspace{1mm} Let $u_0\in S^{-\infty}(\mathbb{R})$ and $K,T,\epsilon>0$ all be given as in lemma 2.7 and let $N,n \in \mathbb{N}$.  Then there exists $0 < k_0 \leq \epsilon$ depending on $N$ and $n$ and there exists $C_{N,n} > 0$ such that if $h \in (0,1)$, $k\in (0,k_0)$, and $u(x_n,t_j)$ is the solution of the difference scheme $(\ref{star2})$ defined on $\mathbb{R}_h \times \left( \mathbb{R}_k \cap \left[0,T\right] \right)$ with initial condition $u_0 $ then we have $\DLN{D^n_+ \left(x^N u_j\right)}<C_{N,n}$ for $t_j \in \left[0, T\right]$.  Moreover, the constant $C_{N,n}$ is independent of the choice of $h$ and $k$. \\

\noindent \textbf{Proof of Corollary 3.4} \hspace{1mm} By the product rule for $D_+$ this statement can be proven by induction on $n$.  We shall omit the necessary details here. $\square$

\subsection{Boundedness of Time-Differentiated Extended Discrete Solutions}

\noindent In this section we will show that a certain time-extension of the discrete solution with domain $\mathbb{R}_h \times \mathbb{R}_k$ remains bounded in the Schwartz semi-norms $\DLN{D_{t,+}^m D_+^n (x^N \cdot)}$ for $n,N \in \mathbb{N}$ and $m\in \left\{0,1,2,3\right\}$. \\
\indent Suppose $u(x_n,t_j)$ is mesh function defined on $\mathbb{R}_h \times \left( \mathbb{R}_k \cap \left[0,T\right] \right)$ for some $h,k,T>0$ where $h\in (0,1)$ and $k \leq T/3$.  We shall define an extension of $u$ to $\mathbb{R}_h \times \mathbb{R}_k$ by the following: \\
Let $\phi(t) \in C^{\infty}_c\left(\mathbb{R}\right)$ such that $\phi(t) = 1$ for $t \in \left[-1,T+1\right]$ and $\phi(t) = 0$ for $t \notin \left[-2,T+2\right]$.  \\
For $t_j > T$ we define recursively 

$$u(x_n,t_j) :=  u(x_n,t_{j-1}) + k D_{t,+} u(x_n,t_{j-2}) + k^2 D_{t,+}^2 u(x_n,t_{j-3}) $$ 

\noindent and similarly for $t_j <0$ we define 

$$u(x_n,t_j) :=  u(x_n,t_{j+1}) - k D_{t,+} u(x_n,t_{j+1})+k^2 D_{t,+}^2 u(x_n,t_{j+1})  $$ 

\noindent Then we define $\hat{u}:= \phi \cdot u $, which is clearly an extension of $u$ to $\mathbb{R}_h \times \mathbb{R}_k$ and is compactly supported in time.  We will now show that a finite time discrete solution extended in this way remains bounded in the discrete Schwartz semi-norms 
$\DLN{D_{t,+}^m D^n_+ \left( x^N u \right)}$ for $m=0,1,2,3$.   \\  

\noindent \textbf{Lemma 3.5} \hspace{1mm} Let $u_0\in S^{-\infty}(\mathbb{R})$ and $K,T,\epsilon >0$ all be given as in lemma 2.7 and let $N,n \in \mathbb{N}$, $m \in \left\{0,1,2,3\right\}$.  Then there exists $0 < k_0 \leq \epsilon$ depending on $m,n$ and $N$ and there exists $C_{m,N,n} > 0$ such that if $h \in (0,1)$, $k\in (0,k_0)$, and $u(x_n,t_j)$ is the solution of the difference scheme $(\ref{star2})$ defined on $\mathbb{R}_h \times \left( \mathbb{R}_k \cap \left[0,T\right] \right)$ with initial condition $u_0 $ then we have 
\begin{equation}
\DLN{D_{t,+}^m D^n_+ \left( x^N \hat{u}(x, t_j) \right)}\leq C  \label{L9.1}
\end{equation}  
\noindent for each $t_j \in \mathbb{R}_k$.\\

\noindent \textbf{Proof of Lemma 3.5} \hspace{1mm} Let $n,N\in \mathbb{N}$ be fixed values.  First we will assume that $m=0$ and prove the lemma in this case. \\
\indent By corollary 3.4 there exists such $0<k_1(n,N)\leq \epsilon$ and $C>0$ such that (\ref{L9.1}) holds if $k\in (0,k_1)$ and $t_j \in \left[0,T\right]$.  If $t_j \notin \left[-2,T+2\right]$ then clearly (\ref{L9.1}) is still true because $\hat{u}(\cdot,t_j) = 0$, so we only need to prove that (\ref{T.1}) holds for $t_j \in \left[-2,0 \right]$ and for $t_j \in \left[T,T+2 \right]$.  We may assume that $k_1$ is sufficiently small so that $t_2 \in \left[0,T\right]$ because otherwise we may simply decrease the value of $k_1$.  \\
\indent If $t_j, t_{j+1},t_{j+2}\in \left[0,T\right]$ then $(\ref{L9.1})$ also holds with $m=1$ and $m=2$ for $k\in (0,k_2)$ for some $0<k_2 \leq k_1$ because $\hat{u}$ is a solution to the difference scheme $(\ref{star2})$ and so we could replace $D_{t,+} \hat{u}(\cdot, t_j)$ by the rest of $(\ref{star2})$ and use some elementary inequalities together with corollary 3.4 to deduce $(\ref{L9.1})$.  Using this fact we will prove that $(\ref{L9.1})$ holds for $m=0$ where $k\in (0,k_2)$ and with $t_j \in \left[-2, 0 \right]$ and $t_j \in \left[T, T+2 \right]$.  \\
\indent If $t_j \in \left[-2, 0\right]$ then by definition of $u(\cdot,t_j)$ we have
\begin{eqnarray*}
u(\cdot,t_j)& = &  u(\cdot,t_{j+1}) - k D_{t,+}u(\cdot,t_{j+1})+k^2 D_{t,+}^2u(\cdot,t_{j+1})   \\
                & = & u(\cdot,t_{j+l}) - k \sum_{p=1}^l D_{t,+}u(\cdot,t_{j+p})+k^2\sum_{p=1}^l D_{t,+}^2u(\cdot,t_{j+p})   \\
                & = & u(\cdot,t_{j+l}) - k \sum_{p=2}^l D_{t,+}u(\cdot,t_{j+p})+k^2 l D_{t,+}^2u(\cdot,t_{j+l})-k D_{t,+}u(\cdot,t_{j+1})  \\
                & = & u(\cdot,t_{j+l}) - k \sum_{p=2}^l D_{t,+}u(\cdot,t_{j+p})+k^2 l D_{t,+}^2u(\cdot,t_{j+l})-k D_{t,+}u(\cdot,t_{j+2}) + k^2 D_{t,+}^2u(\cdot,t_{j+2})  \\
                & = & u(\cdot,t_{j+l}) - k \sum_{p=3}^l D_{t,+}u(\cdot,t_{j+p})+k^2 (l+1) D_{t,+}^2u(\cdot,t_{j+l})-2 k D_{t,+}u(\cdot,t_{j+2}) \\
                & = & u(\cdot,t_{j+l}) - k \sum_{p=4}^l D_{t,+}u(\cdot,t_{j+p})+k^2 (l+1+2) D_{t,+}^2u(\cdot,t_{j+l})-3 k D_{t,+}u(\cdot,t_{j+3}) \\
                & = & u(\cdot,t_{j+l}) - k l D_{t,+}u(\cdot,t_{j+l})+k^2 \frac{l(l+1)}{2} D_{t,+}^2 u(\cdot,t_{j+l}) 
\end{eqnarray*}

\noindent  Therefore by using the above equality we have 

\begin{eqnarray*}
\left| D^n_+ \left( x^N \hat{u}(x,t_j)\right) \right| & = & \left| D^n_+ \left( x^N \phi(t_j) u(x,t_j) \right) \right|  \leq \SupNorm{\phi} \left| D^n_+ \left( x^N u(x,t_j)\right) \right| \\
                                & \leq & C \left( \left| D^n_+ \left( x^N u(x,t_{j+l})\right) \right| +  k l \left| D^n_+  D_{t,+}\left( x^N u(x,t_{j+l})\right) \right| \right. \\
                                & & \left. + k^2 \frac{l(l+1)}{2} \left| D^n_+ D_{t,+}^2 \left( x^N u(x,t_{j+l})\right) \right| \right) \\ 
                                & \leq & C \left( \left|D^n_+ \left( x^N u(x,t_{0})\right) \right| +  k \frac{2}{k} \left|D^n_+ D_{t,+} \left( x^N u(x,t_{0})\right) \right|+k^2 \frac{3}{k^2} \left| D^n_+ D_{t,+}^2 \left( x^N u(x,t_{0})\right) \right| \right) \\ 
                                & \leq & C \left( \left| D^n_+ \left( x^N u(x,0) \right) \right| + \left| D_{t,+} D^n_+ \left( x^N u(x,0)\right) \right|+\left| D_{t,+}^2 D^n_+ \left( x^N u(x,0)\right) \right| \right) 
\end{eqnarray*}

\noindent  By taking $L^2_h$ norm of the above inequality and using the fact that $(\ref{L9.1})$ holds for $t_0 \in \left[ 0,T \right]$ with $m=0,1,2$ we deduce that $(\ref{L9.1})$ holds for any $t_j \in \left[-2,0\right]$.  In a similar way we can prove that $(\ref{L9.1})$ holds for $t_j  \in \left[T,T+2\right]$. \\
\indent If $t_j \in \left[T, T+2\right]$ then by definition of $u(\cdot,t_j)$ we can show, in a similar way that we did for the case $t_j \in \left[-2,0\right]$ that

\begin{eqnarray*}
u(\cdot,t_j) & = &  u(x_n,t_{j-1}) + k D_{t,+}u(x_n,t_{j-2}) + k^2 D_{t,+}^2u(x_n,t_{j-3}) \\
                & = & u(\cdot,t_{j-l}) + k l D_{t,+}u(\cdot,t_{j-l-1})+k^2 \frac{l(l+1)}{2} D_{t,+}^2 u(\cdot,t_{j-l-2}) 
\end{eqnarray*}

\noindent Therefore we have, by similar estimates from the case $t_j \in \left[-2,0\right]$
\begin{eqnarray*}
\left| D^n_+ \left( x^N \hat{u}(x,t_j) \right) \right| & \leq & C \left( \left| D^n_+ \left( x^N u(x,\tilde{T}) \right) \right| + \left| D_{t,+}D^n_+  \left( x^N u(x,\tilde{T}-k) \right) \right|+\left| D_{t,+}^2 D^n_+ \left( x^N u(x,\tilde{T}-2k)\right) \right| \right) 
\end{eqnarray*}

\noindent where $\tilde{T}$ is understood as the largest $t_j\in \left[0,T\right]$.  By taking $L^2_h$ norm of the above inequality and using the fact that $(\ref{L9.1})$ holds for $t_j = \tilde{T},\tilde{T}-k, \tilde{T}-2k \in \left[0,T\right]$ with $m=0,1,2$ respectively, we deduce that $(\ref{L9.1})$ holds for any $t_j \in \left[T,T+2\right]$.  Therefore, we have shown that $(\ref{L9.1})$ holds for any $t_j \in \mathbb{R}_k$ and for $m=0$ with $k \in (0,k_2)$ thus we may take $k_0 = k_2$ in this case.  \\
\indent  Now we shall assume that $m=1$ and prove the lemma in this case for $k \in (0,k_2)$.  We have already established the $m=1$ lemma inside the interval $\left[0,T\right]$.  \\
\indent For $t_j\notin \left[ -3,T+2 \right]$ it is clear that $(\ref{L9.1})$ holds because $D_{t,+} u(\cdot,t_j) = 0$.  If $t_j \in \left[-3, 0 \right]$ then by the product rule for $D_{t,+}$ we have,

\begin{eqnarray*}
\left| D_{t,+} D^n_+ \left( x^N \hat{u}(x,t_j)\right) \right| &  =  & \left| D^n_+ D_{t,+}\left( x^N \phi(t_j) u(x,t_j) \right) \right| \\
                                                         & \leq & C \left( \left|  D^n_+ \left( x^N u(x,t_j) \right) \right| + \left| D^n_+ D_{t,+} \left( x^N u(x,t_j)\right) \right| \right) \\
                                                         & \leq & C \left( \left| D^n_+ \left( x^N u(x,t_j) \right) \right| + \left|  D^n_+ D_{t,+} \left( x^N u(x,t_{j+1}) \right) \right|+ k \left|  D^n_+ D_{t,+}^2 \left( x^N u(x,t_{j+1}) \right) \right| \right) \\
                                                         & \leq & C \left( \left|D^n_+  \left( x^N u(x,t_j) \right) \right| + \left| D^n_+ D_{t,+}\left( x^N u(x,t_{j+l})\right) \right|+ k l \left| D^n_+ D_{t,+}^2\left( x^N u(x,t_{j+l})\right) \right| \right) \\
                                                         & \leq & C \left( \left|D^n_+  \left( x^N u(x,t_j)\right) \right| + \left| D^n_+ D_{t,+}\left( x^N u(x,t_0) \right) \right|+ k \frac{3}{k}\left|  D^n_+ D_{t,+}^2 \left( x^N u(x,t_0) \right) \right| \right) \\
                                                         & \leq & C \left( \left| D^n_+  \left( x^N u(x,t_j) \right) \right| + \left| D_{t,+} D^n_+ \left( x^N u(x,0) \right) \right|+ \left| D_{t,+}^2 D^n_+ \left( x^N u(x,0)\right) \right| \right) 
\end{eqnarray*}

\noindent  If we take the $L^2_h$ norm of both sides of the above inequality and use the fact that $(\ref{L9.1})$ holds with $m=0$ for any $t_j\in \mathbb{R}_k$ and it holds with $m=1,2$ for $t_0 \in \left[0,T\right]$ we then deduce that $(\ref{L9.1})$ holds with $m=1$ for $t_j \in \left[-3,0\right]$.  Similarly for $t_j \in \left[T,T+2\right]$ we may use the product rule for $D_{t,+}$ to write similar estimates and we will obtain,

\begin{eqnarray*}
\left| D_{t,+}  D^n_+ \left( x^N \hat{u}(x,t_j) \right) \right| & = & C \left( \left| D^n_+ \left(x^N u(x,t_j) \right) \right| + \left| D^n_+ D_{t,+} \left( x^N u(x,t_{j-1}) \right) \right|  \right. \\
                                                         & & \left. + k \left| D^n_+ D_{t,+}^2 \left( x^N u(x,\tilde{T}-2k) \right) \right| \right) \\
                                                         & \leq & C \left( \left| D^n_+ \left( x^N u(x,t_j) \right) \right| + \left| D^n_+ D_{t,+} \left( x^N u(x,\tilde{T}-k) \right) \right| \right. \\
                                                         & & \left. + k\frac{2}{k} \left| D^n_+ D_{t,+}^2 \left( x^N u(x,\tilde{T}-2k)\right) \right| \right) \\
                                                         & \leq & C \left( \left|D^n_+ \left( x^N u(x,t_j) \right) \right| + \left| D_{t,+} D^n_+ \left( x^N u(x,\tilde{T}-k) \right)\right| \right. \\
                                                         & & \left. + \left| D_{t,+}^2 D^n_+ \left( x^N u(x,\tilde{T}-2k)\right) \right| \right) 
\end{eqnarray*} 
\noindent where $\tilde{T}$ is understood as before, and in the same way we may take the $L^2_h$ norm of both sides of this inequality and use the boundedness of each term to deduce that $(\ref{L9.1})$ holds for $m=1$ and $t_j \in \left[ T,T+2\right]$.  Thus we may take $k_0 = k_2$ in the case $m=1$ as well.  \\
\indent Now we shall prove the lemma for $m=2$ which has already been established when $t_j, t_{j+1},t_{j+2}\in \left[0,T\right]$.\\
\indent If $t_j <0 $ then since 
\begin{eqnarray*}
D_{t,+}\hat{u}(\cdot,t_j) & = & \frac{1}{k} \left[u(\cdot,t_{j+1}) \phi(t_{j+1}) - u(\cdot,t_j) \phi(t_j) \right] \\
                          & = & \frac{1}{k} \left[ u(\cdot,t_{j+1}) \phi(t_{j+1}) - \left( u(\cdot,t_{j+1}) - k D_{t,+}u(\cdot,t_{j+1}) + k^2 D_{t,+}^2u(\cdot,t_{j+1}) \right) \phi(t_j) \right] \\
                          & = & \frac{1}{k} \left[ u(\cdot,t_{j+1}) \left(\phi(t_{j+1}) - \phi(t_j) \right) + \phi(t_j) \left( k D_{t,+}u(\cdot,t_{j+1}) - k^2 D_{t,+}^2u(\cdot,t_{j+1}) \right) \right] \\
                          & = &  u(\cdot,t_{j+1}) D_{t,+}\phi(t_j) + \phi(t_j)D_{t,+}u(\cdot,t_{j+1}) - k\phi(t_j) D_{t,+}^2u(\cdot,t_{j+1})
\end{eqnarray*}

\noindent and by the product rule for $D_{t,+}$ we have $$D_{t,+}\hat{u}(\cdot,t_{j+1}) = u(\cdot,t_{j+1}) D_{t,+} \phi(t_{j+1}) + \phi(t_{j+2})D_{t,+}u(\cdot,t_{j+1})$$ It then follows that 
\begin{eqnarray*}
D_{t,+}^2 \hat{u}(\cdot,t_j) & = & \frac{1}{k} \left[ D_{t,+}\hat{u}(\cdot,t_{j+1}) - D_{t,+}\hat{u}(\cdot,t_j)  \right] \\
                               & = & \frac{1}{k} \left[u(\cdot,t_{j+1}) \left( D_{t,+} \phi(t_{j+1}) - D_{t,+} \phi(t_j) \right) + \left( \phi(t_{j+2}) - \phi(t_j) \right)D_{t,+}u(\cdot,t_{j+1}) \right. \\
                               &   & \hspace{1cm} \left. + k \phi(t_j)D_{t,+}^2u(\cdot,t_{j+1}) \right] \\
                               & = & u(\cdot,t_{j+1}) D_{t,+}^2 \phi(t_j) + D_{t,+} u(\cdot,t_{j+1}) D_{t,+}\phi(t_{j+1}) \\
                               &   &  \hspace{1cm}+ D_{t,+}u(\cdot,t_{j+1}) D_{t,+} \phi(t_j)+ \phi(t_j) D_{t,+}^2u(\cdot,t_{j+1}) \\
                               & = & u(\cdot,t_{j+1}) D_{t,+}^2 \phi(t_j) + D_{t,+} u(\cdot,t_{j+1}) D_{t,+}\phi(t_{j+1}) \\
                               &   & \hspace{1cm}+ D_{t,+}u(\cdot,t_{j+1}) D_{t,+} \phi(t_j)+ \phi(t_j) D_{t,+}^2u(\cdot,0)
\end{eqnarray*}

\noindent Therefore we have, 
\begin{eqnarray*}
\left| D_{t,+}^2 D^n_+ \left(x^N \hat{u}(x,t_j)\right) \right| & \leq & \left| D^n_+ \left( x^N u(x,t_{j+1}) \right) \right| \cdot \left|  D_{t,+}^2 \phi(t_j) \right| + \left| D^n_+  D_{t,+}\left( x^N u(x,t_{j+1}) \right) \right| \cdot \left| D_{t,+}\phi(t_{j+1}) \right|   \\
                                                            & & + \left| D^n_+ D_{t,+} \left( x^N u(x,t_{j+1})\right) \right| \cdot \left| D_{t,+} \phi(t_j) \right| +\left|\phi(t_j) \right| \cdot \left| D^n_+ D_{t,+}^2 \left( x^N u(x,0) \right) \right|\\
                                                            & \leq & \left| D^n_+ \left( x^N u(x,t_{j+1})\right) \right| \cdot \SupNorm{  \frac{d^2}{dt^2} \phi } + \left| D^n_+  D_{t,+} \left( x^N u(x,t_{j+1}) \right)\right| \cdot \SupNorm{\frac{d}{dt}\phi}  \\
                                                            & & + \left| D^n_+ D_{t,+} \left( x^N u(x,t_{j+1})\right) \right| \cdot \SupNorm{ \frac{d}{dt} \phi } +\SupNorm{\phi}\left| D^n_+ D_{t,+}^2 \left( x^N u(x,0) \right) \right| \\
                                                            & \leq & C \left| D^n_+ \left( x^N u(x,t_{j+1}) \right) \right| + C \left| D_{t,+} D^n_+ \left( x^N u(x,t_{j+1})\right) \right| + C \left| D_{t,+}^2 D^n_+ \left( x^N u(x,0) \right) \right|
\end{eqnarray*}

\noindent By taking $L^2_h$ norm of the above inequality and using the fact that $(\ref{L9.1})$ holds with $m=0,1$ for all $t_j \in \mathbb{R}_h$ and with $m=2$ for $t_0 \in \left[0,T\right]$ we deduce that $(\ref{L9.1})$ holds with $m=2$ for any $t_j <0$.  \\
\indent If $t_{j+2} > T$ then, by equalities similar to those shown in the case $t_j <0 $ we see that 
\begin{eqnarray*} 
D_{t,+}\hat{u}(\cdot,t_{j+1}) & = & u(\cdot,t_{j+1}) D_{t,+} \phi(t_{j+1}) + \phi(t_{j+2}) D_{t_+} u(\cdot,t_j)+ k\phi(t_{j+2}) D_{t,+}^2u(\cdot,t_{j-1})
\end{eqnarray*}

\noindent and by the product rule for $D_{t,+}$ we have $$D_{t,+}\hat{u}(\cdot,t_j) = u(\cdot,t_{j+1}) D_{t,+} \phi(t_j) + \phi(t_j)D_{t,+}u(\cdot,t_j)$$ It then follows that 

\begin{eqnarray*}
D_{t,+}^2 \hat{u}(\cdot,t_j) & = & u(\cdot,t_{j+1}) D_{t,+}^2 \phi(t_j) + D_{t,+} u(\cdot,t_j) D_{t,+}\phi(t_{j+1}) + D_{t,+}u(\cdot,t_j) D_{t,+} \phi(t_j)\\
                               &   &  \hspace{1cm}+ \phi(t_{j+2}) D_{t,+}^2u(\cdot,t_{j-1})
\end{eqnarray*} 

\noindent Therefore we have,

\begin{eqnarray*}
\left| D_{t,+}^2 D^n_+ \left( x^N \hat{u}(x,t_j) \right) \right| & \leq & C \left| D^n_+ \left( x^N u(x,t_{j+1}) \right) \right| + C \left| D^n_+  D_{t,+} \left( x^N u(x,t_j) \right) \right|+ C \left| D^n_+  D_{t,+}^2 \left( x^N u(x,t_{j-1})\right) \right| \\
                                                            & = & C \left( \left| D^n_+  \left( x^N u(x,t_{j+1}) \right) \right| + \left|D_{t,+} D^n_+ \left( x^N u(x,t_j)\right) \right|+ \left| D_{t,+}^2  D^n_+ \left( x^N u(x,\tilde{T}-2k) \right) \right| \right)
\end{eqnarray*}
\noindent where $\tilde{T}$ is again understood as the largest $t_j\in \left[0,T\right]\cap \mathbb{R}_k$.  By taking $L^2_h$ norm of the above inequality and using the fact that $(\ref{L9.1})$ holds for all $t_j \in \mathbb{R}_h$ with $m=0,1$ and for $\tilde{T}-2k \in \left[0,T\right]$ with $m=2$ we deduce that $(\ref{L9.1})$ holds with $m=2$ for any $t_j$ such that $t_{j+2}>T$.  Therefore we have proven that the lemma holds in the case $m=2$ for $k \in (0,k_2)$ thus we may take $k_0 = k_2$.  \\
\indent  Now we shall prove the lemma in the case $m=3$.  We may assume that $k_2$ is sufficiently small so that $t_3\in\left[0,T\right]$.  If $t_j,t_{j+1}, t_{j+2}, t_{j+3}\in \left[0,T\right]$ then $(\ref{L9.1})$ holds for $k\in (0,k_3)$ for some $0<k_3 \leq k_2$ just as in the cases $m=1,2$  because $\hat{u}$ is a solution to the difference scheme $(\ref{star2})$ and so we could repeatedly replace $D_{t,+} \hat{u}(\cdot, t_j)$ by the rest of $(\ref{star2})$ and use some elementary inequalities together with corollary 3.4 to deduce $(\ref{L9.1})$. \\
\indent By the product rule for $D_{t,+}$ we have
\begin{eqnarray*}
D_{t,+}^3 \hat{u}(\cdot, t_j) & = & u(\cdot, t_{j+1}) D_{t,+}^3 \phi(t_j)+ D_{t,+}u(\cdot, t_{j+1})D_{t,+}^2 \phi(t_{j+1})+D_{t,+}^2 u(\cdot, t_{j+1})D_{t,+}\phi(t_{j+2})\\ 
                              &   & + D_{t,+}u(\cdot,t_{j+1})D_{t,+}^2 \phi(t_{j+1}) +D_{t,+}^2 \phi(t_j)D_{t,+}u(\cdot,t_{j+1}) + D_{t,+}\phi(t_j)D_{t,+}^2u(\cdot, t_j) \\
                              &   & + D_{t,+}\phi(t_{j+1}) D_{t,+}^2u(\cdot,t_{j+1}) + \phi(t_{j+1})D_{t,+}^3 u(\cdot,t_j)
\end{eqnarray*}

\noindent If $t_j < 0$ or if $t_{j+3} > T$ then by construction we have $D_{t,+}^3 u(\cdot,t_j) = 0$.  Therefore, after applying $D_+^n \left( x^N \cdot \right)$ to both sides and taking the $L^2_h$ norm we see that $(\ref{L9.1})$ follows for $m=3$ because it holds for $m=0,1,2$ and because $\phi\in C_c^{\infty}(\mathbb{R})$.  Thus we may take $k_0=k_3$ in the case $m=3$.  $\square$ 

\section{Obtaining Smooth Solutions from Discrete Solutions}
\subsection{The Smoothing Operator $I_h$}
\noindent Lemma 4.1 and corollary 4.2 are based on similar statements found in [1,8].  They are the key ingredients allowing us to pass from a discrete function to a continuum function while preserving the necessary estimates for our solution (i.e. boundedness of Schwartz semi-norms).  The proofs of lemma 4.1 and corollary 4.2 can be found in [1] where the author uses ideas from [8].\\
\indent  We will denote by $L^2\left(\mathbb{R}\right)$ (or simply $L^2$) the space of square integrable functions defined on $\mathbb{R}$ with its usual inner product and norm denoted by $\CLIP{\cdot}{\cdot}$ and $\CLN{\cdot}$ respectively (in contrast to $L^2_h$ which the space of square summable functions defined on the mesh $\mathbb{R}_h$ and whose norm is denoted by $\DLN{\cdot}$).  Clearly if we restrict a continuum function $u \in L^2$ to $\mathbb{R}_h$ then we may consider it also as a mesh function in $L^2_h$.  \\

\noindent \textbf{Lemma 4.1} \hspace{1mm} For any $h>0$ there exists a linear isometry $I_h : L^2_h  \rightarrow L^2$ such that if $u \in L^2_h$ then $U:=I_h u$ has the following properties:
\begin{enumerate}
\item  $U \in C^{\infty}\left(\mathbb{R}\right)$ (hence we can think of $I_h$ as a "smoothing operator").
\item For any point $x_n \in \mathbb{R}_h$ we have that $U (x_n) = u(x_n)$.
\item  For each $j>0$ the following inequalities hold: 
$$\left( \frac{2}{\pi}\right)^j \CLN{\frac{\partial^j}{\partial x^j} U} \leq \DLN{ D^j_+ u} \leq  \CLN{\frac{\partial^j}{\partial x^j} U} $$  
\end{enumerate}

\indent   An explicit formula for $U(x)$ is given by

\begin{equation}
U(x)=\sum_{l=-\infty}^{\infty} u(x_l) \left( \frac{ \sin \frac{\pi}{h}(x_l - x)}{\frac{\pi}{h}(x_l-x)} \right) \label{1.1}
\end{equation}

\vspace{2mm}

\noindent \textbf{Corollary 4.2} \hspace{1mm} Let $M\geq 2$ be an integer and let $h > 0$ be a real number.  Suppose $u$ is a mesh function on $\mathbb{R}_h$ such that $x_n^N u(x_n) \in L^2_h$ for each $0 \leq N \leq M$.  Then for each $j \in \mathbb{N}$, and $0 \leq N \leq M-2 $ we have,
$$\left( \frac{2}{\pi}\right)^j \CLN{ \frac{\partial^j}{\partial x^j} \left( x^N U\right) } \leq \DLN{ D^j_+ \left( x^N u\right) } \leq  \CLN{ \frac{\partial^j}{\partial x^j} \left( x^N U \right) } $$  

\vspace{2mm}

\subsection{Schwartz Boundedness of Smoothly Continued Discrete Solutions}
\noindent In this section we will show that a certain smooth continuation of the discrete solution remains bounded in the Schwartz semi-norms $\SupNormtx{\IHN{\cdot}^N \partial_x^n \partial_t^m \cdot} $.  \\
\indent Suppose $u(x_n,t_j)$ is a mesh function defined on $\mathbb{R}_h \times \mathbb{R}_k$ for some $h,k>0$ which is compactly supported in time for $t_j \in \left[T_0,T_1\right]$ and which satisfies the property that there is some $C>0$ such that $\DLN{u(\cdot,t_j)} \leq C$ for each $t_j \in \mathbb{R}_k$.  We will define a smooth continuation of $u$ by the following:  \\
\indent Since $u$ is compactly supported in time we know that for each $x_n \in \mathbb{R}_h$ we have $u(x_n,\cdot)\in L^2_k$, therefore by lemma 4.1 we may apply the operator $I_k$ to $u(x_n, \cdot)$ in $t$ to obtain by (\ref{1.1}) that for any $t\in \mathbb{R}$,

\begin{equation*}
(I_k u)(x_n,t)  =  \sum_{T_0 \leq t_j \leq T_1 } u(x_n,t_j) \left[ \frac{\sin \frac{\pi}{k}(t-t_j) }{\frac{\pi}{k} (t-t_j)} \right]  
\end{equation*}

\noindent and therefore,

\begin{eqnarray}
\DLN{ I_k u(\cdot,t) } & \leq & \sum_{T_0 \leq t_j \leq T_1 } \left| \frac{\sin \frac{\pi}{k}(t-t_j) }{\frac{\pi}{k} (t-t_j)}  \right| \cdot \DLN{ u(\cdot,t_j)} \nonumber \\
                     & \leq & \sum_{T_0 \leq t_j \leq T_1 } \SupNorm{\frac{\sin y }{y}} \cdot \DLN{ u(\cdot,t_j) } < \infty \label{L10.1}
\end{eqnarray}
 
\noindent Hence by lemma 4.1 we can apply the smoothing operator $I_h$ to $I_k u(\cdot,t)$ in the $x$ variable for each $t\in \mathbb{R}$ to obtain a continuum function $I u:=I_h I_k u$.  By linearity of $I_h$ it follows that $I u$ is given by,
\begin{eqnarray*}
I u(x,t) &=& \sum_{T_0 \leq t_j \leq T_1 } I_h u(x,t_j) \left[ \frac{\sin \frac{\pi}{k}(t-t_j) }{\frac{\pi}{k} (t-t_j)} \right]
\end{eqnarray*}
and since for each $t_j$ we have $\frac{\sin \frac{\pi}{k}(t-t_j)}{\frac{\pi}{k}(t-t_j)}$ is smooth in $t$ and also for each $j$ the function $I_h u(x,t_j)$ is smooth in $x$ we see that $ I u\in C^{\infty}\left(\mathbb{R}\times \mathbb{R} \right)$ and by lemma 4.1 $I u (x_n,t_j) = u(x_n,t_j)$ for any $(x_n,t_j)\in \mathbb{R}_h \times \mathbb{R}_k$.  Moreover, it is clear from the above formula that for any $m \in \mathbb{N}$ we have $\partial_t^m I u = I_h \left( \partial_t^m I_k u \right)$.  Given a discrete solution $u$ from lemma 2.7 we may now construct a smooth continuation $I \hat{u}$ and prove the following lemma.\\

\noindent \textbf{Lemma 4.3} \hspace{1mm} Let $u_0\in S^{-\infty}(\mathbb{R})$ and $K,T,\epsilon >0$ all be given as in lemma 2.7 and let $N,n \in \mathbb{N}$, $m \in \left\{0,1,2\right\}$.  Then there exists $0 < k_0 \leq \epsilon$  and $C>0$ both depending on $m,n$ and $N$ such that if $h \in (0,1)$, $k\in (0,k_0)$, and $u(x_n,t_j)$ is the solution of the difference scheme $(\ref{star2})$ defined on $\mathbb{R}_h \times \left( \mathbb{R}_k \cap \left[0,T\right] \right)$ with initial condition $u_0 $ then we have 

\begin{equation*}
\SupNormtx{\IHN{\cdot}^N \partial_x^n \partial_t^m I \hat{u} }\leq C 
\end{equation*}

\vspace{2mm}

\noindent \textbf{Proof of Lemma 4.3} \hspace{1mm} Fix values for $n$ and $N$.  We may assume $N \geq 2$ because then the lemma will automatically follow for the lower values of $N$.  By lemma 3.5 we have for each $m=0,1,2,3$ values for $k_0(m,n,N)>0$ and $C(m,n,N)>0$ such that $(\ref{L9.1})$ holds for $k\in (0,k_0)$ and $t_j \in \mathbb{R}_k$.  We may also assume that $k_0$ is sufficiently small and $C(m,n,N)$ is sufficiently large so that $(\ref{L9.1})$ holds for $D_{t,+}^m  D^i_+ \left( x^l \hat{u}\right)$ for each $i=0,1,\ldots,n+1$, $l=2,3,\ldots,N+2$, and $m=0,1,2,3$.  Since $\hat{u}$ is zero for all $t_j \notin \left[-2,T+2 \right]$ we may take the $L^2_{t,k}$ norm of $(\ref{L9.1})$ (i.e. the discrete $L^2$ norm in time) to obtain for each such choice of $i$, $l$, $m$, and $k$ that

\begin{equation}
\DLNtx{ D_{t,+}^m  D^i_+ \left( x^l \hat{u} \right) }^2 \leq \frac{T+5}{k} \max_{t_j\in \left[-3,T+2 \right]} \DLN{ D_{t,+}^m D^i_+ \left( x^l \hat{u}(x, t_j)\right)}^2\cdot k \leq C \label{L10.2}
\end{equation}

\noindent  Moreover as in (\ref{L10.1}) we have for each $t\in \mathbb{R}$ and for any $m=0,1,2,3$ 

\begin{eqnarray*}
\DLN{ D_+^i \left( x^l \partial_t^m I_k \hat{u}(\cdot,t)\right) } & \leq & \sum_{-2 \leq t_j \leq T+2 } \left|\frac{d^m}{dt^m} \frac{\sin \frac{\pi}{k}(t-t_j) }{\frac{\pi}{k} (t-t_j)}  \right| \cdot \DLN{D_+^i \left( x^l \hat{u}(\cdot,t_j) \right) } \\
                     & \leq & \sum_{-2 \leq t_j \leq T+2 } \SupNorm{\frac{d^m}{dy^m}\frac{\sin y }{y}} \cdot \DLN{D_+^i \left( x^l \hat{u}(\cdot,t_j) \right) } \leq  C < \infty 
\end{eqnarray*} 

\noindent Therefore we may apply corollary 4.2 in $x$, Fubini's theorem (see [6]), lemma 4.1 in $t$, and (\ref{L10.2}) to conclude \\

\begin{eqnarray}
\CLNStx{\partial_t^m \partial_x^i \left( x^l I \hat{u}\right)} & \leq & C \DLNtx{ D_{t,+}^m  D^i_+ \left( x^l \hat{u}\right) }^2 \hspace{.3cm} \leq C \label{L10.4}
\end{eqnarray}

\indent By repeatedly applying the product rule we see from (\ref{L10.4}) that for any $m=0,1,2,3$ we have,

\begin{equation}
\CLNtx{\IHN{\cdot}^N \partial_x^n \partial_t^m I \hat{u}}  \leq  C\CLNtx{\partial_x^n \partial_t^m I \hat{u}} + C\CLNtx{x^N \partial_x^n \partial_t^m I \hat{u}} \leq C \label{L10.5}
\end{equation}

\noindent Then from the continuum versions of Sobolev's inequalities (see the remark in the appendix) applied in $t$ and in $x$ we obtain also for $k \in \left(0,k_0 \right)$ that when $m=0,1,2$ we have
\begin{eqnarray*}
\SupNormtx{\IHN{x}^N \partial_x^n \partial_t^m I \hat{u}}  &\leq&  C\left( \CLNxt{\IHN{x}^N I \hat{u}} +\CLNxt{\IHN{x}^N \partial_t^{m+1}I \hat{u}} +\right.  \\
                                                            &   & + \left. \CLNxt{\IHN{x}^N \partial_x^{n+1}I \hat{u}} +\CLNxt{\IHN{x}^N \partial_t^{m+1}\partial_x^{n+1}I \hat{u}} \right)  \leq  C  
\end{eqnarray*}

\noindent $\square$ 

\subsection{Proof of Local Existence for the Generalized mKdV equation in $S^{-\infty}$}

\noindent By using corollary 4.2 and the Arzela-Ascoli theorem we shall now construct a smooth solution to $(\ref{daggerdagger})$ lying in $S^{-\infty}\left(\mathbb{R}\times \left[ 0,T \right]\right)$ that comes from the discrete solution constructed in lemma 2.7.  Theorem 1.1 and its proof are completely analogous to the corresponding results given by Bondareva for the KdV equation (see [1] theorem 2). \\

\noindent \textbf{Proof of Theorem 1.1} \hspace{1mm} (\textit{existence}) \\  
\indent Since $u_0\in  S^{-\infty}\left(\mathbb{R}\right)$ it follows that there is some $K>0$ such that for any $0<h<1$ we have $\DSN{u_0}\leq K$.  Therefore, by lemma 2.7, there exists $T, L, \epsilon >0$ such that if $h \in (0,1)$ and $k \in (0,\epsilon)$ then there is a solution to the difference scheme $(\ref{star2})$ with initial condition $u_0$ defined on $\mathbb{R}_h \times \left(\mathbb{R}_k \cap \left[0,T\right]\right)$ and we denote this solution by by $u^{h,k}$.  Let $U^{h,k} := I \hat{u}^{h,k}$.  From lemma 4.3 we know that for every $N,n \in \mathbb{N}$, $m \in \left\{0,1,2\right\}$ there exists $0 < k_0(m,n,N) \leq \epsilon$ and there exists $C_{m,N,n} > 0$ such that if $h \in (0,1)$, $k\in (0,k_0)$, then we have 

\begin{equation}
\SupNormtx{\IHN{\cdot}^N \partial_x^n \partial_t^m U^{h,k} }\leq C  \label{T.1}
\end{equation}

\vspace{2mm}

\noindent From the family of functions $\left\{U^{h,k}\right\}_{h\in (0,1),k\in (0,k_0(1,1,0))}$ we now wish to extract a convergent subsequence by using the Arzela-Ascoli theorem (this theorem can be found for example in [6]).\\
\indent Let $(x_0,t_0),(x_1,t_1)$ be points in $\mathbb{R}\times \mathbb{R}$.  By lemma 4.3 and by the intermediate value theorem we have
\begin{eqnarray}
\left|U^{h,k}(x_0,t_0) - U^{h,k}(x_1,t_1)\right| & \leq & \left|U^{h,k}(x_0,t_0) - U^{h,k}(x_0,t_1)\right| + \left|U^{h,k}(x_0,t_1) - U^{h,k}(x_1,t_1)\right| \nonumber \\
                                                 & \leq & \left| \partial_t U^{h,k}(x_0,\tilde{t}) \right| \cdot \left|t_0-t_1 \right| + \left| \partial_x U^{h,k}(\tilde{x},t_1) \right| \cdot \left|x_0-x_1 \right|  \nonumber \\
                                                 & \leq & \SupNormtx{\partial_t U^{h,k}} \cdot \left| t_0-t_1\right| + \SupNormtx{\partial_x U^{h,k}} \cdot \left|x_0-x_1 \right| \nonumber \\
                                                 & \leq & C \left| (x_0,t_0)-(x_1,t_1) \right|\label{T.2}
\end{eqnarray}
\noindent which shows that the family of functions $U^{h,k}$ is equicontinuous on $\mathbb{R}\times \mathbb{R}$. \\
\indent From (\ref{T.1}) it follows that the family $U^{h,k}$ is also bounded uniformly for $h\in (0,1)$, $k\in (0,k_0(1,1,0))$.  Hence, by the Arzela-Ascoli theorem we may construct a subsequence $U^{h_i,k_i}$, where of course $h_i,k_i \searrow 0 $ as $i \rightarrow \infty$, converging uniformly on compact sets to a function $U\in C^0 \left( \mathbb{R}\times \mathbb{R} \right) $. \\
\indent The above argument can also be made for the family of functions $\partial_x U^{h_i,k_i} $for $h\in (0,1)$, $k\in (0,k_0(1,2,0))$.  Namely, estimate (\ref{T.1}) implies that the family is bounded uniformly and also that we may use estimate (\ref{T.2}) with $U^{h,k} $ replaced by $\partial_x U^{h_i,k_i}$ to see that it is also an equicontinuous family.  Thus we conclude that there is some $V\in C^0\left( \mathbb{R}\times \mathbb{R} \right)$ and a subsequence  $\partial_x U^{h_l,k_l}$ converging uniformly on compact sets to $V$.  Since we have uniform convergence on compact sets for $U^{h_l, k_l}$ and $\partial_x U^{h_l,k_l}$ it follows that $U$ is differentiable in $x$ and $\partial_x U = V$ on $\mathbb{R}\times\mathbb{R}$.  \\
\indent By repeating the same argument we conclude by induction that for each $p\in\mathbb{N}$ the function $\partial_x^{p-1} U\in C^0(\mathbb{R}\times\mathbb{R})$ is differentiable in $x$ because the sequence $\partial_x^p U^{h_l,k_l}$ for $h\in (0,1)$ and $k\in (0,k_0(1,p+1,0))$ is bounded uniformly by (\ref{T.1}) and is equicontinuous by (\ref{T.2})and hence it has a subsequence uniformly convergent on compact subsets of $\mathbb{R}\times \mathbb{R}$ to $\partial_x \partial_x^{p-1} U$.  In this way we will obtain a countable array of subsequences, one for each $p\in \mathbb{N}$ and from this array we extract a diagonal subsequence.  From this diagonal subsequence it will follow that for each $p\in \mathbb{N}$ we have $\partial_x^p U^{h_l,k_l} \rightarrow \partial_p U$ uniformly on compact sets.\\
\indent Consider the family of functions $\partial_t U^{h_l,k_l}$ for $h\in (0,1)$ and $k\in (0,k_0(2,1,0))$.  Estimate (\ref{T.1}) for $m=1$ implies that the family is bounded uniformly and also that we may use estimate (\ref{T.2}) with $U^{h,k}$ replaced by $\partial_t U^{h_i,k_i}$ to see that it is also an equicontinuous family.  Hence we may, as before for $x$, conclude that $U$ is differentiable in $t$ and construct a subsequence of $\partial_t U^{h_l,k_l}$ uniformly convergent on compact sets to $\partial_t U$.  From this subsequence of $(h_l,k_l)$ we consider the family $\partial_x \partial_t U^{h_l,k_l}$ for $h\in (0,1)$ and $k\in (0,k_0(2,2,0))$.  Again estimate (\ref{T.1}) implies that the family is bounded uniformly and also we may use estimate (\ref{T.2}) with $U^{h,k}$ replaced by $\partial_x\partial_t U^{h_l,k_l}$ to see that it is also an equicontinuous family.   Thus we may again extract a subsequence $U^{h_l,k_l}$ to see that $\partial_t U$ is differentiable in $x$ and $\partial_x \partial_t U^{h_l,k_l}\rightarrow \partial_x \partial_t U$ uniformly on compact sets.   Continuing inductively we consider the sequence of functions $\partial_x^p \partial_t U^{h_l,k_l}$ for $h\in (0,1)$, $k\in (0,k_0(2,p+1,0))$.  It is equicontinuous by (\ref{T.2}) and from (\ref{T.1}) it is uniformly bounded, thus we conclude that $\partial_x^{p-1}\partial_t U$ is differentiable in $x$ and we may extract a subsequence so that $\partial_x^p \partial_t U^{h_l,k_l}\rightarrow \partial_x^p \partial_t U$.  \\
\indent Continuing in this way we will again obtain an array of subsequences of $U^{h_l,k_l}$, one for each $p\in \mathbb{N}$.  By taking a diagonal subsequence we obtain a subsequence such that for each $p\in \mathbb{N}$ and for $q=0,1$ we have $\partial_x^p\partial_t^q U^{h_l,k_l} \rightarrow \partial_x^p\partial_t^q U$ uniformly on compact subsets of $\mathbb{R}\times \mathbb{R}$.  In addition, it follows that for any $N\in \mathbb{N}$ the sequence $\IHN{x}^N \partial_x^p \partial_t^q U^{h_l,k_l} \rightarrow \IHN{x}^N \partial_x^p\partial_t^q U$ uniformly on compact subsets of $\mathbb{R}\times \mathbb{R}$ because for any compact set $X \subset \mathbb{R}\times\mathbb{R}$ we have the inequality
\begin{eqnarray*}
\left|\IHN{x}^N \partial_x^p \partial_t^q U^{h,k}(x,t) - \IHN{x}^N \partial_x^p \partial_t^q U(x,t)\right| & \leq & \max_{x\in X}\IHN{x}^N \cdot \left| \partial_x^p \partial_t^q U^{h,k}(x,t) - \partial_x^p\partial_t^q U(x,t)\right|
\end{eqnarray*} 

\noindent By construction we can see that $U$ satisfies the following conditions:\\
\begin{enumerate}
\item $\partial_x^p U$ exists for each $p\in \mathbb{N}$ and is continuous (i.e. it lies in $C^0(\mathbb{R}\times\mathbb{R})$).\\
\item $\partial_x^q \partial_t \partial_x^p U$ exists for each $p,q\in \mathbb{N}$ and is continuous (i.e. it lies in $C^0(\mathbb{R}\times\mathbb{R})$).\\
\item If $p+q = p'+q'$ then $\partial_x^q \partial_t \partial_x^p U = \partial_x^{q'} \partial_t \partial_x^{p'} U$. 
\end{enumerate}

\indent To prove that $U$ is a solution to $(\ref{daggerdagger})$ we fix a point $(x,t)\in \mathbb{R} \times \left[0,T\right]$ and show that the equation is satisfied at $(x,t)$.  From our final subsequence of pairs $(h_l,k_l)$ above we first construct points $(x_l,t_l) \in \mathbb{R}_{h_l}\times\left( \mathbb{R}_{k_l}\cap \left[0,T\right] \right)$ to be the nearest points in the grid to $(x,t)$ (note: in this context $x_l \neq l \cdot h_l$ and $t_l \neq l \cdot k_l$).  It then follows that $(x_l,t_l) \rightarrow (x,t)$ as $l\rightarrow \infty$.  By construction we have for each $l\in \mathbb{N}$ that $U^{h_l,k_l}$ satisfies the difference scheme $(\ref{star2})$ at point $(x_l,t_l)$.  Replace the discrete derivatives in $t$ and $x$ of $(\ref{star2})$ by usual derivatives at intermediate points $(\tilde{x}_l,\tilde{t}_l)$ (possibly different intermediate points for each term).  We will then obtain a sum of products of terms of the form $\partial_t U^{h_l,k_l}(x_l,\tilde{t}_l)$, $U^{h_l,k_l}(x_l,t_l)$, $\partial_x U^{h_l,k_l}(\tilde{x}_l,t_l+k_l)$, $\partial_x^3 U^{h_l,k_l}(\tilde{x}_l,t_l+k_l)$, $f(x_l,t_l)$, $\partial_x f(x_l,t_l)$, $U^{h_l,k_l}(x_l,t_l + k_l)$, and $g(x_l,t_l)$.  By continuity of $f$, $\partial_x f$, and $g$ we see that $f(x_l,t_l) \rightarrow f(x,t)$, $\partial_x f(x_l,t_l) \rightarrow \partial_x f(x,t)$ and $g(x_l,t_l)\rightarrow g(x,t)$ as $l\rightarrow \infty$.  Moreover, since $U^{h_l,k_l}(x,t) \rightarrow U(x,t)$ as $l \rightarrow \infty$ and

\begin{eqnarray*}
\left|U^{h_l,k_l}(x_l,t_l)-U(x,t) \right| & \leq & \left|U^{h_l,k_l}(x_l,t_l) - U^{h_l,k_l}(x,t) \right| + \left|U^{h_l,k_l}(x,t) - U(x,t) \right|
\end{eqnarray*}

\noindent it follows by equicontinuity of the family $U^{h_l,k_l}$ for $h,k$ sufficiently small that $U^{h_l,k_l}(x_l,t_l)\rightarrow U(x,t)$.  We may use the same convergence argument for the other terms in the equation to show that as $l\rightarrow \infty$ the equation becomes $(\ref{daggerdagger})$ at the point $(x,t)$. \\
\indent Since $U$ satisfies the equation $(\ref{daggerdagger})$ it follows that $\partial_t U$ is also differentiable in time and its higher time derivatives can be written in terms of the lower $x$ derivatives.  The derivatives also clearly commute as was mentioned above in condition 3, therefore it follows that $U\in C^{\infty}\left(\mathbb{R}\times \left[0,T\right] \right)$.\\
\indent Furthermore, we can show that the limit function $U$ is in $S^{-\infty}\left(\mathbb{R}\times\left[0,T\right]\right)$.  By taking the limit of $\IHN{x}^N \partial_x^n \partial_t^m U^{h_l,k_l}(x,t)$ as $l \rightarrow \infty$ we can see that (\ref{T.1}) also holds for the function $U$ with $n,N \in \mathbb{N}$ and $m=0,1$.  By repeatedly using the equation $(\ref{daggerdagger})$ we may write $\IHN{\cdot}^N \partial_x^n \partial_t^m U$ as a sum of products of terms of the form $\IHN{\cdot}^N \partial_x^n U$ each of which can be bounded by some constant depending on $N,n\in\mathbb{N}$ by using the limiting case of (\ref{T.1}) and this implies that for any $m,n,N\in \mathbb{N}$ we have, \\ 
\begin{eqnarray}
\SupNormtx{\IHN{\cdot}^N \partial_x^n \partial_t^m U}  &\leq&  C_{m,n,N} \label{T.3} 
\end{eqnarray}
\noindent which shows that $U\in S^{-\infty}\left(\mathbb{R}\times\left[0,T\right]\right)$. \hspace{4cm} $\square$ 

\subsection{Proof of Local Existence for the mKdV equation in $S^{\beta}$ when $\beta \leq \frac{1}{2}$}

\noindent Now we shall now construct smooth solutions to $(\ref{dagger})$ lying in $S^{\beta}\left(\mathbb{R}\times \left[ 0,T \right]\right)$ for $\beta \leq \frac{1}{2}$ that come from adding an above solution to $(\ref{daggerdagger})$ to the asymptotic solution constructed in lemma 6.2. \\

\noindent \textbf{Proof of Theorem 1.2} \hspace{1mm} (\textit{existence}) \\
\indent  By lemma 6.2 there exists an asymptotic solution $f(x,t) \in S^{\beta}(\mathbb{R}\times \left[\left.-c,\infty\right.\right))$ of the initial value problem $(\ref{dagger})$ for some $c>0$ whose expansion coefficients satisfy the desired property.  Let $u_0(x) = w_0(x)-f(x,0)$ and let $g: = f_t+f^2f_x + f_{xxx}$.   By construction $u_0 \in S^{-\infty}(\mathbb{R})$.  Moreover $f$ and $g$ satisfy the hypotheses of theorem 1.1. Therefore there exists a $T>0$ and a solution $u(x,t) \in S^{-\infty}(\mathbb{R}\times \left[0,T\right])$ to equation $(\ref{daggerdagger})$.  Let $w(x,t):=u(x,t)+f(x,t)$.  Since $u$ satisfies $(\ref{daggerdagger})$ it follows that $w$ satisfies $(\ref{dagger})$. 
 Moreover, since $u \sim 0$ it follows that $w$ and $f$ have the same asymptotic expansions and in particular the coefficients in the asymptotic expansions of $w$ satisfy the second statement of the theorem.  Finally, since $f(x,t) \in S^{\beta}(\mathbb{R}\times \left[0,T\right])$ it follows that $w(x,t) \in S^{\beta}(\mathbb{R}\times \left[0,T\right])$.  \hspace{2cm}  $\square$ \\

\section{Uniqueness of Solutions}
\subsection{Uniqueness in $S^{-\infty}$ for the Generalized mKdV Equation}
\noindent In this section we shall prove uniqueness of solutions in $S^{-\infty}(\mathbb{R}\times\left[0,T\right])$ for $(\ref{daggerdagger})$ by using Gronwall's Inequality.  We shall state this inequality and prove it in appendix B for completeness.  The below statement also appears in [5]. \\

\noindent \textbf{Lemma 5.1} \hspace{1mm} (Gronwall's Inequality)  \\
\noindent Let $T>0$ and $c_1,c_2 \in \mathbb{R}$ be given and $c_1 \neq 0$.  Suppose $\eta: \left[0,T\right] \rightarrow \mathbb{R}$ is a nonnegative, differentiable function and that for each $t\in \left[0,T\right]$ we have\\
\begin{equation*}
\frac{d\eta}{dt}(t)  \leq  c_1 \eta(t) + c_2
\end{equation*}
\noindent Then for each $t\in \left[0,T\right]$ we have\\
\begin{equation*}
\eta(t)  \leq  e^{c_1 t}\left(\eta(0)-\frac{c_1}{c_2}\right) - \frac{c_1}{c_2}
\end{equation*}

\noindent \textbf{Proof of Theorem 1.1} \hspace{1mm} (\textit{uniqueness}) \\
\indent  Suppose $u(x,t),v(x,t)\in S^{-\infty}(\mathbb{R}\times\left[0,T\right])$ are two solutions of $(\ref{daggerdagger})$ with initial data $u_0\in S^{-\infty}(\mathbb{R})$.  Then,

\begin{equation*}
u_t+u^2 u_x+u_{xxx}+(u^2 f)_x+(uf^2)_x+g  =  0
\end{equation*}

\noindent and

\begin{equation*}
v_t+v^2 v_x+v_{xxx}+(v^2 f)_x+(vf^2)_x+g  =  0
\end{equation*}

\noindent Let $q(x,t) := u(x,t)-v(x,t)$.  By subtracting the above two equations we see that $q$ satisfies the equation

\begin{equation*}
q_t+q u_xu+q u_x v+v^2 q_x+q_{xxx}+ (fqu)_x+(fqv)_x+(f^2q)_x  =  0
\end{equation*}

\noindent Thus if we multiply by $q$ and integrate in $x$ over $(-\infty,\infty)$ we get

\begin{equation*}
\frac{1}{2} \frac{d}{dt} \int_{\mathbb{R}} q^2 dx + \int_{\mathbb{R}} q^2 u_x u dx + \int_{\mathbb{R}} q^2 u_x v dx + \int_{\mathbb{R}} q q_x v^2 dx + \int_{\mathbb{R}} q q_{xxx} dx + \int_{\mathbb{R}} q(fqu)_x dx + \int_{\mathbb{R}} q(fqv)_x dx + \int_{\mathbb{R}} q(f^2 q)_x dx = 0
\end{equation*}

\noindent since $u,v\in S^{-\infty}(\mathbb{R}\times\left[0,T\right])$ and $f^{(n)}(x,t) = O(|x|^{ \frac{1}{2} - n})$ it follows that 
 after integrating by parts we may obtain the estimate

\begin{eqnarray*}
\frac{d}{dt} \int_{\mathbb{R}} q^2 dx & = & -2 \int_{\mathbb{R}} q^2 u_x u dx -2 \int_{\mathbb{R}} q^2 u_x v dx + \int_{\mathbb{R}} q^2 (v^2)_x dx -\int_{\mathbb{R}} q^2(fu)_x dx -\int_{\mathbb{R}} q^2(fv)_x dx - \int_{\mathbb{R}} q^2(f^2)_x dx \\
                                            & \leq &  2\SupNormtx{u_x}\SupNormtx{u} \int_{\mathbb{R}} q^2 dx +2\SupNormtx{u_x}\SupNormtx{v} \int_{\mathbb{R}} q^2 dx +2\SupNormtx{v_x}\SupNormtx{v} \int_{\mathbb{R}} q^2 dx \\
                                            & & +\SupNormtx{(fu)_x} \int_{\mathbb{R}} q^2 dx +\SupNormtx{(fv)_x} \int_{\mathbb{R}} q^2 dx +\SupNormtx{(f^2)_x} \int_{\mathbb{R}} q^2 dx \leq  C_{u,v,f,T} \int_{\mathbb{R}} q^2 dx
\end{eqnarray*}
\noindent and moreover, $q(\cdot,0) = 0$, therefore by lemma 5.1 it follows that $\CLN{q(\cdot,t)} = 0$ for all $t\in \left[0,T\right]$ and since $q$ is smooth this implies that $q(x,t)=0$ for all $(x,t) \in \mathbb{R}\times\left[0,T\right]$. \hspace{2cm} $\square$ 

\subsection{Uniqueness in $S^{\beta}$ for the mKdV Equation When $\beta \leq \frac{1}{2}$}
\noindent In this section we shall prove uniqueness of solutions in $S^{\beta}(\mathbb{R}\times\left[0,T\right])$ for $(\ref{dagger})$ when $\beta \leq \frac{1}{2}$.  First we will need the following lemma. \\

\noindent \textbf{Lemma 5.2} \hspace{1mm} Let $I\subset \mathbb{R}$ be an interval and $\beta \leq \frac{1}{2}$.  Suppose $w(x,t) \in S^{\beta}(\mathbb{R}\times I)$ is a solution to $(\ref{dagger})$ with initial data $w_0\in S^{\beta}(\mathbb{R})$ and that $ w(x,t) \sim \sum_{k=0}^{\infty} a_k^{\pm}(t) x^{\beta_k}$ as $x\to \pm \infty$.  Then $\sum_{k=0}^{\infty} a_k^{\pm}(t) x^{\beta_k}$ is a formal solution to $(\ref{dagger})$.  \\

\noindent \textbf{Proof of Lemma 5.2} \hspace{1mm} By symmetry it suffices to show that the positive $x$ asymptotic expansion satisfies equation $(\ref{dagger})$.   Let $A_0=\{\beta_j \}_{j=0}^{\infty}$ and let $J \subset I$ be a compact interval.  We enlarge $A_0$ to the set $\Gamma$ defined in appendix A having the properties mentioned in lemma 6.1.  Let us re-write the asymptotic expansion as $\sum_{k=0}^{\infty} a_k^+(t) x^{\gamma_k}$ where $a_k^+(t)=0$ if $\gamma_k \notin A_0$.  By definition of being asymptotic it follows that for every $N \in \mathbb{N}$ we may write $$w(x,t)=\sum_{k=0}^N a_k^+(t) x^{\gamma_k} + R_N(x,t)$$ for $x>1$ and $t\in J$ where $\partial_t^i \partial_x^j R_N(x,t) = O\left(|x|^{\gamma_{N+1}-j} \right)$ for every $i,j \in \mathbb{N}$.  Let $f_N(x,t) = \sum_{k=0}^N a_k^+(t)x^{\gamma_k}$.  Then from $(\ref{dagger})$ it follows that $$(f_N)_t + f_N^2 (f_N)_x + (f_N)_{xxx} + (f_N R_N^2)_x + (f_N^2 R_N)_x + (R_N)_t + R_N^2 (R_N)_x + (R_N)_{xxx} = 0$$  From this we obtain, as we see in (A.1), for some $M \leq N$ 
\begin{equation}
\sum_{j=0}^M \left[ \dot{a}_j^+ + \sum_{{k,l,m}\atop {\gamma_k+\gamma_l+\gamma_m-1=\gamma_j}} a_k^+(t)\cdot a_l^+(t) \cdot a_m^+(t) \cdot \gamma_m + \sum_{{p}\atop{\gamma_p-3=\gamma_j}}a_p^+(t) \cdot \gamma_p \cdot (\gamma_p - 1) (\gamma_p - 2) \right] x^{\gamma_j} + O\left(|x|^{2\gamma_0+\gamma_{N+1}-1}\right) =  0 \label{L12.1}
\end{equation}

\noindent We may assume that $N$ is sufficiently large so that $M \geq 1$ and $2\gamma_0+\gamma_{N+1}-1 < \gamma_1$.  Since the above equation must hold for all $x>1$ we may divide by $x^{\gamma_0}$ to obtain from (\ref{L12.1}) that $$\dot{a}_0^++\sum_{{k,l,m}\atop {\gamma_k+\gamma_l+\gamma_m-1=\gamma_0}} a_k^+(t)\cdot a_l^+(t) \cdot a_m^+(t) \cdot \gamma_m + O(|x|^{\gamma_1-\gamma_0}) =  0$$ and hence $$\dot{a}_0^+ +\sum_{{k,l,m}\atop {\gamma_k+\gamma_l+\gamma_m-1=\gamma_0}} a_k^+(t)\cdot a_l^+(t) \cdot a_m^+(t) \cdot \gamma_m = 0$$  
\indent Continuing in the same way we may assume that $N$ is sufficiently large so that $2\gamma_0+\gamma_{N+1}-1 < \gamma_2$.  Dividing (\ref{L12.1}) by $x^{\gamma_1}$ we obtain that $$\dot{a}_1^++\sum_{{k,l,m}\atop {\gamma_k+\gamma_l+\gamma_m-1=\gamma_1}} a_k^+(t)\cdot a_l^+(t) \cdot a_m^+(t) \cdot \gamma_m + O(|x|^{\gamma_2-\gamma_1}) =  0$$ and hence $$\dot{a}_1^+ +\sum_{{k,l,m}\atop {\gamma_k+\gamma_l+\gamma_m-1=\gamma_1}} a_k^+(t)\cdot a_l^+(t) \cdot a_m^+(t) \cdot \gamma_m = 0$$  

\indent This process may be repeated inductively to obtain from $(\ref{L12.1})$ that for any $j\in \mathbb{N}$ we have $$ \dot{a}_j^+ + \sum_{{k,l,m}\atop {\gamma_k+\gamma_l+\gamma_m-1=\gamma_j}} a_k^+(t)\cdot a_l^+(t) \cdot a_m^+(t) \cdot \gamma_m + \sum_{{p}\atop{\gamma_p-3=\gamma_j}}a_p^+(t) \cdot \gamma_p \cdot (\gamma_p - 1) (\gamma_p - 2) = 0$$  and hence $\sum_{k=0}^{\infty} a_k^+(t)x^{\gamma_k}$ is a formal solution to $(\ref{dagger})$. \hspace{2cm} $\square$ \\

\noindent \textbf{Proof of Theorem 1.2} \hspace{1mm} (\textit{uniqueness}) \\
\indent Suppose $w(x,t),r(x,t)\in S^{\beta}(\mathbb{R}\times \left[0,T\right])$ are two solutions of $(\ref{dagger})$ with initial data $w_0(x)\in S^{\beta}(\mathbb{R})$ and that 

$$w_0(x) \sim \sum_{k=0}^{\infty} c_k^{\pm} x^{\beta_k} \hspace{1cm} w(x,t) \sim \sum_{k=0}^{\infty} a_k^{\pm}(t) x^{\alpha_k} \hspace{1cm} r(x,t) \sim \sum_{k=0}^{\infty} d_k^{\pm}(t) x^{\delta_k}$$ as $x\to \pm \infty$.  Let $B_0=\left\{\beta_k\right\}_{k=0}^{\infty}$, $A_0=\left\{\alpha_k\right\}_{k=0}^{\infty}$, and $D_0=\left\{\delta_k\right\}_{k=0}^{\infty}$.  By lemma 5.2 $\sum_{k=0}^{\infty} a_k^{\pm}(t) x^{\alpha_k}$ and $\sum_{k=0}^{\infty} d_k^{\pm}(t) x^{\delta_k}$ are formal solutions to $(\ref{dagger})$ with initial data $\sum_{k=0}^{\infty} c_k^{\pm} x^{\beta_k}$ and hence we may assume that $B_0 \subset A_0$ and $B_0 \subset D_0$.  Let $\Lambda = A_0 \cup D_0$ and $\Gamma=\left\{\gamma_k\right\}_{k=0}^{\infty}$ be the set constructed in appendix A from $\Lambda$ having the properties stated in lemma 6.1.  Then after reindexing we may rewrite the asymptotic expansions for $w_0$, $w(x,t)$, and $r(x,t)$ as  $$w_0(x) \sim \sum_{k=0}^{\infty} c_k^{\pm} x^{\gamma_k} \hspace{1cm} w(x,t) \sim \sum_{k=0}^{\infty} a_k^{\pm}(t) x^{\gamma_k} \hspace{1cm} r(x,t) \sim \sum_{k=0}^{\infty} d_k^{\pm}(t) x^{\gamma_k}$$ where $c_k=0$ if $\gamma_k \notin B_0$, $a_k(t)=0$ if $\gamma_k\notin A_0$, and $d_k(t) = 0$ if $\gamma_k \notin D_0$.  Since $\sum_{k=0}^{\infty} a_k^{\pm}(t) x^{\gamma_k}$ and $\sum_{k=0}^{\infty} d_k^{\pm}(t) x^{\gamma_k}$ are formal solutions with initial data $\sum_{k=0}^{\infty} c_k^{\pm} x^{\gamma_k}$ it follows that the coefficients $a_k(t)$ and $d_k(t)$ both satisfy the same equations (A.1) with the same initial data and hence for all $k\in \mathbb{N}$  and all $t\in \left[0,T\right]$ we have $a_k(t) = d_k(t)$ so that $w(x,t)-r(x,t) \in S^{-\infty}(\mathbb{R} \times \left[0,T\right])$.  \\
\indent Let $u(x,t) = w(x,t)-r(x,t)$.  Then $u(x,t)$ satisfies $(\ref{daggerdagger})$ with initial condition $u_0(x)=0$ and where $f(x,t) = r(x,t)$ and $g(x,t) = 0$.  By uniqueness of solutions to $(\ref{daggerdagger})$ in $S^{-\infty}(\mathbb{R}\times \left[0,T\right])$, which was proven in theorem 1.1, it follows that $u(x,t) = 0$ for all $(x,t) \in \mathbb{R}\times \left[0,T\right]$. \hspace{2cm} $\square$ \\

\noindent \textbf{ \large{Acknowledgements}} \\

\indent  The author wishes to express his sincere gratitude to P. Topalov for his guidance during the work.  He would also like to thank M. Shubin for providing inspiration to pursue this research.  \\

\noindent \textbf{ \large{Appendix A:  Existence of an Asymptotic Solution}} \\

\noindent In this section we will prove existence of an asymptotic solution to $(\ref{dagger})$.  First we will need the following lemma.  \\
\indent  Let $A_0=\{\beta_j \}_{j=0}^{\infty}$.  Where $\frac{1}{2} \geq \beta_0> \beta_1 > \cdots$, and $\lim_{j\to \infty}\beta_j =-\infty$.  We enlarge $A_0$ to the set $\Gamma$ given by,
\begin{equation*}
\Gamma := \left\{ \sum_{p=1}^k\beta_{i_p} - 3l - \frac{k-1}{2}: k\geq 1, l\geq 0, k,l \in \mathbb{Z}, \beta_{i_p}\in A_0\right\}
\end{equation*}

\vspace{2mm}

\noindent \textbf{Lemma 6.1} \hspace{1mm} The set $\Gamma$ has the following properties:
\begin{enumerate}
\item $A_0 \subset \Gamma$
\item $\Gamma$ is countable.
\item $\Gamma$ is bounded above by $\beta_0$.
\item If $\gamma_l$, $\gamma_m$, $\gamma_n$ are all in $\Gamma$ then $\gamma_l+\gamma_m+\gamma_n - 1$ is in $\Gamma$.
\item If $\gamma_p$ is in $\Gamma$ then $\gamma_p - 3$ is also in $\Gamma$.
\item $\Gamma$ is discrete, i.e. all points are isolated.  
\end{enumerate}

\vspace{2mm}

\noindent \textbf{Proof of Lemma 6.1} \hspace{1mm} Statements 1 to 5 follow easily from the definition of $\Gamma$ so we shall only prove discreteness here.\\
\indent  Let $A:=\left\{\frac{1}{2}-\beta_j:\beta_j\in A_0 \right\}$.  Since $A_0$ is discrete it follows that $A$ is also discrete, hence $\Sigma A = \left\{\sum_{p=1}^k \delta_{p}:\delta_{p}\in A, k\geq 1 \right\}$ is discrete.  Therefore $\frac{1}{2}-\Big(3\mathbb{N}+\Sigma A \Big) = \Gamma $ is discrete. $\square$ \\

\noindent \textbf{Lemma 6.2} \hspace{1mm} For any $\beta \leq \frac{1}{2}$ and for any initial condition $w_0 \in S^{\beta}(\mathbb{R})$ there exists an asymptotic solution $f(x,t) \in S^{\beta}(\mathbb{R}\times I)$ of the initial value problem $(\ref{dagger})$ where $I=\mathbb{R}$ if $\beta<\frac{1}{2}$ and $I=\left[\left.-c,\infty\right.\right)$ for some $c>0$ if $\beta=\frac{1}{2}$.  Moreover, if $w_0 \sim \sum_{k=0}^{\infty} a_k^{\pm} x^{\beta_k}$ and $j$ is the smallest index such that $a_j^{+} \ne 0$ (resp. $a_j^{-} \ne 0$) then the coefficient $a_j^{+}(t)$ (resp. $a_j^{-}(t)$) in the asymptotic expansion of the solution is a nonvanishing continuous function of $t$ and all preceeding coefficients are identically zero.  \\

\noindent \textbf{Proof of Lemma 6.2} \hspace{1mm} First we will show how to construct a formal solution $\sum_{k=0}^{\infty} a_k^{\pm}(t) x^{\beta_k}$.
By symmetry it suffices to construct only the positive $x$ formal solution.  For simplicity we shall omit the superscript $+$ sign in the coefficients $a_j(t)$. \\
\indent First we enlarge the exponent set $A_0=\{\beta_j \}_{j=0}^{\infty}$ to the set $\Gamma$ as defined above.  From lemma 6.1 it follows that we may write the set $\Gamma$ as a decreasing sequence $\Gamma=\left\{\gamma_j\right\}_{j=0}^{\infty}$ where $\frac{1}{2}\geq \gamma_0$, $\gamma_j >\gamma_{j+1}$, and $\gamma_j \rightarrow -\infty$ as $j\rightarrow \infty$, and we may rewrite the positive $x$ asymptotic expansion of $w_0$ as $\sum_{j=0}^{\infty}a_j x^{\gamma_j}$ where $a_j=0$ if $\gamma_j \notin A_0$.  In order to construct the formal solution we need to solve for the coefficients $a_j(t)$ of $x^{\gamma_j}$.  If $\sum_{j=0}^{\infty}a_j(t) x^{\gamma_j}$ is the positive $x$ formal solution to $(\ref{dagger})$ then,
\begin{equation*}
\sum_{j=0}^{\infty}\dot{a}_j(t) x^{\gamma_j} =  - \Big(\sum_{j=0}^{\infty} a_j(t) x^{\gamma_j}\Big)^2 \cdot\Big(\sum_{j=0}^{\infty}a_j(t)\cdot \gamma_j \cdot x^{\gamma_j-1} \Big) - \Big(\sum_{j=0}^{\infty}a_j(t)\cdot \gamma_j \cdot( \gamma_j -1)\cdot (\gamma_j -2) x^{\gamma_j-3}\Big)
\end{equation*}
from which we deduce that the coefficients $a_j(t)$ satisfy the equations,
\begin{equation}
\dot{a}_j = \sum_{{k,l,m}\atop {\gamma_k+\gamma_l+\gamma_m-1=\gamma_j}} - a_k(t)\cdot a_l(t) \cdot a_m(t) \cdot \gamma_m - \sum_{{p}\atop{\gamma_p-3=\gamma_j}}a_p(t) \cdot \gamma_p \cdot (\gamma_p - 1) (\gamma_p - 2) \tag{A.1}
\end{equation}
\indent First we will consider the case when $\gamma_0 < \frac{1}{2}$.  Notice first that for $j=0$ the second sum is nonexistent since $\gamma_0 \geq \gamma_p$ for all $p\geq 0$ and hence there is no $p\geq 0 $ such that $\gamma_p-3 = \gamma_0$.  Also for $j=0$ the first sum is nonexistent because $\gamma_k+\gamma_l+\gamma_m-1 \leq 3\gamma_0 - 1$ and if $\gamma_0 = \gamma_k+\gamma_l+\gamma_m-1$ then $\gamma_0 \leq 3\gamma_0-1$ and hence $\gamma_0 \geq \frac{1}{2}$ which is a contradiction to our assumption that $\gamma_0 <\frac{1}{2}$.  Thus we have $\dot{a}_0=0$ and hence $a_0(t) = a_0$ is constant.  Moreover, for $j\geq 0$ we can see that both sums only contain indices less than $j$.  To see this let us first consider the second sum.  If $\gamma_p-3 = \gamma_j$ then $\gamma_p = \gamma_j+3>\gamma_j$ and hence $p<j$.  For the first sum, if $\gamma_j = \gamma_k+\gamma_l+\gamma_m-1$ and $k \geq j$ then $0\leq \gamma_j-\gamma_k = \gamma_l+\gamma_m-1 $ so that $\gamma_l+\gamma_m\geq 1$, but $\gamma_l,\gamma_m <\frac{1}{2}$, so this is a contradiction, thus $k < j$.  The same argument shows that $l< j$ and $m< j$.  Therefore we may solve for $a_j$ recursively by integrating the right side of the equation to obtain a polynomial in $t$.  By construction the polynomial will be identically zero for the first few indices until we reach $a_j \ne 0$, then it will be a constant $a_j(t) = a_j$, and for all larger indices $a_j(t)$ is polynomial and hence each $a_j(t)$ is defined for all $t\in \mathbb{R}$.  \\
\indent Now let us assume that $\gamma_0 = \frac{1}{2}$.  When $j=0$ the second sum is again nonexistent for the same reason given above however the first sum is nonzero.  If $\gamma_0 = \gamma_k+\gamma_l+\gamma_m - 1$ then $\gamma_k+\gamma_l+\gamma_m = \frac{3}{2}$ and hence $\gamma_k = \gamma_l = \gamma_m = \frac{1}{2}$.  Therefore $\dot{a}_0 = -\frac{1}{2}a_0^3$ which implies that $a_0 = \pm (t+c)^{-\frac{1}{2}}$ for some $c>0$ and hence $a_0$ is continuous for $t\in \left[-c\left. , \infty \right) \right.$.  Furthermore, when $j \geq 1$ the second sum always contains indices less than $j$ however for the first sum at least two of the three indices must be less than $j$.  To see this, suppose that $\gamma_j = \gamma_k+\gamma_l+\gamma_m -1$ and suppose that $k,l \geq j$.  Then $0 \leq \gamma_j-\gamma_k = \gamma_l+\gamma_m - 1$ and hence $\gamma_l+\gamma_m \geq 1$ which implies that $\gamma_l=\gamma_m = \frac{1}{2}$ and hence $l=m=0<j$ which is a contradiction.  For the terms where exactly one index is equal to $\gamma_j $ it must follow that the other two indices are $\gamma_0$ because $\gamma_j = \gamma_j+\gamma_k+\gamma_l - 1$ implies that $\gamma_k+\gamma_l = 1$ and hence $\gamma_k = \gamma_l= \frac{1}{2} = \gamma_0$.  Therefore when $j \geq 1$ the coefficient $a_j(t)$ satisfies the equation $$\dot{a}_j= -\frac{3}{2}a_0^2a_j + P_j(a_0,\ldots,a_{j-1})$$ where $P_j$ is a polynomial and hence $a_j(t)$ is a continuous function existing for $t\in \left[-c\left. , \infty \right) \right.$.     
This concludes the construction of the formal solution $\sum_{k=0} ^{\infty} a_k ^{\pm}(t) x^{\beta_k}$.  \\
\indent  Now let $f(x,t)$ denote any smooth function which is asymptotic to the formal solution $\sum_{k=0} ^{\infty} a_k ^{\pm}(t) x^{\beta_k}$ (by proposition 3.5 in [7] there exists such a function).  By plugging in $f(x,t)$ to $(\ref{dagger})$ we will now show that one obtains a function $g(x,t)\in S^{-\infty}(\mathbb{R}\times I)$ where $I=\mathbb{R}$ if $\beta<\frac{1}{2}$ and $I=\left[-c\left. , \infty \right) \right.$ if $\beta=\frac{1}{2}$.  Let $S_N(f) = \sum_{k=0} ^N a_k ^{\pm}(t)(\pm x)^{\beta_k}$.  Suppose $\left|x\right|\geq 1$ and $J\subset I$ is a compact subset and  $i,j,N\geq 0$ are integers.  Then we have

\begin{eqnarray*}
\partial^i_t\partial^j_x \left( f_t+f^2 f_x+ f_{xxx}\right) & = & \partial^i_t\partial^j_x \left[S_N(f) +\left(f-S_N(f)\right) \right]_t + \partial^i_t\partial^j_x \left[\left( S_N(f) +\left(f-S_N(f)\right) \right)^2\left(S_N(f) +\left(f-S_N(f)\right) \right)_x\right] + \\
                     &    &  + \partial^i_t\partial^j_x\left[S_N(f) +\left(f-S_N(f)\right) \right]_{xxx} 
\end{eqnarray*}

\noindent After expanding the right side we will obtain the expression $\partial^i_t\partial^j_x \left[ S_N(f)_t+S_N(f)^2S_N(f)_x + S_N(f)_{xxx}\right]$ and some terms of the form $\partial^m_t \partial^n_x \left(f-S_N(f) \right)^p \partial^k_t \partial^l_x S_N(f)^q$.  Since the coefficients $a_j(t)$ satisfy the equations above and since $f$ is asymptotic to $\sum_{k=0} ^{\infty} a_k ^{\pm}(t) x^{\beta_k}$ we have for any $k\in \mathbb{N}$ there exists a sufficiently large $N$ such that these two terms are bounded by $C_{k,J,i,j}\left|x \right|^{-k}$.$\square$ \\

\noindent \textbf{ \large{Appendix B:  Additional Proofs}} \\

\textbf{Proof of Proposition 2.1}
\begin{enumerate}
\item 
From the assumptions on $\rho$ the Sobolev inequality in lemma 2.3 implies that we also have $D_+\rho, D_+^2 \rho \in L^2_h$.  Therefore,
\begin{eqnarray*}
\DLIP{D_+^2D_- \rho}{\rho} = \frac{1}{2} \DLIP{\rho}{ D_+^2 D_- \rho - D_+D_-^2 \rho} = \frac{1}{2}\DLIP{\rho}{D_+D_- (hD_+D_- \rho)} = \frac{h}{2}\DLIP{D_+D_-\rho}{D_+D_-\rho} \geq 0
\end{eqnarray*}
\item  From the product formula for $D_+$ it follows by induction on $n$ that $$D_+^n (\rho \cdot \xi) = \sum_{i=0}^n c_i (E^iD_+^{n-i} \rho)( D_+^i \xi)$$
for some constants $c_i\in\mathbb{N}$.  Therefore we have for some constants $c_{i,j,k}$ that,
\begin{eqnarray*}
D_+^n \left(\rho \cdot \nu \cdot \xi\right) = \sum_{i=0}^n c_i \left(E^i D_+^{n-i} \rho \right) \sum_{j=0}^i \tilde{c}_j \left( E^j D_+^{i-j}\nu \right) D_+^j \xi =  \sum_{i_1+i_2 + i_3 = n}c_{i_1,i_2,i_3} \left(E^{i_2 + i_3}D_+^{i_1}\rho \right)\cdot \left(E^{i_3}D_+^{i_2}\nu \right)\cdot\left(D_+^{i_3}\xi\right)
\end{eqnarray*}
  
\item 
\begin{eqnarray*} 
\DLIP{\rho \nu}{D_0 \nu} = \lim_{N \to \infty}\frac{1}{2} \left[ \sum_{-N}^N \rho_n \nu_n \nu_{n+1} -\sum_{-N}^N \rho_n \nu_n \nu_{n-1} \right] = \lim_{N \to \infty}\frac{1}{2} \sum_{-N}^N \rho_n \nu_n \nu_{n+1} - \rho_{n+1} \nu_{n+1} \nu_n = -\frac{1}{2}\DLIP{\nu}{E\nu D_+ \rho}
\end{eqnarray*}
\item 
From the Sobolev inequality in lemma 2.3 the assumptions on $\rho$ and $\nu$ imply that we also have $D_+^i \rho,D_+^i \nu \in L^2_h$ for $i=0,1,2$ and the same thing holds for any such lower derivatives involving a mix of $D_+$, $D_-$, and $D_0$.  Therefore,  
\begin{eqnarray*}
\DLIP{\rho \nu}{D_+^2D_- \nu} & = & \frac{1}{h^2}\sum_{-\infty}^{\infty} \rho_n  \left(\nu_n \nu_{n+2} - 3\nu_n \nu_{n+1} + 3\nu_n^2 - \nu_n \nu_{n-1} \right)  \\
                      & = & \lim_{N \to \infty} \left[\frac{1}{2h^2} \sum_{-N}^N \rho_n (\nu_{n+2}^2 - 3\nu_{n+1}^2 + 3\nu_n^2 - \nu_{n-1}^2) + \right. \\
                      &   & + \left. \frac{1}{2h^2} \sum_{-N}^N \rho_n \left(-(\nu_{n+2} - \nu_n)^2+ 3(\nu_{n+1}-\nu_n)^2 + (\nu_{n-1} - \nu_n)^2 \right) \right] \\
                      & = & \lim_{N \to \infty}\Big[ \frac{1}{2h^2}\sum_{-N}^N(\rho_{n-2}-3\rho_{n-1}+3\rho_n-\rho_{n+1})\nu_n^2 + h.o.t + \frac{1}{2h^2}\sum_{-N}^N \rho_n \left(-(\nu_{n+2}-\nu_{n+1})^2 - \right.  \\
                      &   &  \left. - (\nu_{n+1}-\nu_n)^2  + 2(\nu_{n+2}-\nu_{n+1})(\nu_n - \nu_{n+1}) + 3(\nu_{n+1}-\nu_n)^2 + (\nu_n - \nu_{n-1})^2 \right) \Big] \\
                      & = &-\frac{1}{2}\DLIP{\nu}{\nu D_-^2D_+ \rho} + \lim_{N \to \infty} \left[ \frac{1}{h^2}\sum_{-N}^N \rho_n (\nu_{n+2} - \nu_{n+1})(\nu_n - \nu_{n+1}) \right. + \\
                      &   & + \left. \frac{1}{2h^2}\sum_{-N}^N \rho_n \left(-(\nu_{n+2}-\nu_{n+1})^2 + 2(\nu_{n+1}-\nu_n)^2 + (\nu_n - \nu_{n-1})^2 \right) \right] \\
                      & \geq & -\frac{1}{2}\DLIP{\nu}{\nu D_-^2D_+ \rho} + \frac{1}{2h^2} \sum_{-\infty}^{\infty} \rho_n \left[-2(\nu_{n+2}-\nu_{n+1})^2 + (\nu_{n+1}-\nu_n)^2 + (\nu_n - \nu_{n-1})^2 \right] \\
                      & = &  -\frac{1}{2}\DLIP{\nu}{\nu D_-^2D_+ \rho} + \frac{1}{2} \sum_{-\infty}^{\infty} \left(D_+ \nu_n\right)^2 \left( \frac{\rho_n - \rho_{n-1}}{h} + \frac{\rho_{n+1} - \rho_{n-1}}{h} \right) h \\
                      & = & \frac{1}{2}\DLIP{\nu}{\nu D_-^2D_+ \rho} + \frac{1}{2}\DLIP{D_-\rho D_+\nu}{D_+\nu} + \DLIP{D_0 \rho D_+\nu}{D_+\nu}
\end{eqnarray*}

\end{enumerate}

\hspace{2mm}

\noindent \textbf{Proof of Proposition 2.2} \hspace{1mm} Since
\begin{equation*}
\int_0^1 \frac{d}{dr} \big[g(x+rh,t)\big]dr = h \int_0^1 \frac{dg}{dx}(x+rh,t) dr
\end{equation*}
\noindent we may apply the fundamental theorem of calculus to the left side, divide both sides by $h$, and use the definition of $D_+$ to obtain
\begin{equation*}
D_+ g(x,t) = \int_0^1 \frac{dg}{dx}(x+rh,t)dr
\end{equation*}
\noindent Since $g\in S^{-\infty}(\mathbb{R}\times \left[-c\left. ,\infty \right) \right.)$ it follows that $g$ is Schwartz 'uniformly' on the interval $\left[ 0, T \right]$ so that we may define constants $$A_{l,m} := \sup_{t\in \left[0,T \right]} \SupNorm{x^l \frac{d^m}{dx^m}g(x,t)} < \infty$$   
\noindent  For $r,h \in \left[0,1\right]$ we then obtain that
\begin{eqnarray*}
\SupNorm{x^l \frac{d^m g}{dx^m}(x+rh,t)} & =    & \SupNorm{(x-rh)^l \frac{d^m}{dx^m}g(x,t)} \leq \sum_{q=0}^l \binom{l}{q} \SupNorm{x^q \frac{d^m}{dx^m}g(x,t)} (1 \cdot 1)^{l-q} = \sum_{q=0}^l \binom{l}{q} \cdot A_{q,m} 
\end{eqnarray*}
and hence $g(x+rh,t)$ is Schwartz 'uniformly' for $t \in \left[ 0,T \right]$, and $r,h\in\left[0,1\right]$ - and clearly $\frac{dg}{dx}(x+rh,t)$ is too.  Moreover the map $x\mapsto \int_0^1 \frac{dg}{dx}(x+rh,t)dr$ remains Schwartz uniformly for $t\in \left[0,T\right]$ and $h\in\left[0,1\right]$ because
\begin{eqnarray*}
\left| x^l \frac{d^m}{dx^m} \int_0^1 \frac{dg}{dx}(x+rh,t)dr \right| \leq \int_0^1 \SupNorm{x^l \frac{d^{m+1}g}{dx^{m+1}}(x+rh,t)} dr \leq  C_{l,m+1} 
\end{eqnarray*}
and we may take supremum over $x$ on the left side to obtain
$\SupNorm{ x^l \frac{d^m}{dx^m} D_+ g(x,t)}  \leq C_{l,m+1}$
and hence $D_+ g(x,t)$ is Schwartz uniformly for $t\in\left[0,T\right]$ and $h\in\left[0,1\right]$. \\
\indent If we continue inductively on $n$ then we see that since $D_+^{n-1}g(x,t)$ is Schwartz uniformly for $t\in\left[0,T\right]$ and $h\in\left[0,1\right]$ we have also that $D_+^n g(x,t)$ is Schwartz uniformly for $t\in\left[0,T\right]$ and $h\in\left[0,1\right]$ and moreover $(1+x^2)^N (D_+^n g(x,t))^2$ is too. Thus there exists $C_{N,n}>0$ independent of $t\in\left[0,T\right]$ and $h\in\left[0,1\right]$ so that
\begin{equation*}
(1+x^2)^N\big(D_+^n g(x,t) \big)^2 \leq \frac{C_{N,n}}{1+x^2}
\end{equation*}
and from this it follows that for any $t\in \left[0,T\right]$ (in particular $t_j$) we have,
\begin{eqnarray*}
\DLNS{\IHN{x}^N D_+^n g(x,t)} &  \leq & \sum_{p=-\infty}^{\infty} \frac{C_{N,n}h}{1+x_p^2} \leq C_{N,n}\cdot 1 + 2\sum_{p\geq 1} \int_{p-1}^{p} \frac{C_{N,n}h}{1+x_p^2} ds \\
                              & \leq & C_{N,n} + 2\sum_{p\geq 1} \int_{h(p-1)}^{h p} \frac{C_{N,n}}{1+y^2} dy \leq C_{N,n} + 2\int_0^{\infty} \frac{C_{N,n}}{1+y^2} dy \leq C_{N,n}
\end{eqnarray*} $\square$ \\

\noindent \textbf{Proof of Lemma 2.3} \hspace{1mm} If $u \notin L^2_h$ then the right side of both inequalities is infinity so they clearly hold in that case, thus we may assume that $u \in L^2_h$.  First we will prove inequality 1.  By lemma 4.1, it suffices to prove that  
\begin{equation}
\CLN{\frac{d^k}{dx^k} U } \leq C\left( \CLN{ U} + \CLN{\frac{d^n}{dx^n}U} \right)  \tag{B.1}
\end{equation}  
for $0 \leq k \leq n$ and for some $C>0$ where $U = I_h u$.  Since $\frac{\left\langle \xi \right\rangle^n}{1+\left|\xi\right|^n}$ is bounded by some constant we have that 
\begin{eqnarray*}
\CLN{ \frac{d^k}{dx^k} U }  = \CLN{ \xi^k (FU)(\xi) } \leq \CLN{ \IHN{ \xi }^k (FU)(\xi) } \leq  C\left(\CLN{ FU } + \CLN{ \xi^n (FU)(\xi) }\right) = C\left( \CLN{ U } + \CLN{ \frac{d^n}{dx^n}U } \right)
\end{eqnarray*}
Therefore inequality 1 holds.  \\
\indent For property 2, we first note that for any $x_n$ in the mesh there exists some intermediate point $\tilde{x}_n $(not necessarily in the mesh) such that $x_n \leq \tilde{x}_n \leq x_{n+k}$ and
\begin{eqnarray*}
\left| D_+^k u(x_n) \right| &  = & \left| D_+^k U(x_n) \right| = \left| \frac{d^k}{dx^k} U(\tilde{x}_n) \right| \leq \SupNorm{ \frac{d^k}{dx^k} U }
\end{eqnarray*}

\noindent Therefore $\SupNorm{ D_+^k u } \leq \SupNorm{ \frac{d^k}{dx^k}U }$.  Moreover, for any $x\in \mathbb{R}$ we have, by the Cauchy-Schwarz inequality, the below estimates  \\

\begin{eqnarray*}
\left| \frac{d^k}{dx^k} U (x) \right| & \leq & \frac{1}{2\pi} \left( \int_{\mathbb{R}} \left( \left\langle  y\right\rangle^{-1} \right) \left(\left\langle  y\right\rangle^{k+1} \left|F_h u(y) \right| \right)  dy \right) \leq \frac{1}{2\pi} \left( \int_{\mathbb{R}}\frac{dy}{1+\left|y \right|^2} \right)^{\frac{1}{2}} \cdot \left( \int_{\mathbb{R}}\left(1 + \left|y\right|^2 \right)^{k+1} \left|(F_h u)(y) \right|^2 dy \right)^{\frac{1}{2}} \\ 
                                        &  \leq  & C \left( \CLN{ (F_h u)(y) } + \CLN{ y^{k+1}(F_h u)(y) }\right) = C \left( \CLN{ U } +  \CLN{\frac{d^{k+1}}{dx^{k+1}}U }   \right)
\end{eqnarray*}

\noindent and hence 
\begin{equation}
\SupNorm{ \frac{d^k}{dx^k}U } \leq C \left( \CLN{ U } +  \CLN{ \frac{d^{k+1}}{dx^{k+1}}U }   \right) \tag{B.2}
\end{equation} 
  Therefore, by lemma 4.1 we have 
\begin{eqnarray*}
\SupNorm{ D_+^k u } & \leq & \SupNorm{ \frac{d^k}{dx^k}U } \leq C \left( \CLN{ U } +  \CLN{\frac{d^{k+1}}{dx^{k+1}}U }   \right) \leq C \left( \DLN{ u} + \DLN{ D^{k+1}_+u} \right)
\end{eqnarray*}

\noindent where $C$ is clearly independent of $h$.  Since the above inequality holds for any $k$, we may apply the above inequality to $D_+^{k+j} u$ for each $j\in \mathbb{N}$ satisfying $k+j < n$ so that we will obtain inequality 2.  $\square$ \\

\noindent \textbf{Proof of Corollary 2.4} \hspace{1mm} We will first prove inequality 1.  For $N=0$ the result follows from lemma 2.3.  We proceed by induction on $N$.  Assume the result is true for each $M \leq N-1$.  We will prove it for $M=N$.  For $j=0$ the result is trivial with $C=1$.  Consider the case $j=1$.  Then we have, by the product rule for $D_+$,
\begin{equation*}
D_+\big(\IHN{x_n}^N u(x_n)\big) = D_+\big(\IHN{x_n}^N\big)\cdot u(x_{n+1})+\IHN{x_n}^N D_+ u(x_n)
\end{equation*}
therefore, 
\begin{eqnarray*}
\DLN{\IHN{x}^N D_+ u} & \leq & \DLN{D_+ \big(\IHN{x_n}^N\big) \cdot u(x_{n+1})} + \DLN{D_+\big( \IHN{x}^N u(x_n) \big)} \leq C_{N,k} \Big( \DLN{\IHN{x}^N u} + \DLN{D_+^{k+1}\big(\IHN{x}^N u\big)} \Big) \\
											& \leq & C_{N,k} \Big( \DLN{\IHN{x}^N u} +\sum_{0\leq l \leq N-1} \DLN{\IHN{x}^{N-l} D_+ ^{k+1-l}u} + \sum_{N\leq l \leq k+1}\DLN{\IHN{x}^{N-l}D_+^{k+1-l} u } \Big) \\
											& \leq & C_{N,k}\left[ \DLN{\IHN{x}^N u} + \DLN{\IHN{x}^{N} D_+^{k+1}u} + \sum_{1\leq l \leq N-1}\Big( \DLN{\IHN{x}^{N-l} u}+\DLN{\IHN{x}^{N-l} D_+^{k+1}u} \Big) \right. \\
											& & \left. + \sum_{N\leq l \leq k+1}\DLN{1\cdot D_+^{k+1-l} u } \right] \leq C_{N,k} \Big( \DLN{\IHN{x}^N u} + \DLN{\IHN{x}^{N} D_+^{k+1}u} \Big)
\end{eqnarray*}
Hence we have proven the case $M=N$ for $j=0,1$.  Assume it is true for $p=j-1$.  We will prove it for $p=j$.  By applying the above inequality to $D^{j-1}u$ we obtain,
\begin{eqnarray*}
\DLN{\IHN{x}^N D_+^j u} & \leq & C_{N,k} \Big( \DLN{\IHN{x}^N D_+^{j-1}u} + \DLN{\IHN{x}^N D_+^{j+k}} \Big) \leq C_{N,k} \Big( \DLN{\IHN{x}^N u} + \DLN{\IHN{x}^N D_+^{j+k}} \Big)
\end{eqnarray*}
which concludes both inductions.  \\
\indent Now we prove inequality 2.  By inequality 1 and lemma 2.3 we have for $k\geq 1$,
\begin{eqnarray*}
\SupNorm{\IHN{x}^N D_+^j u} & \leq & C \Big( \DLN{\IHN{x}^N D_+^j u} + \DLN{D_+\big( \IHN{x}^N D_+^j u \big)} \Big) \\
                            & \leq & C\Big( \DLN{\IHN{x}^N D_+^j u} + \DLN{\IHN{x}^{N-1} D_+^j u} + \DLN{\IHN{x}^N D_+^{j+1} u } \Big) \leq C \Big( \DLN{\IHN{x}^N u} + \DLN{\IHN{x}^N D_+^{j+k}} \Big)
\end{eqnarray*} $\square$ \\

\noindent \textbf{Remark} \hspace{1mm} By using (B.1) and (B.2) we see that the proof of corollary 2.4 can also be applied to $U = I_h u$ with $D_+$ replaced by $\partial_x$ and we would show that for every $N,k,j \in \mathbb{N}$ there exists $C_{N,k,j}>0$ such that for any $h \in (0,1)$ and for all mesh functions $u(x_n)$ defined on $\mathbb{R}_h$ we have the inequalities,
\begin{enumerate}
\item   $ \CLN{ \IHN{x}^N \partial_x^j U } \leq C_{N,j,k} \left( \CLN{ \IHN{x}^N U} + \CLN{\IHN{x}^N \partial_x^{j+k} U} \right) $  \hspace{1cm}  for $k\geq 0$
\item  $\SupNorm{ \IHN{ x }^N \partial_x^j U } \leq C_{N,j,k} \left( \CLN{ \IHN{ x }^N U} + \CLN{\IHN{ x }^N \partial_x^{j+k} U} \right) $   \hspace{1cm} for $k\geq 1$
\end{enumerate}

\hspace{2mm}
 
\noindent \textbf{Proof of Lemma 5.1} \hspace{1mm} By assumption we have,
\begin{eqnarray*}
\frac{d}{dt}\left(\eta(t) e^{-c_1 t} \right) = \frac{d \eta}{dt} \cdot e^{-c_1 t} - c_1 \eta(t) e^{-c_1 t} \leq c_1 \eta(t) e^{-c_1 t} + c_2 e^{-c_1 t}- c_1 \eta(t) e^{-c_1 t}= c_2 e^{-c_1 t}
\end{eqnarray*}
\noindent Integrating the above inequality from zero to $t$ and applying the fundamental theorem of calculus yields the desired inequality.
$\square$ 

\noindent \textbf{\Large{References}} \\

\noindent [1] I. N. Bondareva,\emph{ The Korteweg-De Vries Equation in classes of increasing functions with prescribed asymptotics as $|x| \rightarrow \infty$}, Math. USSR-Sb. $\mathbf{50}$ (1) (1985) 125-135.\\

\noindent [2] I. N. Bondareva, M. A. Shubin, \emph{Increasing Asymptotic Solutions of the Korteweg-De Vries Equation and its Higher Analogues}, Soviet Math. Dokl. $\mathbf{26}$ (3) (1982) 716-719.\\ 

\noindent [3]  T. Kappeler, P. Perry, M. Shubin, P. Topalov,\emph{ Solutions of mKdV in classes of functions unbounded at infinity} , J. Geom. Anal. $\mathbf{18}$ (2) (2008) 443-477.\\

\noindent [4]  C. Kenig, G. Ponce, L. Vega, \emph{Global solutions for the KdV equation with unbounded data}, J. Differential Equations $\mathbf{139}$ (2) (1997) 339-364. \\

\noindent [5] A. Menikoff, \emph{ The Existence of Unbounded Solutions of the Korteweg-De Vries Equation}, Comm. Pure Appl. Math. $\mathbf{25}$ (1972) 407-432.\\

\noindent [6] W. Rudin, \emph{Real and Complex Analysis}, third ed.,  McGraw-Hill, New York, 1987. \\

\noindent [7] M. Shubin, \emph{Pseudodifferential Operators and Spectral Theory}, second ed., Springer-Verlag, Berlin, 2001, translated from the 1978 Russian original by S. I. Andersson.\\
 
\noindent [8] F. Stummel,\emph{ Elliptische Differenzenoperatoren unter Dirichletrandbedingungen}, Math. Z. $\mathbf{97}$ (1967) 169-211.

\end{document}